\documentclass[12pt]{article}
\usepackage{hyperref}
\usepackage{amssymb}
\usepackage{latexsym}
\usepackage{amsmath}
\usepackage{amscd}
\usepackage{theorem}
\usepackage{xy}
\xyoption{all}
\CompileMatrices

\usepackage{graphicx}

\makeatletter

\def\@@seccntfont{\bfseries\slshape}
\def\@@secheaderfont{\bfseries\upshape}
\def\@@precnt{{\upshape}}
\def\@@postcnt{{\upshape}}
\def\@@startsection#1#2#3#4#5#6{\if@noskipsec \leavevmode \fi\par
  \@tempskipa #4\relax\@afterindenttrue
  \ifdim \@tempskipa <\z@\@tempskipa -\@tempskipa
  \@afterindentfalse\fi
  \if@nobreak\everypar{}\else\addpenalty\@secpenalty\addvspace\@tempskipa\fi
  \@ifstar{\@dblarg{\@sect{#1}{#2}{#3}{#4}{#5}{#6}}}%
    {\@dblarg{\@@sect{#1}{#2}{#3}{#4}{#5}{#6}}}}
\def\@@sect#1#2#3#4#5#6[#7]#8{\ifnum
  #2>\c@secnumdepth\let\@svsec\@empty
  \else\refstepcounter{#1}\protected@edef\@svsec{\@seccntformat{#1}\relax}\fi\@tempskipa #5\relax
  \ifdim
  \@tempskipa>\z@\begingroup#6{\@hangfrom{\hskip#3\relax\@svsec}\interlinepenalty 
  \@M #8\@@par}\endgroup\csname #1mark\endcsname{#7}\else
  \def\@svsechd{#6{\hskip #3\relax\@svsec #8}\csname
    #1mark\endcsname{#7}}\fi\@xsect{#5}}
\def\setseccntfmt{\renewcommand{\@seccntformat}[1]{\S
  \csname the##1\endcsname\hspace{1ex}}}
\def\@@seccntfmt{\renewcommand{\@seccntformat}[1]{%
  {\@@seccntfont\@@precnt\csname
  the##1\endcsname}\@@postcnt\hspace{1ex}}}
\newcommand{\@@secnostar}[1][]{\def\tmpa{}\def\tmpb{#1}%
  \ifx\tmpa\tmpb\def\tmpa{-1ex}\else\def\tmpa{-\labelsep}\fi\@@seccntfmt\@@startsection
  {\@@name}{\@@level}{0mm}{-\baselineskip}{\tmpa}{\@@secheaderfont}{#1}}
\newcommand{\@@secstar}[2][]{\def\tmpa{}\def\tmpb{#1}%
  \ifx\tmpa\tmpb\def\tmpa{#2}\else\def\tmpa{#1}\fi\@@seccntfmt\@startsection
  {\@@name}{\@@level}{0mm}{-\baselineskip}{-\labelsep}{\@@secheaderfont}[\tmpa]{#2}}
\def\@@section{\def\@@name{section}\def\@@level{1}\@ifstar{\@@secstar}{\@@secnostar}}
\def\@@subsection{\def\@@name{subsection}\def\@@level{2}\@ifstar{\@@secstar}{\@@secnostar}}
\def\@@subsubsection{\def\@@name{subsubsection}\def\@@level{3}\@ifstar{\@@secstar}{\@@secnostar}}
\def\@@paragraph{\def\@@name{paragraph}\def\@@level{4}\@ifstar{\@@secstar}{\@@secnostar}}
\def\@@subparagraph{\def\@@name{subparagraph}\def\@@level{5}\@ifstar{\@@secstar}{\@@secnostar}}
\let\@@latexsection=\section\let\mysection=\@@section
\let\@@latexsubsection=\subsection\let\mysubsection=\@@subsection
\let\@@latexsubsubsection=\subsubsection\let\mysubsubsection=\@@subsubsection
\let\@@latexparagraph=\paragraph\let\myparagraph=\@@paragraph
\let\@@latexsubparagraph=\subparagraph\let\mysubparagraph=\@@subparagraph
\newif\if@@lastcharstar
\def\@@xxlastcharstar#1{\gdef\@@prevchar{}\@@lastcharstarfalse\@@yylastcharstar#1\end}
\def\@@yylastcharstar#1{\ifx#1\end\def\tmpa{*}\ifx\@@prevchar\tmpa\@@lastcharstartrue\fi
  \let\@@next=\relax\else\def\@@prevchar{#1}\let\@@next=\@@yylastcharstar\fi\@@next}
\gdef\@thm#1#2{\@@xxlastcharstar{#1}\if@@lastcharstar\else\refstepcounter{#1}\fi
  \trivlist\@topsep \theorempreskipamount\@topsepadd
  \theorempostskipamount
  \@ifnextchar [{\@ythm{#1}{#2}}{\@begintheorem{#2}{\csname
  the#1\endcsname}\ignorespaces}}
\gdef\th@nonumplain{\normalfont\itshape
  \def\@begintheorem##1##2{\item[\hskip\labelsep\theorem@headerfont
  ##1]}%
  \def\@opargbegintheorem##1##2##3{%
    \item[\hskip\labelsep \theorem@headerfont ##1\ (##3)]}}
\gdef\th@change{
  \def\@begintheorem##1##2{\item[\hskip\labelsep{\@@seccntfont
  \@@precnt##2\@@postcnt\hskip 1ex}\theorem@headerfont ##1]}%
  \def\@opargbegintheorem##1##2##3{\item[\hskip\labelsep{\@@seccntfont
  \@@precnt##2\@@postcnt\hskip 1ex}\theorem@headerfont ##1\ (##3)]}}
\makeatother   
\def\newtheoremset[#1][#2]{{\theoremstyle{change}\newtheorem{-#1}[section]{#2}%
  \newtheorem{#1}[subsection]{#2}\newtheorem{*#1}[subsubsection]{#2}%
  \newtheorem{**#1}[paragraph]{#2}\newtheorem{***#1}[subparagraph]{#2}}
  {\theoremstyle{nonumplain}\newtheorem{#1*}{#2}}}
{\theoremheaderfont{\normalfont\bfseries}
{\theorembodyfont{\normalfont\itshape}
\newtheoremset[theorem][Theorem.]\newtheoremset[proposition][Proposition.]
\newtheoremset[lemma][Lemma.]\newtheoremset[corollary][Corollary.]\newtheoremset[properties][Properties.]\newtheoremset[property][Property.]
\newtheoremset[conjecture][Conjecture.]}\newtheoremset[sublemma][Sublemma.]\newtheoremset[assumption][Assumption.]
{\theorembodyfont{\normalfont\rmfamily}
\newtheoremset[definition][Definition.]
\newtheoremset[example][Example.]
\newtheoremset[facts][Facts.]\newtheoremset[fact][Fact.]
\newtheoremset[examples][Examples.]\newtheoremset[claim][Claim.]
\newtheoremset[notation][\itshape\mdseries Notation.]
\newtheoremset[remark][\itshape\mdseries Remark.]
\newtheoremset[remarks][\itshape\mdseries Remarks.]}}
\def\XYmatrix{\xymatrix@M=5pt} 
\def\ncmd{\newcommand}
\ncmd{\xysubset}[1][r]{\ar@<-2.5pt>@{^(-}[#1]\ar@<2.5pt>@{_(-}[#1]}
\ncmd{\XYmatrixc}[1]{\vcenter{\XYmatrix{#1}}}
\ncmd{\xyto}[1][r]{\ar@{->}[#1]}      \ncmd{\xyinj}[1][r]{\ar@{^(->}[#1]}
\ncmd{\xysurj}[1][r]{\ar@{->>}[#1]}   \ncmd{\xyline}[1][r]{\ar@{-}[#1]}
\ncmd{\xydotsto}[1][r]{\ar@{.>}[#1]}  \ncmd{\xydots}[1][r]{\ar@{.}[#1]}
\ncmd{\xyleadsto}[1][r]{\ar@{~>}[#1]} \ncmd{\xyeq}[1][r]{\ar@{=}[#1]}
\ncmd{\xyequal}[1][r]{\ar@{=}[#1]}    \ncmd{\xyequals}[1][r]{\ar@{=}[#1]}
\ncmd{\xymapsto}[1][r]{\ar@{|->}[#1]}\ncmd{\xyimplies}[1][r]{\ar@{=>}[#1]}
\ncmd{\xytofrom}[1][r]{\ar@{<->}[#1]} 
   
\def\XYTOTO[#1]^#2_#3{\xyto[#1]<0.5ex>^{#2}\xyto[#1]<-0.5ex>_{#3}}
\ncmd{\xytoto}[1][r]{\XYTOTO[#1]}

\newenvironment{proof}
{\noindent {\em Proof.}}
{\hfill $\Box$\medbreak}

\newcommand{\eqdef}{\;{:=}\;}

\newcommand{\C}{{\mathbb C}}
\newcommand{\Q}{{\mathbb Q}}
\newcommand{\R}{{\mathbb R}}
\newcommand{\Z}{{\mathbb Z}}
\newcommand{\op}{\operatorname}

\newcommand{\sign}{\op{sign}}

\newcommand{\End}{\op{End}}

\newcommand{\Ker}{\op{Ker}}
\newcommand{\cok}{\op{coker}}

\newcommand{\ind}{\op{ind}}
\newcommand{\Nov}{\op{Nov}}
\newcommand{\tensor}{\otimes}

\newcommand{\qed} {\hfill$\Box$}
\newcommand{\cm}{{\mathcal M}}

\newcommand{\Tr}{\op{Tr}}
\newcommand{\p}{\partial}
\newcommand{\ls}{{\mathcal L}}
\newcommand{\lc}{{\mathcal C}}

\begin{document}
\title{Reidemeister Torsion in Floer-Novikov Theory\\ and \\Counting
  Pseudo-holomorphic Tori, I.}
\author{\scshape Yi-Jen Lee \\  \small Department of Mathematics, Purdue
  University\\ \small West Lafayette, IN 47907 U.S.A. \\ \small {\tt yjlee@math.purdue.edu}}
\date{{\small 
This version: April 2005 
}}

\begin{titlepage}
\maketitle

\abstract{This is the first part of an article in two parts,
  which builds the foundation of a Floer-theoretic invariant, \(I_F\).

The Floer homology can be trivial in many variants of the
  Floer theory; it is therefore interesting to consider more refined
  invariants of the Floer complex. We consider one such
  instance---the Reidemeister torsion \(\tau_F\) of the Floer-Novikov complex
  of  (possibly {\em non-hamiltonian}) symplectomorphisms. \(\tau_F\) turns
  out {\em not} to be invariant under hamiltonian isotopies, but this
  failure may be fixed by introducing certain ``correction term'':
We define a Floer-theoretic zeta function
  \(\zeta_F\), by counting perturbed pseudo-holomorphic tori
  in a way very similar to the genus 1 Gromov invariant. The main
  result of this article states that under suitable monotonicity
  conditions, the product \(I_F:=\tau_F\zeta_F\) is
  invariant under hamiltonian isotopies. In fact, \(I_F\) is 
invariant under general symplectic isotopies 
when the underlying symplectic manifold
  \(M\) is monotone.  

Because the torsion invariant we consider is a {\em
    not a homotopy} invariant, the continuation method used in typical
  invariance proofs of Floer theory does not apply; instead, 
  the detailed bifurcation analysis is worked out. This is the first time such
  analysis appears in the Floer theory literature in its entirety.

Applications of \(I_F\), and the construction of \(I_F\) in different
versions of Floer theories are discussed in sequels to this article
\cite{L2, LS, ogw}.}
\newpage

\tableofcontents
\end{titlepage}

\section{Introduction}
This is the first of a series of papers dealing with torsion
invariants in Floer theories. In this paper and its companion, Part II
\cite{part2}, we concentrate on
establishing the foundation and invariance of the torsion invariant
\(I_F\), and thus our main purpose here is to develope a general method for
proving such foundational results. Examples, applications, and
adaptations to various other versions of Floer theories will be
discussed in subsequential papers (e.g \cite{L2, LS, ogw}).

Part I contains the construction
of the proposed invariant, the main framework of the invariance proof,
and the structure theorems for the relevant moduli spaces. The heavy
analysis required for the proof of expected bifurcation behavior is
postponed to Part II.
 
\subsection{Background and Motivations}
The two original versions of symplectic Floer homologies 
have been initially introduced as tools for proving the Arnold
conjecture, which gives a lower bound on the number of fixed points of
a hamiltonian symplectomorphism by the total betti number of the symplectic
manifold \cite{floer.cmp, floer.jdg}. 
Very roughly speaking, the Floer homology is the homology of a ``Morse function''
(the action functional) on an infinite dimensional space. In the first
version (the ``absolute'' version), this infinite dimensional space is
the free loop space of the compact symplectic manifold
\((M, \omega)\). In the second version (the ``relative'' version), 
it is the space of paths ending at two transversely-intersecting
lagrangian submanifolds in \(M\). 

To make the above heuristics work actually requires highly nontrivial
transversality and compactness results for the relevant moduli spaces,
which make use the assumptions of monotonicity or \(\pi_2=0\) in Floer's
original papers. These assumptions have been subsequently weakened
through the efforts of many people; recently it has been completely removed
in the absolute version via the virtual moduli technique by several groups of
people, see for example \cite{liu-tian,fukaya-ono, ruan, salamon.park}. For the relative
version, traditionally the success has been more limited,
though the recent paper \cite{FOOO} is able to deal with fairly
general settings, whose implication includes the general version of Arnold
conjecture proved in \cite{liu-tian, fukaya-ono}.

The basic strategy in these works is the same as Floer's, namely
proving the invariance of the Floer homology under hamiltonian isotopies. This
enables one to compute the Floer homology at a small,
\(t\)-independent hamiltonian, where \(HF_*(M)=H_*(M)\), the usual homology of the
symplectic manifold.

Floer's proof guarantees the existence of fixed points for any
hamiltonian symplectomorphism of \(M\), because for closed \(M\),
\(H_*(M)\) is never trivial. This is however no longer true in many
variants of Floer theory. Here are some examples.

\subsubsection{Examples of vanishing Floer homology.}

{\bf Example A. (Symplectic manifolds with boundary)}
There are various ways of defining Floer homologies for symplectic
manifolds with contact-type boundaries \cite{Vit, CFHW}. The Floer homology
may be trivial in this situation. E.g. $HF_*(D^2)=0$ according to Viterbo's definition.
\medbreak

\noindent
{\bf Example B. (Space of loops in a nontrivial homotopy class)}
To find {\it noncontractible} hamiltonian orbits, one might use the
space of noncontractible loops (in a fixed homotopy class) to define
Floer homology. However, this version of Floer homology vanishes by
its invariance under hamiltonian isotopies, because for small
hamiltonians, there is obviously no closed orbit \cite{BH}. 
\medbreak

\noindent
{\bf Example C. (Floer theory of lagrangian intersections)}
There are many examples where the relative version---the
lagrangian intersection Floer homology vanishes. For example, the
Floer homology of compact lagrangian submanifolds $L\subset\C^n$ vanishes
once it is defined and invariant under hamiltonian isotopies, because it
is easy to find hamiltonian isotopies disengaging one from the
other. According to \cite{FOOO}, such \(L\) can be any
simply-connected lagrangian submanifold in \(\C^n\) for \(n>2\).
\medbreak

\noindent
{\bf Example D. (Twisted versions of Floer theories)}
This is the our prime example. The Floer theory considered in this
article will always be the twisted version.
 
The action functional is in general not globally
defined (the usual construction only defines a
closed 1-form on the loop space, which is the differential of a
real-valued action function only in special cases). Therefore
one often needs to consider twisted versions of Floer
homology, modeling on the Morse theory of closed
1-forms introduced by S.P. Novikov. A typical example is the Floer homology of a {\it non-hamiltonian}
symplectomorphism $f$. 

For any symplectic isotopy connecting the identity with a
non-hamiltonian symplectomorphism \(f\), Le-Ono \cite{le-ono}
defined such a twisted version of Floer homology and showed following
Floer's strategy that when \(M\)
is monotone, its total betti number equals the total betti number of the Novikov homology
\(HN(M, \theta_f)\). (Le-Ono actually had a slightly
weaker assumption on \(M\).) $\theta_f\in H^1(M)$ above is the 
``flux'' or ``Calabi invariant'' of $f$, 
and the Novikov homology \(HN(M, \theta_f)\) is a twisted
version of Morse homology.

Le-Ono's result guarantees the existence of symplectic fixed points in
many cases; however, there are also many examples where this version of
Floer homology vanishes. The reason is that {\it the ``twisting'' procedure often reduces
the rank of the homology}: different twisted versions of homologies
correspond to homologies of different coverings, and the larger the
covering group is, the smaller the rank of the corresponding homology is.
\bigbreak

\mysubsubsection[Torsion invariants in Floer and Morse theories.]
Naturally, one is interested in more refined invariants of the Floer
complex when the Floer homology vanishes. 

It is known that the finite dimensional Morse theory captures much
more than homological informations of the underlying manifold; in
fact, Cohen-Jones-Segal showed that Morse theory recovers the entire
homeomorphism type of the underlying manifold \cite{CJS}. One
therefore expects that Floer theory similarly has much more to offer
than merely the Floer homology. The path from Morse theory to Floer
theory is strewn with heavy technical difficulties (see \S1.2 below
for some discussion); here are however some instances of progress in 
this direction, in the context of symplectic Floer theories:

Fukaya proposed to study 
the \(A_\infty\) structure of Floer theory which will capture the full
rational homotopy type of the Floer complex. In \cite{fukaya-oh}, Fukaya and Oh showed that
the \(A_\infty\) category of the lagrangian intersection Floer theory for
cotangent bundles is equivalent to the \(A_\infty\) category of Morse theory.

As the first {\it non-homotopy} invariant historically, the
Reidemeister torsion is the obvious next candidate to consider.
In 1994, Fukaya proposed studying torsion and even higher torsions  
in symplectic Floer theories. 
 In \cite{fukaya}, He sketched the definition
of the Whitehead torsions for symplectomorphisms and lagrangian
intersections, claimed that they are invariant under hamiltonian isotopies, and conjectured that the Floer 
homology together with the torsion give complete obstructions of hamiltonian 
isotopying a symplectomorphism to one without fixed points, or a lagrangian submanifold
to one without intersection with another fixed lagrangian submanifold
(the ``symplectic s-cobordism conjecture''). There were however no details.

More recently, rigorous works in this direction were done by
Eliashberg-Gromov \cite{Gromov-Eliashberg} and M.G. Sullivan
\cite{Sul}. However, to sidestep the substantial technical difficulties
in the invariance proof (see \S1.2, II.1
for more discussion), these results require
restrictive assumptions or ad hoc methods: Eliashberg-Gromov consider
only graphical Lagrangian submanifolds in cotangent bundles, so that it
may be reduced to purely finite-dimensional method of
generating functions. Sullivan made several restrictive assumptions to ensure 
that the action functional is globally-defined, which do not hold in 
general (cf. discussion in Example D above). 
He also used a
stabilization trick to sidestep the analysis of death-birth bifurcations, which
applies only to the untwisted version of Lagrangian-intersection Floer theory. 
In both \cite{Gromov-Eliashberg} and \cite{Sul}, the Lagrangian
submanifolds have to be non-compact for their versions of
torsion to be nontrivial. (Compare with \cite{L2}, where an adaptation of
this article gives torsion invariants that are often nontrivial for
compact Lagrangian submanifolds).

One goal of this article is to provide, once and for all, a direct and
general invariance proof that should work for any version of Floer
theory, by overcoming the analytical difficulties.   
In this paper, we consider a general version of the Floer theory of symplectomorphisms (cf. \S3.1), which covers Examples B
and D above. Analogous results for the versions of Floer theories in
Examples A and C can be obtained by simple adaptations of this
paper. (See e.g. \cite{L2, LS}). We study the Reidemeister torsion of the Floer
complex, denoted \(\tau_F\). The torsion we use is an abelian (and
hence weaker) version of the Whitehead
torsion. On the other hand, since this definition uses a standard set of 
bases singled out by the Morse-theoretic context, it is also 
in a sense more refined. See \S2.1 for the precise definition.

However, \(\tau_F\) is {\em not} an invariant for the twisted versions
of Floer theories in Example D. Nevertheless, we find that one may 
construct an invariant by taking into account (perturbed)
pseudo-holomorphic tori. More precisely, we define the product
\[
I_F:=\tau_F \zeta_F,
\]
where the ``zeta function'' \(\zeta_F\) is a generating series
counting perturbed pseudo-holomorphic tori. Our main result is that
\(I_F\) is a symplectic invariant. 

In fact, this result should be viewed as an infinite-dimensional
analog of our previous results
on the torsion invariants of Morse 1-forms \cite{hl1,hl2, Hut}. 
Given a Morse 1-form \(\theta\) on a closed, finite dimensional
manifold \(M\), we considered the dynamics of the flow generated by the
vector field dual to \(\theta\). The Reidemeister torsion \(\tau\) of
the associated Morse-Novikov complex counts  flow lines
ending in critical points in a sense; whereas the dynamical zeta
function \(\zeta\) counts the closed orbits. We showed that the
product \(I=\tau\zeta\) is independent of the metric or \(\theta\),
though neither \(\tau\) nor \(\zeta\) alone is invariant. Actually, \(I\) is equivalent
to the combinatorially-defined Reidemeister torsion of the manifold \(M\). 
In keeping with the picture of Floer theory as an infinite dimensional
Morse theory, the flow lines in the loop space correspond to perturbed
pseudo-holomorphic curves. Thus,
\(\tau_F\) and \(\zeta_F\) above are respectively the infinite-dimensional analogs of the 
Reidemeister torsion of the Morse-Novikov complex and the dynamical
zeta function. 

\mysubsubsection[Relation with Gromov Theory.]
Due to the nature of the definition of \(I_F\) as a product, one may
either regard it as a refinement of Floer homologies, or as a counting invariant of
pseudo-holomorphic curves. These two perspectives lead to different
types of applications. One incentive for our choice of the
abelian version of torsion over the non-abelian Whitehead torsion is so that
\(I_F\) counts perturbed pseudo-holomorphic curves with {\it homology
  class}, which is more in keeping with the usual definition of curve-counting
invariants (Gromov invariants). In
fact, \(I_F\) much resembles the Gromov-Taubes
invariant. The perturbed pseudo-holomorphic tori here has the
interpretation as the perturbed pseudo-holomorphic sections of the
symplectic mapping tori \(\Big(M\times[0,1]/(x, 0)\sim(f(x),
1)\Big)\times S^1\), and the definition of \(\zeta_F\) is very similar
to the zeta function introduced by Ionel-Parker in
\cite{IP}, which computes an averaged version of the genus 1 Gromov
invariant of symplectic mapping tori. 
A fascinating problem is to better understand the precise relation
between \(I_F\) and genus 1 Gromov invariants: such a relation would 
provide a link between Floer theory and Gromov theory, which may be
applied in both directions. Some
simplest cases of this relation are discussed 
for two variants of \(I_F\) in \cite{L2}, which lead to interesting results in
symplectic topology that are beyond the reach of either Floer theory
or Gromov theory alone. 

More generally, ``physical reasoning'' leads one to expect a
Floer-theoretic interpretation of genus 1 Gromov invariants along the
line of our construction of \(I_F\). Philosophically, mirror symmetry
is a correspondence between symplectic (A-model) and complex (B-model)
geometries. Since
the complex side is typically easier to compute, mirror symmetry
produces many interesting conjectural formulae for the curve-counting
invariants on the symplectic side. Recently, much progress has been
made in mathematically rigorous formulation and verification of 
0-loop mirror symmetry, but the higher genera case remains little understood.
In their famous paper in 1993,
the four physicists Bershadsky, Cecotti, Ooguri, Vafa extended
mirror symmetry to higher genera curves \cite{BCOV}.  At the 1-loop level,
this says that the generating series of 
genus 1 Gromov-Witten invariants from the A-model
side should correspond to certain holomorphic analytic torsion from the
B-model side. 
On the other hand, Kontsevich formulated the 0-loop mirror symmetry
as an equivalence between the Fukaya-Floer category
from the A model side, and the category of coherent sheaves on the B
model side. In particular, the Floer cohomologies should correspond to
certain sheaf cohomologies on the B model side. A naive generalization
of Kontsevich's proposal leads one to expect higher-loop mirror
symmetry as an equivalence of secondary invariants on symplectic and
complex sides. In particular,
\begin{eqnarray*}
\text{generating series of genus 1 Gromov invariants} &\text{``=''} &\text{
Floer-theoretic torsion \(\zeta_F\)}\\ &\text{``=''}&\text{
torsion on symplectic side}; \\
\text{holomorphic
analytic torsion} &\text{``=''}&\text{torsion on the complex side}.
\end{eqnarray*} 
Our construction of using \(\zeta_F\) to ``correct'' the Reidemeister
torsion of the Floer complex also has a better-known analog on the
complex side: the
holomorphic analytic torsion can be viewed as the ``correction term'' to the
\(L^2\)-metric on the determinant line bundle of sheaf cohomologies,
and the combination of the two defines the invariant Quillen metric.
(This above observation was due to Fukaya \cite{FOOO}.) 

\mysubsubsection[Higher torsions in Floer theories.] In \cite{igusa},
Igusa defined higher Franz-Reidemeister torsions via parametrized
Morse theory. The bifurcation
analysis carried out in this paper enables one to lift Igusa's work
to the Floer-theoretic setting. In principle, these Floer-theoretic
higher torsions may be useful for investigating homotopy groups of the
symplectomorphism groups. We hope to return to this subject in the future.

\subsection{Bifurcation Analysis in Floer Theories}
Our method of invariance proof is very different from the traditional
one. 

Hitherto the only available tool for invariance proof in Floer
theory is the ``continuation method'' (cf. e.g. \cite{floer.cmp}),
which relies on the construction of chain homotopy equivalences
between Floer complexes. The drawback of this method is that it is only
useful for proving the invariance of {\em homotopy} invariants such as
homology or \(A_\infty\) categeory. It does not apply for {\em
  non-homotopy} invariants such as the Reidemeister torsion considered
in this paper. 

A natural and straightforward approach, which would apply in more general
situations, is direct bifurcation analysis. Namely, consider a generic
1-parameter family of Floer complexes, classify the possible
bifurcations, and study how the proposed invariant change at these
bifurcations.

In his first paper on Floer homologies \cite{floer.jdg}, Floer
gave a very brief outline of this approach in the case when the action
functional is globally defined. The details were however largely
missing, and this method was subsequently superceded by the
continuation method mentioned above. 
During the long gestation period of the present article, some recent papers
following Floer's original approach have appearred, see
e.g. \cite{Sul, EES}. As noted earlier, these papers, including Floer's, rely
ad hoc methods perculiar to the special
cases considered. (In particular, they do not apply in the situation
considered here). 

In this article, we present a fairly general invariance proof which
carries out this direct approach in its entirety. This proof is
amenable to adaptations to other versions of Floer theories, see e.g. \cite{L2}.
In addition to openning 
the door for the study of other more refined Floer-theoretic
invariants, it also offers a viable alternative to the conventional
continuation method in more intricate variants of Floer homology,
as \cite{EES} demonstrates.

Since we work in the general Morse-Novikov situation, where 
the action functional may {\em not be
globally defined}, there are several important new features which did
not appear in previous works by Floer and other
authors, making the analysis substantially harder
than the situation previously considered. Further remarks follow in
\S1.2.1, 1.2.2.

\mysubsubsection[Finite dimensional model: the Morse-theoretic
picture.] 
It is helpful to first understand 
the simpler case of Morse theory on a finite dimensional manifold.

Given a generic 1-parameter family of Morse functions,
there are two basic types of
bifurcations: a ``handle-slide'', namely a flow line between two
nondegenerate critical points of the same index, or a ``death-birth'', namely two
critical points cancelling each other. 

Moreover, it is not hard to understand how the Morse complexes change
at a bifurcation point, thanks to the geometric interpretation of 
flow lines as intersection points (of
ascending/descending manifolds, or in the case of closed orbits, as 
the fixed points of a local Poincar\'{e} return map). 

For real-valued
Morse theory, the Morse complex changes by an elementary
transformation at a handle-slide, and at a death-birth, the Morse
complex changes by an expansion/collapse modulo elementary
transformations. (See Remark 4.4.3(b) for terminologies). 
This shows that not only the homotopy type of 
the Morse complex is invariant,
but its {\em simple-homotopy} type is invariant. 

For Morse theory of
of closed 1-forms, the bifurcation behavior is more complicated, but
has been partially worked out in \cite{Hut}.
We list 
some major differences between the real-valued Morse theory and the
Morse-Novikov theory of a non-exact closed 1-form below.
\begin{enumerate}
\item Unlike the case with real-valued Morse theory, in this case the
  moduli spaces of flows between two fixed critical points are {\em
    noncompact}. To obtain a well-defined Morse complex in this
  setting, one needs to work with Novikov coefficients. Of crucial
  importance here is a filtration, which depends on the cohomology
  class of the Morse 1-form, modulo rescaling by positive numbers. 
(The truncated moduli space of flows with an upper bound on energy 
is compact). The 1-parameter family
of Morse 1-forms \(\{\theta_\lambda\}_{\lambda\in \Lambda}\)
used for the invariance proof should have the same
filtration, namely, their cohomology classes should be ``co-directional'' in
the sense that \([\theta_\lambda]=\alpha_\lambda[\theta]\) for a fixed
\([\theta]\in H^1\), and a family of non-negative numbers \(\{\alpha_\lambda\}_{\lambda\in\Lambda}\).
\item The flow lines may form closed orbits. At a bifurcation point,
  infinitely many new isolated flow lines between critical points, or
  closed orbits, may be generated by gluing arbitrarily many 
flow lines simultaneously. 
\item There can be infinitely many handleslides in a family
of Morse-Novikov theory parametrized by a compact interval \(\Lambda\subset
\R\). As the Morse complex is undefined on a possibly dense subset
of \(\Lambda\), even the notion of the ``change'' of the Morse complex
at a bifurcation point requires careful definition. 
\end{enumerate}

At a handleslide bifurcation where a flow line starts and ends at the
same critical point (called ``type II handleslide'' in this article),
the higher order (in terms of energy filtration) bifurcation behavior
is not understood. In \cite{Hut}, Hutchings side-stepped this problem by
considering the induced flow on finite-cyclic coverings of the
manifold, using the relation between the torsion of \(M\) and of its
finite-cyclic coverings to reduce the higher order problem to the low
order problem. We shall follow his approach.

\mysubsubsection[From Morse theory to Floer theory.] While
Morse theory gives a nice model for Floer theory, there are several
important differences. For example: 

\begin{enumerate}
\item The grading of a Floer complex often takes values not in the
  group \(\Z\), but rather in \(\Z/N\Z\) for some \(N\in \Z^+\). In order
  for the Reidemeister torsion to be well-defined, \(N\) has to be
  even. It is not hard to see that this is the case for the version of
  Floer theory considered in this article. 
\item The geometric interpretation of flow lines mentioned in \S1.2.1
is no longer available in the infinite dimensional context of Floer
theory. Instead, the spaces of flow
lines need to be described as moduli spaces of certain elliptic PDEs.
The verification of the bifurcation behaviors predicted by Morse theory
requires detailed understanding of the behaviors of moduli spaces under
bifurcation, which relies on certain gluing
theorems. Proving these gluing theorems constitutes the major
difficulty in this direct approach of invariance proof.
(In comparison, the traditional continuation method only requires certain
characteristic numbers of the relevant moduli
spaces to be well-defined). Part II will deal with this problem. 
\item While a Morse function admits rather flexible perturbations,
  only very restrictive perturbations to a Floer theory make good moduli
problems. Consequently, some trivial facts in Morse theory become
highly nontrivial in a Floer theory. As a well-known instance, the
transversality proofs in Floer theory can be rather involved, unlike
the case of Morse theory. Two instances relevant to this article are:

(a) To apply Hutchings's argument for type II handleslide bifurcation,
one needs nonequivariant perturbations to the induced flow on
finite-cyclic covers of the loop space. (See
\S4.4.5 for a more precise statement). In the Floer theory context, it
is highly nontrivial to construct such perturbations satsifying the
desired properties. Unless one restricts to special classes of
symplectic manifolds, the perturbation will be {\em nonlocal}, and new
arguments are required to establish the usual transversality and
compactness properties of the moduli spaces. See section 6 of Part II.

(b) As explained in item 1 of \S1.2.1 above, the 1-parameter family of Morse
1-forms used in the invariance proof must be ``co-directional''. 
Nonetheless, the finite dimensional invariant \(I\) is independent of the 
Morse 1-form used for its definition, {\em regardless of its cohomology
  class}. This is because any Morse 1-form may be connected to an
exact one via a path of co-directional closed 1-forms. In Floer theory,
it is often difficult to find a path which both satisfies the
co-direction condition and makes a good moduli problem.  
In our Floer-theoretic context, the analog of the Morse 1-form  is the
``action 1-form'' \({\cal Y}_X\). Perturbing the symplectic vector
field \(X\) by Hamiltonian vector fields does not change the
cohomology class \([{\cal Y}_X]\), but perturbing by a general
symplectic vector field does. While one expects \(I_F\) to be invariant under
Hamiltonian isotopies, whether it is invariant under
general symplectic isotopies depends on whether one may find a
symplectic vector field \(X_0\) such that \([{\cal Y}_{X_0}]=0\). 
The existence of such \(X_0\) is guaranteed by
imposing strong monotonicity conditions on the symplectic manifold;
see Theorem \ref{thm:A} (b) for a precise statement. It is not known
whether \(I_F\) is invariant under symplectic isotopies in general.
\end{enumerate}

\mysubsubsection[Adaptability of the method.]
We now briefly indicate which portion of the arguments in this article is
independent of the specific Floer theory considered. By a `Floer
theory', we mean one which satisfies the properties outlined in
\S2.1 below.

The following {\em depend} on the specific Floer theory: 
part of the transversality and compactness arguments, such as 
the structure theorems of the moduli spaces proven in section 3,  
and the arguments in section 6, and the construction of nonequivariant
perturbations in \S6.3 of Part II; the orientability of moduli spaces, as discussed in \S7.2 of Part II.

The following are {\em universal} for all Floer theories: 
the main framework of proof for the general invariance theorem, 
Theorem \ref{thm:gen-inv} below; the analysis for
moduli space of flows ending at a degenerate critical point contained
in section 5;
the main outline of the proofs of gluing theorems in sections 2--5 in
Part II applies to any version of Floer theory; the coherentness of 
orientations of moduli spaces discussed in \S7.3 of Part II, which
rely on linear versions of gluing theorems. 
The estimates in the proofs of these gluing theorems require the
following additional conditions, which hold for most of the existent
Floer theories:
\begin{itemize}
\item Recall that a Floer theory is constructed from the formal flow of a
  (densely defined) vector field \({\cal V}\) on a Banach manifold
  \({\cal C}\) modelling on a function space. In local coordinates, 
we need \({\cal V}\) to
  depend only on the function itself, not its higher derivatives, at
  least in a neighborhood of the degenerate critical points. 

In symplectic Floer theories, this condition holds by the assumption
of the degenerate critical points being in ``standard d-b
neighborhoods'' (cf. \S5.3). In gauge theories, it follows directly
from the defining formula for \({\cal V}\) that this condition holds
anywhere on \({\cal C}\).

\item We use the Sobolev embedding theorem to obtain \(C_0\) estimates
  from \(L^p\)-estimates on several occasions. The applicability of
  Sobolev embedding depends on the dimension \(1+l\), when \(T\lc\) is
  a space of functions/sections over an \(l\)-dimensional manifold. In
  this article, as well as in other symplectic Floer theories,
  \(l=1\). In a typical gauge theory, \(l=3\). In this case we need to
  require \(p\geq 4\) for these specific estimates to work.
\end{itemize}

We shall also consider in \cite{L2} an equivariant version of Floer
theory, that does not entirely fit in the framework of \S2.1.2. 
Nevertheless, a minor modification of the methods of this article
still serves the purpose. 

\mysubsubsection[Outline of Part I.] In section 2,
we summarize the abstract framework for defining the invariant
\(I_F\), and give precise statements of the main results. 
The concrete Floer theory considered in this article is set up in
section 3, where structure theorems of the relevant moduli spaces 
needed for the construction of \(I_F\) are established. 
In section 4, we state a general invariance theorem, Theorem
\ref{thm:gen-inv}, and show how it implies the theorems stated in section 2.
Based on the finite-dimensional Morse theoretic picture, we summarize
the expected bifurcation behavior of a generic 1-parameter family of
Floer systems into the notion of a ``regular homotopy of Floer
systems'' (RHFS). Following \cite{Hut}, it is shown that the existence
of an RHFS with the extra property (NEP) implies the general invariance
theorem \ref{thm:gen-inv}. Part II and the rest of Part I are devoted
to establishing the existence of such an RHFS. 

Since a degenerate critical point appears at a death-birth
bifurcation, we need to describe the structure of the
  moduli spaces of flows ending at such
  degenerate critical points. This is done in section 5.
Section 6 introduces the notion of an ``admissible \((J,
X)\)-homotopy'', which is equipped with most of the properties of an
RHFS. We then establish the existence of an
admissible \((J, X)\)-homotopy.

Basing on a series of hard gluing theorems, Part II will show that an
admissible \((J, X)\)-homotopy also possess the remaining properties
of an RHFS. There the issue of orientation will also be addressed, thus
concluding the proof of Theorem \ref{thm:gen-inv}. See section 1 of
Part II for a summary of these results, and an outline of
the methods of proof.

The main technical components of Part I are contained in sections---
\begin{description}  
\item[3.3,] in which we establish transversality for moduli
  spaces of closed orbits in homology class \(A\),
  \(\hat{\cm}_O^0(A)\), {\em simultaneously} for all \(A\). This
  depends on an adaptation of Taubes's argument in \cite{T2}, which
  translates the problem of transversality for multiple covers into
  the simultaneous surjectivity of a sequence of differential
  operators over simple orbits. The knowledge of the kernels and
  cokernels of these differential operators as representation spaces
  of finite-cyclic groups helps to solve the latter problem. 
\item[5,] in which we establish the polynomial decay of flows
  ending at a degenerate critical point, and the Fredholm theory for
  the moduli of such flows. The decay estimates employs center manifold
  theory for flows on Hilbert manifolds, modelling on the discussion
  in \cite{MMR}. To find the suitable Fredholm framework, 
the appropriate spaces for the domain and range of the
  deformation operator turn out to be somewhat complicated, 
due to the following multiple constraints:\begin{itemize}
\item the deformation operator must be Fredhom between these spaces;
\item the domain must be big enough to contain the moduli space
  considered;
\item the nonlinear part of the PDE must satisfy the expected
  quadratic bound.
\end{itemize}
They are polynomially weighted Sobolev spaces, but the weights on the
  ``transversal directions'' and the ``longitudinal direction'' are
  different. The rationale for the choice is that the deformation
  operator is modelled on the operator \(d/ds +C\), where \(C\) is a
  constant that is nonzero for the transversal directions, and zero
  for the longitudinal directions.
\item[6.2,] where we show how to perturb a 1-parameter family
  of Floer systems into an admissible \((J, X)\)-homotopy, basing on a
  refinement of the standard transversality arguments. 
\end{description}
Due to the noncompactness of moduli spaces (cf. e.g. \S1.2.1 items 2
and 3), careful filtration arguments are frequently called for
throughout the length of this article. See in particular \S3.3, 6.2.

\subsection{Notation and Conventions}
The following notation and conventions are adopted in both parts of
this article. References to section or equation numbers in part II
will be denoted as II.*.
{\it We suggest the reader to first browse through this subsection and
  to return later for reference of the notations.}
\begin{itemize}\itemsep -1pt
\item \(C, C', C_i\) etc. usually denote positive constants depending
  on the context. \(\varepsilon, \varepsilon'\) usually denote small
  positive numbers. In contrast, the plain \(\epsilon\) usually
  denotes a fixed parameter; for example, in section 3 and 5 
it parameterizes the
  weight in the exponentially weighted Sobolev spaces; in II.2--II.5 it
  denotes a fixed small number associated with the choice of partition
  of \(\Theta\).
\item\(\Pi_V\) in general denotes a projection. It denotes the
  projection to the space \(V\) if \(V\) is a space; or to the
  direction of \(V\) if \(V\) is a vector.
\item \(S^1_T=\R/T\Z\) denotes the circle of length \(T\).
\item The top exterior power of a vector space \(V\) is denoted by
  \(\det V\).
\item\(\|\xi\|_{p,k}\) denotes the Sobolev norm obtained by summing
  \(L^p\) norms of derivatives of \(\xi\) up to order \(k\). \(L^p_k\)
  denotes the associated Sobolev space. Throughout this article, \(p\)
  is an integer \(p>2\).
\item \(C_\epsilon^{\infty}=C_\epsilon \subset C^{\infty}\) denotes
  the Banach space endowed with the
  \(\|\cdot\|_\epsilon\)-norm defined in \cite{floer.cpam}.
\item  Let \(L\) be any of the norms used in this article, and \(F\)
  be an Euclidean/hermitian vector bundle over \(M\). 
\(L(M)\) denotes the space of functions on the manifold \(M\) with
finite \(L\)-norm. \(L(M; F)=L(F)\) denotes the space of sections of
\(F\) with finite \(L\)-norm.
\item In this paper, the inner product \(\langle\cdot, \cdot\rangle_2\),
  and the norms \(\|\cdot\|\) are usually for functions or sections
  over the cylinder or torus \(\Theta=\{(s,t)| t\in S^1\}\), \(s\in \R\) or \(S^1_T\). We
  will occasionally encounter restrictions of these functions/sections
  to a circle of fixed \(s\) in \(\Theta\). The inner product or norms
  of such restrictions are denoted by \(\langle\cdot,
  \cdot\rangle_{2,t}\), \(\|\cdot\|_{p,t}\) etc.
\item Let \(\xi\) be a section over \(\Theta\). \(\p_s, \p_t\) denote
  the covariant derivative with respect to a natural connection,
  e.g. that induced from the Levi-Civita connection. We also often use
  the short hand \(\dot{\xi}:=\partial _t\xi\); \(\xi':=\partial_s \xi\).
\item Many functions in this paper depend on both the position of a
  point in the symplectic manifold \(M\), and other variables such
  as \(s, t, \lambda\). However, the notation \(\nabla\) is always
  taken to mean the gradient as a function of \(M\), holding other
  variables constant. Derivatives in the directions of \(s, t,
  \lambda\) are denoted \(\partial_s, \partial_t, \partial_\lambda\) respectively.
\item Let \(R\) be a commutative ring. 
\(Q(R)\) denotes the total ring of fractions (also called quotient
ring) of \(R\). Namely, the localization at all non- zero-divisors of \(R\).
\item Let \(V_1, V_2\) be two isomorphic oriented spaces. Then
  \(V_1/V_2\in \{1, -1\}\) denotes the relative sign of them. 
If \(V_1/V_2=-1\), we write \(V_1=-V_2\). If
  \(V_1\) is 1-dimensional and \(a\) is a nontrivial element of \(V_1\),
  then \([a]/V_1\) denotes the sign of \(a\)
  with respect to the orientation of \(V_1\). \(\sign (u):=u/|u|\) for
  \(u\in \R\).
\item We follow the convention of calling a point in a moduli space
  (Zariski) {\em smooth} or {\em
    nondegenerate} if the relevant deformation operator at that point
  is surjective. A moduli space is said to be (Zariski) {\em smooth} or {\em
    nondegenerate} if it consists of smooth points. 
\item Following a well-known convention, a vector field on a Banach
  space refers to one that is only densely defined. 
\end{itemize}

More notation and conventions will be introduced along the way in later sections.

\section{Statements of Main Results.}
Some preliminaries are required for the precise statements of our
results. 

\subsection{The Floer Theory Package.}
As our construction of Floer-theoretic torsions applies to many
versions of Floer theories, we give here a general outline of the abstract
Floer-theoretic framework needed for this construction. The concrete
content of the specific version of Floer theory considered in this
article will be described in section 3. 
The main difference between the following discussion and those in the
existent literature consists of the requirement of an absolute
\(\Z/2\Z\)-grading, the attention to the moduli spaces of closed
orbits, and the emphasis on ``grading-compatible'' orientation for the
moduli spaces of closed orbits.

\subsubsection{Basic ingredients: \((\lc, \mathfrak{H}, \ind;
  \mathcal{Y}_\chi, \mathcal{V}_\chi)\).}
In a typical Floer-Novikov theory, one has an (infinite-dimensional)
Banach manifold \(\lc\), and a closed 1-form \({\cal Y}_\chi\) on \(\lc\), which
may depend on certain parameter \(\chi\). Together with an
\(L^2\)-metric on \(\lc\) (also possibly dependent on \(\chi\)),
\({\cal Y}_\chi\) determines a (densely defined) dual vector field
\(\mathcal{V}_\chi\) on \(\lc\). This defines a {\em formal flow} on
\(\lc\) in the following sense: {\em The moduli space of critical points}, \({\cal
  P}_\chi\subset \lc\) is defined as the zero locus of
\(\mathcal{V}_\chi\). A {\em flow line} \(u\) is a solution to the {\em
  elliptic} PDE, \(\partial_s u+\mathcal{V}_\chi(u)=0\), where \(s\in
\R\) or \(S^1\). It is called a {\em closed orbit} when \(s\in
S^1_T\), for some \(T>0\), and \(T\) is said to be the {\em period} of
the closed orbit. 
When \(s\in\R\) and \(u(s)\) approaches critical points as \(s\to
\pm\infty\), it is called a {\em connecting flow line}. The moduli spaces
of non-constant 
closed orbits and the moduli spaces of connecting flow lines will
be denoted, respectively, by \(\cm_O(\chi)\) and \(\cm_P(\chi)\).
A closed orbit is typically written as either \(u\) or \((T, u)\),
when one wants to emphasize its period.

In order for the Floer-theoretic torsion to be defined, 
we consider Floer theories endowed with an {\em absolute
  \(\Z/2\Z\)-grading}: Let \({\cal P}_{\chi, ndg}\subset {\cal
  P}_\chi\) be the set of nondegenerate elements. There is a map 
\[\ind: {\cal P}_{\chi , ndg}\to \Z/2\Z \quad \forall \chi\] defined by the spectral flow from certain standard
operator to the linearization of \({\cal
  V}_\chi\) at \(x\in {\cal P}_{\chi, ndg}\), denoted \(A_x\). 

Let \(\tilde{\lc}\) denote the regular covering of \(\lc\) with
an abelian covering group \(\mathfrak{H}\) and
monodromy homomorphism \(\op{im}:
\pi_1(\lc)\to \mathfrak{H}\).
\begin{notation*}
We shall frequently use the following notation for elements in a
covering space. Let \(\gamma_0\in \lc\) be a base point. 
Identify $\tilde{\lc}$ with the set \(\tilde{\lc}^{\gamma_0}\) of equivalence classes
$(x,[w])$, where $x\in \lc$, $\mu$ is a map $[0,1]\to\lc$,
$\mu(0)=\gamma_0; \mu(1)=x$,
and two pairs $(x,[w]), (x',[w'])$ are equivalent iff $x=x'$, and
$\op{im}[\mu-\mu']=0$. 
\end{notation*}

We consider Floer theories in which \(\mathfrak{H}\) can be chosen so
that for {\em any} \(\chi\),
\begin{itemize}
\item \(\mathcal{Y}_\chi\) lifts to an exact form \(d{\cal
    A}_\chi\). Namely, the cohomology class
  \([\mathcal{Y}_\chi]=\underline{\op{im}}^*\langle\mathcal{Y}_\chi\rangle\), where \(\langle\mathcal{Y}_\chi\rangle\in \op{Hom}(\mathfrak{H}, \R)\), and \(\underline{\op{im}}\) is the abelianization of \(\op{im}\).
\item The absolute \(\Z/2\Z\)-valued grading lifts to a relative
  \(\Z\)-valued grading: 
Let \(\tilde{\cal P}_\chi\subset \tilde{\lc}\) denote the lift of
\({\cal P}_\chi\). 
There is a map 
\(\op{gr}: \tilde{\cal P}_{\chi, ndg}\times \tilde{\cal P}_{\chi,ndg}\to \Z\),
so that for any \((x, [w]), (y, [v]),
(z, [r])\in \tilde{\cal P}_\chi\), 
\begin{gather*}
\op{gr}((x, [w]), (y, [v]))+\op{gr}((y, [v]), (z,
[r]))=\op{gr}((x, [w]), (z, [r])); \\
\op{gr}((x, [w]), (y,
[v])) \mod 2=\ind (x)-\ind (y), \quad \text{and}
\end{gather*}
there is a homomorphism (which we call {\em the SF-homomorphism})
\[\psi: \mathfrak{H}\to 2\Z\]
with the following property:
Let \((x, [w])\), \((x,
[w'])=A\cdot(x, [w])\) be different lifts of the same \(x\in {\cal
  P}_\chi\), where \(A\in \mathfrak{H}\) acts by deck transformation.
Then \begin{equation}\label{eq:ind}
\op{gr}((x, [w']), (x, [w]))=\psi(A).\end{equation}
\end{itemize}
Via index theory, \(\psi \) is typically computable from the
topology of the underlying manifold and is independent of the
parameter \(\chi\). However, in a general Floer theory it is only
expected to take value in \(\Z\) instead of \(2\Z\), as we require here.

Sometimes weighted versions of \(\ind\) and \(\op{gr}\) are also
needed: For \(\sigma, \sigma_1, \sigma_2\in \R\), \(\ind^\sigma(x)\) and \(\op{gr}^{(\sigma_1, \sigma_2)}((x_1,
[w_1]), (x_2, [w_2]))\) are defined similarly to \(\ind\) and
\(\op{gr}\), but with the role of \(A_{x_*}\) replaced by
\(A_{x_*}+\sigma_*\).

We shall often denote a Floer theory by the 5-tuple of its basic
ingredients: \((\lc, \mathfrak{H}, \ind; \mathcal{Y}_\chi,
\mathcal{V}_\chi)\). Note that the first three items are independent
of the parameters, while the last two do.

\subsubsection{Variants of moduli spaces.}
We now introduce some notation and terminologies frequently used in this
article.  Throughout this subsubsection, we work with a fixed
parameter \(\chi\). Therefore, though all the moduli spaces below depend on
the parameter, we shall omit \(\chi\) from the notation.

Let \(\cm_P(x, y)\) denote the moduli space of connecting flow lines
starting from the critical point \(x\) and ending at the critical
point \(y\). 
Given \((x,[w])\), \((y,[v])\subset\tilde{\cal P}_\chi\),
\(\cm_P((x,[w]), (y,[v]))\subset\cm_P(x, y)\) denote the subset of
elements which lift to a path in
\(\tilde{\lc}\) starting from \((x,[w])\) and ending in \((y,[v])\).

Given \(A\in \mathfrak{H}\), \(\cm_O(A)\subset \cm_O\) denotes the
set of closed orbits with homotopy class in \(\op{im}^{-1}(A)\).

The {\em grading} of a closed orbit \(u\in \cm_O(A)\) is given by
\(\psi(A)\). 
For a connecting flow line \(u\in \cm_P((x,[w]),
(y,[v]))\), it is given by \(\op{gr}((x,[w]), (y,[v]))=k\).
Let \(\cm_O^{k+1}\subset \cm_O\), \(\cm_P^k(x, y)\subset \cm_P(x, y)\)
denote the subset consisting of
elements with grading \(k\).

Weighted versions of moduli spaces are occasionally needed. Roughly
speaking, \(\cm_P^{(\sigma_1, \sigma_2)}\subset \cm_P\) consists of
connecting flow lines decaying exponentially at \(\pm\infty\) ends,
with the exponent \(>\sigma_1\) in the negative end, and \(>-\sigma_2\)
in the positive end. \(\cm_P^{k,(\sigma_1, \sigma_2)}\) consists of
such connecting flow lines with the weighted grading \(\op{gr}^{(\sigma_1, \sigma_2)}=k\).

The {\em energy} of an element \(u\in \cm_O(A)\) is given by
\(-\langle{\cal Y}_\chi\rangle(A)\); for an element 
\(u \in\cm_P((x,[w]), (y,[v]))\), it is given by 
\(\tilde{\cal A}_\chi(x, [w])-\tilde{\cal A}_\chi(y,[v])\).

We often need truncated versions of moduli spaces: 
Given \(\Re\in \R\), let \(\cm_O(\chi)^{\Re}\subset \cm_O(\chi)\),
\(\cm_P(\chi)^\Re\subset \cm_P(\chi)\) be the subsets consisting of elements with energy \(\leq\Re\).
Similarly for other variants of \(\cm_P, \cm_O\).

Note that there is a free \(\R\)-action on \(\cm_P\) and a semi-free
\(S^1\)-action on \(\cm_O\) by translation. The {\em reduced moduli
  spaces} \(\hat{\cm}_P\), \(\hat{\cm}_O\), are respectively the
quotient spaces of \(\cm_P\), \(\cm_O\) under this action; similarly
for other variants of \(\cm_P\) and \(\cm_O\) introduced above. In
particular, \(\hat{\cm}_P^k=\cm_P^{k+1}/\R\);
\(\hat{\cm}_O^k=\cm_O^{k+1}/S^1\). An element
in the reduced moduli space is called a {\em reduced flow line}. 
Given \(u\in \cm_P\) or \(\cm_O\), we use \(\hat{u}\) to denote the
corresponding element in the reduced moduli space; conversely, \(u\)
is said to be a representative of \(\hat{u}\).

Given an element \(u\in \cm_O\), \(\op{mult}(u)=\op{mult}(\hat{u})\)
denotes the {\em multiplicity} of \(u\). Namely, the order of the stabilizer
at $u$ under the \(S^1\) action. A closed orbit of multiplicity 1 is
said to be {\em simple}; otherwise it is said to be a {\em multiple cover}.

It is often convenient to identify the reduced moduli spaces with a
slice of representatives in the un-reduced version. In particular, we
often identify \(\hat{\cm}_P\) with the
set of centered elements in \(\cm_P\), where 
an element \(u\in \cm_P((x, [w]), (y, [v]))\) is {\em centered} if 
\[
2\tilde{\cal A}_\chi(u(0))= \tilde{\cal A}_\chi(x, [w])-\tilde{\cal A}_\chi(y, [v]).
\]
The moduli spaces \(\cm_P(x, y)\), \(\cm_O\) embed respectively into
certain Banach manifolds \({\cal B}_P(x, y)\), \({\cal B}_O\). They
and their reduced versions are thus endowed with the ambient topology,
which we call {\em \({\cal B}\)-topology}.
The reduced moduli spaces \(\hat{\cm}_P\), \(\hat{\cm}_O\) are
included in the following bigger spaces: the moduli space of broken
trajectories, denoted \(\hat{\cm}_P^+\), and the moduli space of
broken orbits, denoted \(\hat{\cm}_O^+\). The latter consist of the
following more general objects:

A (\(k\)-th) {\em broken trajectory} is an ordered set \(\{\hat{u}_0,
\hat{u}_1, \ldots, \hat{u}_k\}\), where \(\hat{u}_i\) are reduced
connecting flow lines, with the end point of \(\hat{u}_i\) identified
with the starting point of \(\hat{u}_{i+1}\) \(\forall i\in \{0,
\ldots, k-1\}\). These critical points are said to be the {\em
  connecting rest points} of the broken trajectory.  
The starting point of the broken trajectory is defined to be the
starting point of \(\hat{u}_0\), while its end point is defined to be
the end point of \(\hat{u}_k\). With the notion of starting and end
points of a broken trajectory clarified, the space \(\hat{\cm}^+_P(x, y)\) now
makes sense. The connecting flow line 
\(\hat{u}_i\) are said to be the \(i\)-th {\em components} of the
broken trajectory. In particular, any connecting flow line is a 0-th
broken trajectory. 

Similarly, when \(k\in \Z^+\), 
a \(k\)-th {\em broken orbit} is a cyclically ordered set \(\{\hat{u}_1,
\hat{u}_2, \ldots, \hat{u}_k\}\), where \(\hat{u}_i\) are reduced
connecting flow lines, with the end point of \(\hat{u}_i\) identified
with the starting point of \(\hat{u}_{i+1}\) \(\forall i\in \Z/k\Z\);
a 0-th broken orbit is by definition a closed orbit. 
Given representatives \(u_i\) of the components \(\hat{u}_i\),
The homology class of the k-th broken orbit above is defined to be
the homology class of the 1-cycle \(\sum_{i=1}^ku_i\) in \(\lc\).
With this explained, the definition of \(\hat{\cm}_O(A)\) easily
extends to give \(\hat{\cm}_O^+(A)\).

The spaces of broken trajectories/orbits are endowed with the {\em
chain topology}. They are stratified
spaces, with the strata indexed by the set of rest points (regarded as
an ordered set for broken trajectories, and as an cyclically ordered
set for broken orbits). A stratum
indexed by \(S_1\) is in the closure of the stratum indexed by \(S_2\)
iff \(S_2\subset S_1\) as ordered or cyclically ordered
subsets. 

The notion of grading and energy extend naturally to the moduli spaces
of broken trajectories or broken orbits. Thus, we also have
\(\hat{\cm}_P^{k, +}(\chi)\), \(\hat{\cm}_P^{+}(\chi)^\Re\),
\(\hat{\cm}_P^{k, +}(\chi)^\Re\), \(\hat{\cm}_O^{k, +}(\chi)^\Re\)
etc---the moduli spaces of broken trajectories or orbits with the
indicated grading and/or energy bound.  

Generically, \(\hat{\cm}^+_P(x, y)\) or \(\hat{\cm}^+_O\)
are expected to be manifolds with corners. As the construction of
\(I_F\) involves
only moduli spaces of expected dimension no more than 1, 
we limit ourselves to the following special case:
\begin{definition*}[lmb]
Let \(\hat{\cm}^+\) be a stratified space 
consisting of smooth strata of dimension no more
than 1, and \(S\subset \hat{\cm}^+\) is a stratum of dimension
0. \(\hat{\cm}^+\)
said to be {\em locally a 1-manifold with boundary along \(S\)}
(abbreviated {\em lmb}), if there
is a homeomorphism \(S\times [0, 1)\) to a neighborhood of
\(S\subset\hat{\cm}^+\), which restricts to a diffeomorphism over
\(S\times (0,1)\).
\end{definition*}

\subsubsection{Floer systems.}

For generic parameter \(\chi\), the moduli spaces are expected to
satisfy the following regularity and compactness properties:
\begin{description}
\item[(FS1)] {\sl (Structure of \({\cal P}_\chi\))} \({\cal P}_\chi\) consists of finitely many
nondegenerate points. 
\item[(FS2)] {\sl (Structure of \(\cm_P(\chi)\))} For any integer \(k\leq 2\) and any \(x, y\in {\cal
    P}_\chi\), \(\cm_P^k(x,y; \chi)\) is a (Zariski) smooth manifold of
  dimension \(k\).
Furthermore, for any real constant \(\Re\), \(\hat{\cm}_P^0(\chi)^\Re\) consists
of finitely many (Zariski) smooth points, \(\hat{\cm}_P^{1, +}(\chi)^\Re\)
is compact, and is lmb along the strata of 1-th broken trajectories.

\item[(FS3)] {\sl (Structure of \(\cm_O(\chi)\))} For any \(k\leq 1\), the
  space \(\cm_O^k\) is a (Zariski) smooth manifold of dimension k.
Moreover, for any real \(\Re\), \(\hat{\cm}^0_O(\chi)^\Re\) consists
of finitely many elements.
\item[(FS4)] {\sl (Orientation)} The moduli spaces \(\cm_P^k(x,y)\),
  \(\cm_O^k\) are orientable; and hence so are their reduced
  versions. Furthermore, the orientations of \(\{\cm_P^k(x,y)\}_{k\in
    \Z; x, y\in {\cal P}}\) can be chosen coherently, and the
  orientation of \(\cm_O^1\) can be chosen to be
  grading-compatible. With respect to this choice,
  \[\partial\hat{\cm}_P^{1, +}(x, y)=\coprod_{z\in {\cal P}}\hat{\cm}_P^0(x, z)\times \hat{\cm}_P^0(z,y).\]
\end{description}
The precise definition of ``coherent orientation'' and
``grading-compatible orientation'' will be postponed to
II.7.2. Roughly speaking, a coherent orientation is a choice of
orientations consistent with gluing, and a grading-compatible
orientation is one which allows a spectral-flow interpretation
compatible with the spectral-flow interpretation of the absolute
\(\Z/2\Z\)-grading \(\ind\). Note also that the definition of
orientability here means the orientability of the relevant determinant
line bundle; hence the moduli space might be orientable even if it is
not a manifold. See II.7 for details.  

\begin{remark*}
The issue of orientation is simpler in the context of Morse theory,
when \(\lc\) is a finite dimensional oriented manifold. In this case, 
the orientation of
\(\lc\) gives a standard choice of the orientations of \(\cm_P\) and
\(\cm_O\), via their interpretation as spaces of intersections (see
e.g. \cite{hl1}). As emphasized in \S1.2.2, item 2, this
interpretation is unavailable in Floer theory.
\end{remark*}

The formal flow associated to \(\mathcal{V}_\chi\)
is said to be a {\em Floer system} if (FS1--FS4) hold.  
An {\em oriented Floer system} is a Floer system together with a
choice of coherent orientation on \(\cm_P\), a grading-compatible
orientation on \(\cm_O^1\), and the induced orientation on the reduced
moduli spaces.

\subsection{Constructing $I_F$: the Algebraic Framework.}

Below is a straightforward adaptation of the construction of \cite{hl2} to
the Floer-theoretic settings.

\subsubsection{The Novikov ring}
Let $G$ be an abelian group, $R$ a ring, and $N:G\to\R$ a homomorphism.
The {\em Novikov ring} $\Nov(G,N;R)$ is the set of formal sums
\(
\sum_{g\in G}a_g\cdot g,
\)
with $a_g\in R$, such that for every $C\in\R$, the set
\(
\{g\in G \mid \mbox{$N(g)<C$ and $a_g\neq 0$}\}
\)
is finite.  $\Nov(G,N;R)$ is a ring with the obvious addition and the
convolution product.  (See e.g. \cite{HS}.) 

Notice that \(\Nov(G, 0, R)=R[G]\), and there is an inclusion
\[i_N: R[G]\hookrightarrow\Nov(G,N;R).\] The Novikov ring should thus be
viewed as a completion of the group ring.

The {\em degree} of \(a\), denoted \(\op{deg}(a)\), is defined to be the minimum of \(N(g)\)
among \(g\) such that \(a_g\neq 0\). (Such minimum exists by the
definition of Novikov rings).

The notion of {\it limit} is defined for Novikov rings: \(\lim_{n\to\infty}
a_n=a\) for a sequence \(\{a_n\in \Nov (G, N;R)\}\) if \(\forall
\Re\in\R\), \(\exists \Lambda(\Re)\) such that \(\deg (a_n-a)\geq \Re\) \(\forall
n\geq \Lambda(\Re)\).

Given \(a=\sum_g a_g g\in\Nov(G,N;R) \), the ``leading term'' of \(a\)
is defined to be
\[
\op{lt}(a):=\sum_{N(g)=\op{deg}(a)}a_g g.
\]
\(a-\op{lt}(a)\) is called ``higher order terms''.  
Notice that \(\op{lt}\) defines a homomorphism 
\[
\Nov(G, N;R)/ (\pm G) \to R[\ker N]/(\pm \ker N).\]

In this article, 
\begin{equation}\label{sp-nov}
\text{\(G\) is a finitely generated abelian group;
\(R=\Z\) or \(\Q\).} 
\end{equation}
In this case, the Novikov ring is commutative.

We shall often need to consider the rings of fractions of Novikov
rings. 

First, observe that a splitting 
\[
G=\ker N\oplus G/\ker N
\]
induces an embedding:
\begin{equation}\label{Q-emb}
Q(\Nov (G,N; R))\hookrightarrow \Nov (G/\ker N, N; Q(R[\ker N])),
\end{equation}
and different embeddings are related by 
the natural action of the space of splittings \(\ker N\) on the right
hand side. 

Furthermore, in the case of (\ref{sp-nov}), both sides of 
(\ref{Q-emb}) are finite sums of fields. (As a special case,
\(Q(R[G])\) is a finite sum of fields, see e.g. \cite{turaev:spinc} \S3.1,
\cite{Hut} Lemma A.4).
The embedding (\ref{Q-emb}) is compatible with the decompositions
on both sides as sums of fields. 

\begin{remark*}
In comparison with the first ring, the second ring in (\ref{Q-emb})
has nicer properties (e.g. existence of the notions of degree, order,
and limit), which the invariance proofs in \cite{hl2, Hut} made use of.
There is a confusion between the two rings in \cite{hl2},
which also propagate to later papers. In these
papers, the notion of order for elements in \(Q(\Nov (G, N, R))\)
should be understood in terms of the larger ring above through the
embedding (\ref{Q-emb}). 
\end{remark*}

Through the embedding (\ref{Q-emb}), we may extend the notion of leading-term to:
\[
\op{lt}: Q(\Nov (G, N; R))/\pm G \to Q(R[\ker N])/\pm \ker N.\]
Note that the above map is independent of the choice of splitting
since we mod out \(\ker N\).
The embedding \(i_N\) also extends to the ring of fractions, 
\[i_N: Q(R[G])\hookrightarrow Q(\Nov(G,N;R)).\]
Later we shall also use the same notation \(i_N\) to denote the induced map from
\(Q(R[G])/(\pm G)\) to \(Q(\Nov (G,N;R))/(\pm G)\).

Let \(\tilde{\Q}\) be a ring extension of \(\Q\), e.g,
\(\tilde{\Q}=\Q\) or the ring \(Q(R[\ker N])\) in (\ref{Q-emb}). 
Let $\Nov^+(G,N;\tilde{\Q})\subset\Nov(G,N;\tilde{\Q})$ denote the subset of elements
of  positive degree.

Let \(\Nov^1(G,N;\tilde{\Q})\subset\Nov(G,N;\tilde{\Q})\) be the
subgroup consisting of elements of the form \(1+c\), \(c\in
\Nov^+(G,N;\tilde{\Q})\). The exponential
\[
\exp:\Nov^+(G,N;\tilde{\Q})\to\Nov^1(G,N;\tilde{\Q})\hookrightarrow Q(\Nov(G,N;\tilde{\Q}))
\] is well defined via the usual
power series. Conversely, the logarithm \[\ln: \Nov^1(G,N;\tilde{\Q})\to
\Nov^+(G,N;\tilde{\Q})\] also makes sense formally.

Novikov rings arise naturally in Morse-Novikov theory as the
coefficient rings.

\subsubsection{The Floer-Novikov complex.}

Each oriented Floer system \((\lc, \mathfrak{H}, \ind; \mathcal{Y}_\chi,
\mathcal{V}_\chi)\) is associated with a (twisted) Floer complex as follows.

Let the chain groups \(\tilde{\op{CF}}\)
be the free $\Nov (\mathfrak{H}, -\langle{\cal Y}_\chi\rangle; \Z)$-module generated by elements in $\tilde{\cal
  P}_\chi$, on which  \(\mathfrak{H}\) acts by
deck-transformation. 

Let the boundary map \(\tilde{\p}_F: \tilde{\op{CF}}\to
\tilde{\op{CF}}\) be defined by 
\[
\tilde{\p}_F (x, [w])=\sum_{(y,[v])\in\tilde{\cal P}}\chi(
    \hat{\cm}_P^0((x, [w]), (y, [v]))) \, (y, [v]).
\]
By (FS2), $\tilde{\p}_F$ is a well-defined $\Nov (\mathfrak{H}, -\langle{\cal Y}_\chi\rangle; \Z)$-linear transformation, and
$\tilde{\p}_F^2=0$.

Via the relative \(\Z\)-grading on \(\tilde{\cal P}_\chi\), 
\(\tilde{\op{CF}}\) splits as:
\[
(\tilde{\op{CF}}, \tilde{\p}_F)=
\bigoplus_{k} (\tilde{\op{CF}}_k, \tilde{\p}_{F, k}),
\]
where \(k\) is in an affine space under \(\Z\).
Each summand $\tilde{\op{CF}}_k$ is a free \(\Lambda_F\)-module of finite rank, where
\[
\Lambda_F
:=\Nov (\ker \psi, -\langle{\cal Y}_\chi\rangle; \Z) \subset \Nov (\mathfrak{H}, -\langle{\cal Y}_\chi\rangle; \Z).
\]
Furthermore, they satisfy the periodicity condition
\[
(\tilde{\op{CF}}_k, \tilde{\p}_{F, k})=(\tilde{\op{CF}}_{k+2{\mathbb
    N}_\psi}, \tilde{\p}_{F, k+2{\mathbb N}_\psi}),
\]
$2{\mathbb N}_\psi$ being the gcd of the values of $\psi$. 

As the relative \(\Z\)-valued index on \(\tilde{\cal P}_\chi\) reduces
to an {\em absolute} \(\Z/2\Z\)-valued grading \(\ind\), there is also
a reduced version of Floer complex, 
\[
(\op{CF}, \p_F)=\bigoplus_{i}(\op{CF}_i, \p_{F,i}),
\]
where \(i\) takes values in the {\em group} \(\Z/2\Z\), and
\(\op{CF}_i\) is the free \(\Lambda_F\)-module generated by all
elements in \({\cal P}_\chi\) of index \(i\).

\subsubsection{The Reidemeister torsion.}\label{def:torsion}
We now specify the version of torsion used in this paper.
\paragraph{(1)} Suppose first for simplicity that the coefficient ring \(F\) is a
field. Let $(C_i,\partial_i),\, i\in \Z/2\Z$ be a complex of finite dimensional
\(F\)-vector spaces, and \(Z_i, B_i\) be respectively the subspace of
cycles and boundaries in \(C_i\). The standard short exact sequences
$0\to Z_i\to C_i\to B_{i-1}\to 0$ and $0\to B_i\to Z_i\to H_i\to 0$
induce a canonical isomorphism 
\[
{\cal T}:\bigotimes_i\det(C_i)^{(-1)^i}\longrightarrow \bigotimes_i\det(H_i)^{(-1)^i}.
\]

Let $e$ be an ordered basis for $C_*$, i.e. an ordered basis $e_i$ for
each $C_i$.  Let $h$ be an ordered basis for $H_*$. 
Let $[e]\in \bigotimes_i\det(C_i)^{(-1)^i}$ and
$[h]\in\bigotimes_i\det(H_i)^{(-1)^i}$ denote the resulting volume forms.

In this simplest case, when the coefficient ring is a field, define the Reidemeister torsion
\[\begin{split}
\tau(C_*;e)& \eqdef\left\{\begin{array}{cl}{\cal T}([e])/[h]\in F^\times&
    \mbox{if $H_*=0$,}\\ 0 & \mbox{otherwise.}\end{array}\right.
\end{split}\]

 \paragraph{(2)} For our applications, the coefficient ring of the
complex is a Novikov ring of the type specified in (\ref{sp-nov})
(including group rings). We saw that in this case it is 
in general not a field, but its total ring of fractions is
a direct sum of fields. \begin{definition*}
\cite{turaev:spinc}
Let $R$ be a ring, and assume that its total ring of fractions $Q(R)$
is a finite sum of fields, $Q(R)= \bigoplus_jF_j$.  Let $(C_i,\partial_i), \,
i\in \Z/2\Z$ be a complex of
finitely generated free $R$-modules with
an ordered basis $e$. Then
\[
\tau(C_*,e)\eqdef\sum_j\tau(C_*\tensor_RF_j,e\tensor 1)\in\bigoplus_jF_j=Q(R).
\]
\end{definition*}

\paragraph{(3)} Associated to an oriented Floer system, the {\em Reidemeister torsion of the Floer complex} is defined as
\begin{equation}
\tau_F:=\tau(\op{CF}, e_{\cal P})\in Q(\Lambda_F)/(\pm \ker\psi),
\end{equation}
where $e_{\cal P}$ is an ordered basis of $CF$ given by an ordering of
elements in ${\cal
  P}$, and a lift \(\op{Lf}: {\cal P}\to \tilde{\cal P}\). 
In the case when \(\mathbb{N}_\psi\neq0\), we require further:
\begin{itemize}
\item that the lift \(\op{Lf}\) is such that \(|\op{gr}(\op{Lf}(x),
  \op{Lf}(y))|<2\mathbb{N}_\psi\) \(\forall x, y\in {\cal P}\), and
\item a choice of an \(A_0\in \mathfrak{H}\), such that \(\psi(A_0)=2\mathbb{N}_\psi\).
\end{itemize}
The set of all \((nA_0)\cdot \op{Lf}(x)\), \(n\in\Z\), \(x\in {\cal
  P}\) then forms a basis of \(\tilde{\op{CF}}\) as a graded free \(\Lambda_F\)-module.

Different ordering of the critical points results in a possible
change of sign for $\tau (CF, e_{\cal P})$, and different such lifts
and \(A_0\) result in a
multiplication of \(\tau\) by an element in \(\ker\psi\);  so by modding out $\pm \ker\psi$ 
in the definition we obtain an invariant independent of these choices.

\paragraph{(4)} Another version of torsion that is important in topology is the {\em
  Reidemeister torsion of a manifold}, denoted \(\tau (M)\). Let \(M\) be a
manifold with a cell-decomposition, such that the cell chain complex 
\(C_*(M)\) is a finite complex of finite-rank \(\Z\)-modules. The
universal abelian covering \(\tilde{M}\) is endowed with an induced 
equivariant cell-decomposition, and \(C_*(\tilde{M})\) is a
\(\Z[H_1(M;\Z)]\)-module. 
\[
\tau(M):=\tau (C(\tilde{M}), e_M)\in Q(\Z[H_1(M;\Z)])/\pm H_1(M;\Z), 
\]
where \(e_M\) is an ordered basis consisting of lifts of cells in \(M\).

\subsubsection{The zeta function and the counting invariant $I_F$.}
Fix an oriented Floer system \((\lc, \mathfrak{H}, \ind; {\cal Y}_\chi, {\cal
  V}_\chi)\).

To count the closed orbits, we imitate 
the definition of the dynamical zeta function, and define the
the Floer-theoretic {\em zeta function} as
\begin{equation}\label{zeta}
\begin{split}
\zeta_F &:=\op{exp}\left(\sum_{A\in
    \ker \psi, A\neq 0} \chi(\hat{\cm}_O(A))\, A\right)\\
&= \op{exp}\left(\sum_{A\in
    \ker \psi, A\neq 0} \sum_{u\in \hat{\cm}_O(A)}\frac{\sign (u)}{\op{mult}(u)}\, A\right)\\
&\quad \in \Nov^1 (\ker \psi, -\langle{\cal Y}_\chi\rangle; \Q)\subset Q(\Nov(\ker \psi, -\langle{\cal Y}_\chi\rangle; \Q)).
\end{split}\end{equation}
In (\ref{zeta}), \(\chi\) is the ``orbifold Euler number'', 
and we shall denote the exponent therein by \(\eta_F\).
By (FS3), \(\eta_F\in Nov^+(\ker\psi,
-\langle{\cal Y}_\chi\rangle; \Q)\), and the exponential is well-defined.

Finally, viewing both \(\zeta_F\) and \(\tau_F\) as elements in
\(Q(\Nov(\ker \psi, -\langle{\cal Y}_\chi\rangle; \Q))/\pm \ker \psi\), \(I_F\) is
simply defined as the product
\begin{equation}
I_F:=\zeta_F \tau_F\in Q(\Nov(\ker \psi, -\langle{\cal Y}_\chi\rangle; \Q))/\pm \ker \psi.
\end{equation}
\begin{remark*}
Ideally, \(I_F\) should be defined in \(Q(\Lambda_F)/\pm\ker\psi\)
instead of the above larger monoid. However, this would require
proving a product formula similar to \cite{hl2} equation (2) to ensure
that \(\zeta_F\in \Nov^1(\ker\psi, -\langle{\cal Y}_\chi\rangle;\Z)\). 
\end{remark*}

\subsection{Statement of the Main Results}

Let $(M,\omega)$ be a closed, connected
symplectic manifold of dimension $2n$,
and let $f:M\to M$ be a smooth, (possibly non-hamiltonian)
symplectomorphism.  

Let \(c_1:=c_1(TM)\). 
We shall need the notion of weak monotonicity and  a stronger
companion version.
\begin{*definition}\label{def:weak-mon}
A symplectic manifold $M$ is {\em  weakly monotone} if $\omega(A)>0$
for all $A\in H_2(M)$
in the image of $\pi_2(M)$ under the Hurewicz map satisfying 
$0<c_1(TM)(A) <n-2$. 

It is said to be {\em \(w^+\)-monotone} if $\omega(A)>0$
for all $A\in H_2(M)$
in the image of $\pi_2(M)$ under the Hurewicz map satisfying 
$0<c_1(TM)(A) \leq n-2$. 
\end{*definition}
It is clear from the definition that monotonicity implies
\(w^+\)-monotonicity, which in turn implies weak-monotonicity. 
Examples of \(w^+\)-monotone symplectic manifolds include Fano and 
Calabi-Yau manifolds, and any symplectic manifold of dimension less
than six.
 
Our main result concerns a ``\(w^+\)-monotone'' symplectic manifold.
The assumption of \(w^+\)-monotonicity is mainly here for simplicity:
In view of \cite{liu-tian,fukaya-ono, ruan}, we expect the result to
hold for general symplectic manifolds.
In fact, most of the proofs contained in this article work for weakly
monotone manifolds. The stronger \(w^+\)-monotonicity assumption is
imposed in section 6 to shorten the discussion on the structure of
parameterized moduli spaces.  

It is convenient to introduce the {\em mapping torus} of $f$:
\[
T_f:=\Big\{(x,t): x\in M, t\in [0,1]\Big\}/(x,0)\sim (f(x),1),
\]
which fibers over the circle of unit length $S^1_1$:
\[
M\stackrel{\iota_f}{\longrightarrow}T_f\stackrel{\pi_S}{\longrightarrow} S^1_1.
\]
\begin{notation*}
By restricting to the fibers of \(\pi_S\), a function \(S\) on \(T_f\)
or a section of a bundle over \(T_f\) corresponds to a path of
functions/sections over \(M\) with matching conditions on the two ends
determined by \(f\). We shall denote these functions/sections over
\(M\) by \(S_t\), \(t\in [0,1]\)
\end{notation*}
\begin{*definition}
The {\em $f$-twisted loop space} $\ls_fM$ (usually just called
the ``loop space'') is the space of $L^2_1$
sections of $T_f$. 

Given \(\gamma_0\in \ls_fM_{\gamma_0}\),
let $\ls_fM_{\gamma_0}=:\lc$ be the path component
containing \(\gamma_0\). 
\end{*definition}
There is a natural map 
\[\op{im}: \pi_1(\lc)\to H_2(T_f;\Z)\] 
sending each 1-cycle $\gamma: S^1\to \ls_fM_{\gamma_0}$ representing a homotopy
class, to the corresponding 
map from $S^1\times S^1$ to $T_f$.

Let $K$ be the subbundle of $TT_f$ consisting of tangent vectors
to the fibers. 

A complex structure \(J\) of the bundle \(K\) is said to be {\it \(\omega\)-compatible}
if \(J_t\) are \(\omega\)-compatible for all \(t\in [0,1]\). Namely,
\(\omega(J_t \cdot,  \cdot)\) gives a riemannian metric on
\(M\)\footnote{This definition differs by a sign from some literature,
  and is convenient for constructing Floer {\it homology} instead of
  Floer cohomology.}. Let
\({\cal J}_K\) denote the Banach manifold of \(C_\epsilon\)
\(\omega\)-compatible complex structures on \(K\). 

Let \({\cal X}=\{X|\, X\in C_\epsilon (T_f, K), \, \text{\(X_t\) is
a symplectic vector field \(\forall t\)} \}\).

Let \({\cal H}=C_\epsilon (T_f)\). Given \(H\in {\cal H}\), let
\(\chi_H\in {\cal X}\) be defined by
\[
(\chi_H)_t=\chi_{H_t}; \quad \omega(\chi_{H_t}, \cdot)=dH_t.
\]

We shall describe in section 3 a Floer theory \((\lc, \mathfrak{H},
\ind; {\cal Y}_\chi, {\cal V}_\chi)\) satisfying the framework
outlined in \S2.1, in which \(\lc=\ls_fM_{\gamma_0}\),
\(\mathfrak{H}=\op{Image}(\op{im})\subset H_2(T_f)\), \(\ind\) is a variant of the
Conley-Zehnder index, and \({\cal Y}_\chi\), \({\cal V}_\chi\) are
parameterized by \(\chi=(J,
X)\in {\cal J}_K\times {\cal X}\). For certain generic \((J, X)\)
which we call ``regular
pairs'', the associated formal flow is an oriented Floer system (see
section 3 and II.7.2).
Thus, for each regular pair \((J, X)\) the construction of \S2.2 
defines an \(I_F\), which we denote by \(I_F^{f, [\gamma_0]}(M; J,
X)\). 

For fixed \(f, [\gamma_0]\), the cohomology class \([{\cal
  Y}_\chi]\) only depends on the flux of \(X\). In particular, \(I_F^{f,
  [\gamma_0]}(M; J, X)\) and \(I_F^{f, [\gamma_0]}(M; J', X')\) take
values in the same monoid if \(X-X'\) is Hamiltonian.  
Our main result concerns the invariance properties of \(I_F^{f,[\gamma_0]}(M; J, X)\).
\begin{*theorem}\label{thm:A}
(a) Let $(M, \omega)$ be a \(w^+\)-monotone symplectic manifold, and let
\(f\) be a symplectomorphism of \(M\) and \(\gamma_0\in \ls_fM\) be
chosen as above. Suppose \((J_1, X+\chi_{H_1}), (J_2, X+\chi_{H_2})\in
{\cal J}_K\times{\cal X}\) are two regular pairs.
Then 
\begin{equation}\label{ham-inv}
I_F^{f,[\gamma_0]}(M; J_1,X+\chi_{H_1})=I_F^{f,[\gamma_0]}(M;
J_2, X+\chi_{H_2}).
\end{equation}

(b) Suppose $(M, \omega)$ is monotone. Then there exists an \(I_F^{f, [\gamma_0]}(M)\in Q(\Z[\ker\psi])/\pm\ker\psi\)
such that 
\begin{equation}\label{sym-inv}
I_F^{f,[\gamma_0]}(M; J,X)=i_{-[{\cal
    Y}_X]}I_F^{f,[\gamma_0]}(M) \quad \text{for any regular pair \((J,
  X)\). }
\end{equation}
\end{*theorem}
We shall review in \S3.1.4 a well-known principle that relates varying
\(f, [\gamma_0]\) by symplectic isotopies and varying the symplectic
vector field \(X\). According to this principle, the above theorem implies
that \(I_F\) is invariant under Hamiltonian or symplectic isotopies of
\(f, [\gamma_0]\), in the case of (a) or (b) respectively. In
particular, under the assumptions of Theorem \ref{thm:A} (b), this
implies that \(I_F\) depends only on the symplectic mapping class of
\(f\) and the conjugacy class of \([\gamma_0]\in\pi_0(\ls_fM)\) under the standard action of
\(\pi_1(\op{Symp}_0(M))\) on \(\pi_0(\ls_fM)\). (\(\op{Symp}_0(M)\) denotes the path
component of the symplectomorphism group of \(M\) 
containing the identity). 

Due to this principle, we shall say that \(I_F\) is {\it invariant
  under Hamiltonian isotopies} when (\ref{ham-inv}) holds; similarly,
we say that \(I_F\) is {\it invariant under symplectic isotopies} when
(\ref{sym-inv}) holds. 

Combining this result with the main theorem in \cite{hl2, Hut},
we can compute $I_F^{f, [\gamma_0]}(M)$ in the following basic case:
\begin{*corollary}\label{cor:C}
Suppose $M$ is monotone, and $f$ is connected to the identity via the
symplectic-isotopy \(f_t\), \(t\in [0,1]\), with \(f_0=\op{Id}\), \(f_1=f\).
Let \(\gamma_0\) be the path \(\gamma_0(t)=f_t(p_0)\) for a base point
\(p_0\in M\). In this case \(\ker\psi=H_1(M;\Z)\oplus \ker
c_1(TM)\Big|_{\pi_2(M)}\); let \(\iota: H_1(M;\Z)\to \ker \psi\)
denote the inclusion, and let \(\iota_*: Q(\Z[H_1(M; \Z])/\pm H_1(M;
\Z)\to Q(\Z[\ker\psi])/\pm \ker \psi\) denote the induced map. Then
\[I^{f, [\gamma_0]}_F(M)=\iota_*\tau(M).\] 
\end{*corollary}
In section 4, we shall state a more general invariance theorem and
explain how Theorem \ref{thm:A} follows as a consequence.

\mysubsubsection[Some Immediate Application.]
Since critical points in the Floer theory correspond to symplectic
fixed points, nontriviality of the Reidemeister torsion \(\tau_F\) will
guarantee the existence of symplectic fixed points. Here is a sample
result in this direction:
\begin{corollary*}
Let \(M, f, \gamma_0\) be as in Corollary \ref{cor:C}. In addition,
suppose that \(\op{lt} \tau(M)\neq 1\). Then \(f\) has a fixed point. 
\end{corollary*}
This is a direct consequence of the computation in Corollary \ref{cor:C} and the
observation that \(\op{lt}(\zeta_F)=1\) by definition.

More applications shall be discussed in subsequential articles
(e.g. \cite{L2, LS, ogw}), including applications exploiting
higher order terms of \(I_F\).
Here we shall content us with the following few remarks.
\begin{remarks*}
(a) The abelian version of torsion introduced in this paper is not
amenable for obtaining quantitative results (e.g. getting
fixed-point-number bounds). A possible refinement is to work with the
nonabelian Whitehead torsion instead: one would then obtain a lower
bound on fixed-point-numbers in terms of the minimal rank of
representations of the relevant Whitehead torsion. This is however
often very difficult to compute, and thus does not seem an effective
way for obtaining such results. Some results in this direction for the
Floer theory of lagrangian intersections (i.e. Example 1.1.3) may be
found in \cite{Gromov-Eliashberg}.

(b) The Morse-theoretic picture leads one to expect a converse of such
existence theorems: that the critical
points of the Floer theory (i.e. symplectic fixed points or Lagrangian
intersection points) can be removed via a Hamiltonian isotopy when 
both the Floer homology and the Floer-theoretic Whitehead torsion
vanish. This is the outstanding ``symplectic s-cobordism conjecture''
posed by Fukaya \cite{fukaya}. 

(c) The requirement that \(M\) is both closed and symplectic places a
very strong constraint on the topology of \(M\), and it is not easy to
find closed symplectic manifolds with interesting Reidemeister torsion
(i.e. the field components of \(\tau (M)\) are not all 0 or 1). We
shall show in the sequels to this paper that interesting examples are
easier to produce in other versions of Floer theories, e.g. those in 
Examples 1.1.1--1.1.3. See also an \(S^1\)-equivariant version in
\cite{L2}. For instance, one may take \(M\) to be a Stein manifold in
Example 1.1.1. By the combinatorial description of Stein manifolds
(\cite{elias}), \(M\) has very flexible topology and it is easy to
produce \(M\) with interesting torsion accordingly.
\end{remarks*}

\section{Floer Theory of Symplectomorphisms.}

Via a simple modification of Floer's work, 
the Floer homology of a general symplectomorphism was first introduced
by Dostoglou and Salamon in \cite{dostoglou-salamon} for
monotone symplectic manifolds. The virtual moduli method being
available today, it should be defined for general symplectic
manifolds. In this paper we however take the middle ground by imposing the
weakly monotonicity condition, to avoid the virtual moduli method. 

Though it is more-or-less standard, there does not seem to be a good
reference that deals directly with the situation we need. We shall
therefore go into some details on the construction of this version
of Floer theory, by indicating how the available literature should be modified.

\subsection{Setup and Basics.}

Recall the definitions and notation from \S2.3.

\subsubsection{Topology of the Loop Space.}

Fix a base point $p_0\in M$, and
let $\gamma_0\in\lc$ be a path with $\gamma_0(0)=p_0$.
We gather some basic facts about $\ls_fM_{\gamma_0}=:\lc$ below.

First, notice that the fiber-bundle structure of \(T_f\) 
\[
M\stackrel{\iota_f}{\longrightarrow}T_f\stackrel{\pi_S}{\longrightarrow} S^1_1
\]
gives rise to the following useful Mayer-Vietoris type sequence:
\begin{equation}\label{M-V}
\cdots H_2(M;\Z)\stackrel{1-f_*}{\longrightarrow}H_2(M;\Z) \stackrel{\iota_{f*}}{\longrightarrow}
H_2(T_f) \stackrel{\p_f}{\longrightarrow} H_1(M;\Z)\stackrel{1-f_*}{\longrightarrow}H_1(M;\Z) \to
H_1(T_f) \cdots.
\end{equation}

Let $c_1^f:=c_1(K)$; since $f$ is a symplectomorphism, $\omega$ defines
a two-form on $T_f$, which we denote by $\omega^f$. The restrictions
of the cohomology classes \(c_1^f\), \([w^f]\) define two
homomorphisms
\(\psi_c: \mathfrak{H}\to \Z\), \(\psi_\omega: \mathfrak{H}\to \R\) respectively.
They are related to the SF-homomorphism and the
cohomology class \([{\cal Y}_\chi]\) respectively.

Next, note the fibration
\begin{equation}\label{weak-fib}
\Omega M\to \ls_fM\stackrel{e_f}{\to} M, 
\end{equation}
where \(\Omega M\) is the loop space (with based point) of \(M\), and
\(e_f(\gamma)=\gamma(0)\) is the end-point map.
Thus, we have the associated homotopy sequence:
\begin{equation}\label{htpy-sq}
\cdots \pi_2(M; p)\to \pi_1(\ls_fM;
\gamma_0)\to\pi_1(M;p)\stackrel{\delta_{f, \gamma_0}}{\longrightarrow}
\pi_1(M; p)\to \pi_0(\ls_fM)\cdots.
\end{equation}

The two exact sequences (\ref{M-V}) and (\ref{htpy-sq}) fit into the
following commutative diagram with \(\op{im}\):
\begin{equation}\label{c-dgrm}
\begin{CD}
\pi_2(M; p_0)@>\rho_f>> \pi_1(\ls_fM;{\gamma_0})@>e_{f*}>> \pi_1(M;p_0)\\
@VV\op{Hurewicz}V @VV\op{im}V @VV\op{Hurewicz}V\\
H_2(M;\Z) @>\iota_{f*}>> H_2(T_f;\Z) @>\p_f>> H_1(M;\Z),
\end{CD}\end{equation}
where ``Hurewicz'' denotes the Hurewicz map.

 \subsubsection{The space of almost complex structures.}

Given a \(\omega\)-compatible almost complex structure $J_0$ over
\(M\), let ${\cal S}(A, J_0)$ be the space of simple
(i.e. nonconstant and not multiply-covered) smooth $J_0$-holomorphic spheres with
homology class $A$ in $M$. The automorphism group of $\C P^1$,
$G:=\mbox{PSL}(2,\C)$, acts freely on ${\cal S}(A, J_0)$ when $A$ is
nontrivial. Let ${\cal S}(J_0):=\coprod_A {\cal S}(A,J_0)$.

A \(J\in {\cal J}_K\) is said to be {\em regular} if the space
\[
S(A, J)/G:=\bigcup_t S(A, J_t)/G
\]
is a smooth manifold of expected dimension (i.e. of dimension \(2n+2c_1(A)-5\)). 

We denote by \({\cal J}_K^{reg}\subset {\cal J}_K\) the subset
consisting of regular elements, and summarize some useful facts on
such regular elements as follows. 
\begin{lemma*}
Let \(M\) be weakly monotone.

(a) If \(J\in {\cal J}_K^{reg}\), then
for all \(t\in S^1\), \(J_t\) is semi-positive in the sense that any
\(J_t\)-holomorphic sphere has non-negative Chern number. 

(b) The subset \({\cal J}_K^{reg}\subset {\cal J}_K\) is Baire. 

(c) Part (b) above may be refined as follows. 
Let \(y: S^1\to T_f\) be a section, and let \(U_y\) be a small tubular
neighborhood of the image of \(y\) in \(T_f\). Let \(\underline{J}\in
{\cal J}_K\), and let \[{\cal J}_K(\underline{J}, U_y):=\Big\{J\, |\, J\in
{\cal J}_K, \, J\Big|_{U_y}=\underline{J}\Big|_{U_y}\Big\},\]
endowed with the obvious Banach manifold structure. 
Then the space of regular elements in \({\cal J}_K(\underline{J},
U_y)\) is Baire. 
\end{lemma*}
\begin{proof}
(a) follows from the proof of \cite{MS} Lemma 5.1.3.
(b) is a result of a standard transversality argument (see e.g. \cite{MS} Theorem
3.1.3). 
To see (c), one simply augments the proof of \cite{MS} Theorem 3.1.3
by the observation that if \(U_y\) is small enough, every
\(J_t\)-holomorphic sphere exits the small neighborhood of \(y(t)\)
obtained by restricting \(U_y\) to the fiber over \(t\).
(See e.g. Remark 5.2 of \cite{trans} for an argument of the same type.)
\end{proof}

Finally, it is useful to introduce the following notation:

Let $M_c(J_t)\subset M$ be the image under the evaluation map
\[
\mbox{ev}: \coprod_{A,\, c_1(A)\leq c}{\cal S}(J_t, A)\times_G S^2\to M, \quad
\mbox{ev}(u, \theta)=u(\theta).
\]
Let $M_c(J)\subset T_f$ be the subspace such that its intersection with
the fiber over \(t\) is \(M_c(J_t)\).

\subsubsection{The Action 1-form and its Associated Formal Flow. }
\begin{notation*}
For \(X\in {\cal X}\), let
\(\theta_{X,t}:=\iota(X_t)\omega\). \(\theta_X:=\int\theta_{X, t} \,
dt\); \([\theta_X]\in H^1(M)\) is the {\it flux} of \(X\). 
In addition, given \(J\in {\cal J}_K\), let \(\check{\theta}_{X,t}^{J_t}\)
be the dual vector field of \(\theta_{X,t}\) with respect to the
metric associated to \(J_t\), and let \(\check{\theta}_{X}^J\) denote
the section \(\{\check{\theta}_{X,t}^{J_t}\}_{t\in S^1}\) of
\(K\). The superscript \(J_t\) or \(J\) will often be omitted when
\(J\) is fixed. 

Conversely, given a closed 1-form \(\theta\) on \(M\), \(X_\theta\)
denotes the symplectic vector field such that \(\theta=\iota(X_\theta)\omega\).
\end{notation*}
In this article, \({\cal Y}_\chi\) and \(\tilde{\cal A}_\chi\) will be the
{\em action 1-form} \({\cal Y}_X\) and {\em action functional}
\(\tilde{\cal A}_X\), given as follows.

For any $X\in {\cal X}$, 
\begin{equation}\label{cal-Y}
{\cal Y}_X (\gamma)(\xi):=-\int_{S^1}\omega (\p_t\gamma(t), \xi(t))\,
dt+\int_{S^1} \theta_{X,t}(\xi(t))\, dt\quad \text{for \(\xi\in T_\gamma\lc=L^2_1(\gamma^*K)\)}.
\end{equation}
Since \(H_2(T_f;\Z)\) is abelian, 
the map \(\op{im}\) factors as a composition of the
Hurewicz map and a homomorphism \(\underline{\op{im}}: H_1(\lc;\Z)\to
H_2(T_f;\Z)\). Similarly, we may define the homomorphism
\(\underline{e}_{f*}: H_1(\lc;\Z)\to H_1(M;\Z)\) so that \(\underline{e}_{f*}\circ
\op{Hurewicz}=\op{Hurewicz}\circ e_{f*}\).
The cohomology class \([{\cal Y}_X]\in H^1(\lc)\) may be expressed in
these homomorphisms as
\[
[{\cal Y}_X]=-\underline{\op{im}}^*[\omega^f]
+\underline{e}_f^*[\theta_X]=\underline{\op{im}}^*
(-[\omega^f]+\partial_f^*[\theta_X]), 
\]
where (\ref{c-dgrm}) was used for the second equality. We see that
\([{\cal Y}_X]: H_1(\lc; \Z)\to \R\) factors through
\(\underline{\op{im}}\) and a \(\langle{\cal Y}_X\rangle\in
\op{Hom}(\mathfrak{H}, \R)\):
\[
\langle{\cal Y}_X\rangle=-\psi_\omega+\partial_f^*[\theta_f].
\]
Hence \({\cal Y}_X\) lifts to an exact
form on the \(\mathfrak{H}\)-cover \(\tilde{\lc}\), meeting the
requirement in \S2.1.

On the other hand, given \(J\in {\cal J}_K\), we have a riemannian
metric on $T_f$ by the $\omega$-compatibility of $J_t$, which in turn defines a metric for $\ls_fM$.

The vector field dual to \(-{\cal Y}_X\) with respect to
this metric is \(-{\cal V}_X\), 
\begin{equation}
{\cal V}_X(\gamma):=J_t(\gamma) \frac{\p \gamma}{\p t}+\check{\theta}_{X,t}^{J_t}(\gamma)
\qquad \text{for \(\gamma\in \lc\).}
\end{equation}
Thus, a formal flow line of \({\cal V}_X\) is a solution to the perturbed Cauchy-Riemann equation
\begin{equation}\label{pcr}
\bar{\p}_{JX} u:=\frac{\p u}{\p s}+J_t(u) \frac{\p u}{\p
  t}+\check{\theta}_{X,t}^{J_t}(u)=0,
\end{equation}
where \(t\in S^1_1\), and \(u\) is a smooth \(M\)-valued
function on \((s,t)\) such that \(u(s, \cdot)\in \lc\, \, \forall s\);
\(s\in \R\) when \(u\in \cm_P\), while \(s\in S^1_T\) for some \(T>0\)
if \(u\in \cm_O\).

The {\em energy} of a solution \(u\) to (\ref{pcr}) is 
\[
{\cal E}(u)=\int|\p_s u|^2 ds\, dt.
\]
A straightforward computation confirms that this agrees with the
definition for energy of closed orbits or connecting flow
lines given in terms of \({\cal Y}_\chi\) and \(\tilde{\cal A}_\chi\)
in \S2.1.

The various moduli spaces in this Floer theory depend on \(\gamma_0\),
\(M\), \(f\), and \(J, X\). All these shall be incoporated in the full
notation for the moduli spaces, e.g. \({\cal M}_P^{\gamma_0}(M, f;J,
X)\). However, when
\(\gamma_0, M, f\) are fixed, they will be often omitted from the
notation. So is the parameter \((J, X)\), if it is clear from the context.

We now introduce the chain topology on
\(\hat{\cm}^+_P\) and \(\hat{\cm}^+_O\). First, let \(p_2: \R\times S^1\to S^1\) denote the
projection to the second factor. For a point \(a\in S^1\) and a
representative \((T, u)\) of \((T,
\hat{u})\in \hat{\cm}_O\), we define a \(u_a\in L^p_{1, loc}(\R\times
S^1; p_2^*T_f)\) so that:
\begin{itemize}
\item Over \([-T/2, T/2]\times S^1\), \(u_a\) agrees with a
  representative of \((T, \hat{u})\);
\item Over \(\R\times S^1\backslash ([-T/2, T/2]\times S^1)\),
  \(u_a(s, t)=u (a,t)\).
\end{itemize}
Let \(\tau_L\) denote translation by \(L\):
\[\tau_Lw(s):=w(s-L).\]

\begin{definition*}[Chain topology on \(\hat{\cm}_P^+\),
  \(\hat{\cm}_O^+\).]
The chain topology of \(\hat{\cm}_P^+\) and 
  \(\hat{\cm}_O^+\) are respectively the topology with neighborhood
  base given by the subsets
  \(N_P(U_1, \ldots, U_k; \Lambda)\) and \(N_O(U_1, \ldots, U_k;
  \Lambda)\), \(N_O(U)\) defined below:
\smallbreak

\noindent{\bf (a)} Given \(\{\hat{u}_1, \ldots, \hat{u}_k\}\in
\hat{\cm}^+_P(x, y)\), let
\(u_i\) be the centered representative of \(\hat{u}_i\), and
\(U_i\) be a neighborhood of \(u_i\) in \(L^p_{1,loc}(\R\times
S^1; p_2^*T_f)\). Let \(\Lambda\in \R^+\), and \(N_P(U_1, \ldots, U_k; \Lambda)\subset
\hat{\cm}_P^+(x, y)\) be the subset consisting of broken trajectories
\(\{\hat{v}_1, \ldots, \hat{v}_m\}\) such that:
\begin{itemize}
\item \(\exists  L_i\in \R \, \forall i\), and a surjective map \(j: \{1,
  \ldots, k\}\to \{1, \ldots, m\}\), so that \(\tau_{-L_i}v_{j(i)}\in
  U_i\), where \(v_j\) is the centered representative of
  \(\hat{v}_j\);
\item The map \(j\) preserves the partial ordering. Namely,
  \(j(i+1)=j(i)+1\) or \(j(i)\), and in the latter case, \(L_{i+1}-L_i>\Lambda\).
\end{itemize}

\noindent{\bf (b)} When \(\{\hat{u}_1, \ldots, \hat{u}_k\}\) is a
\(k\)-th broken orbit with \(k>0\), and \(U_1, \ldots, U_k, \Lambda\)
be defined as in part (a), the subset \(N_O(U_1, \ldots, U_k; \Lambda)\subset
\hat{\cm}_O^+\) consists of: 
\begin{itemize}
\item broken orbits \(\{\hat{v}_1, \ldots,
\hat{v}_m\}\) with \(m>0\) satisfying similar conditions as in part {\bf (a)}
above, but with \(j\) now mapping from the indexing set of
\(\{\hat{u}_1, \ldots, \hat{u}_k\}\) to the indexing set of \(\{\hat{v}_1, \ldots,
\hat{v}_m\}\), which preserves partial {\em cyclic} ordering, or
\item closed orbits \((T, \hat{v})\) with the property that \(\exists
  L_i\in \R \, \forall i\) and \(a\in S^1\), so that \(\tau_{-L_i}v_{a}\in
  U_i\) \(\forall i\), \(T>k\Lambda\), and \(L_{i+1}-L_i>\Lambda\) for
  \(i=1, \ldots, k-1\).
\end{itemize}
When \((T_u, \hat{u})\) is a closed orbit, and \(U\) is a neighborhood
of \(u_b\) in \(L^p_{1,loc}(\R\times
S^1; p_2^*T_f)\) for some \(b\in S^1\), the subset \(N_O(U)\subset
\hat{\cm}_O^+\) consists of closed orbits \((T, \hat{v})\) such that
\(v_a\in U\) for some \(a\in S^1\).
\end{definition*}

The stronger \({\cal B}\)-topology on \(\cm_P(x, y)\) and \(\cm_O\) in
this Floer theory will be introduced in \S3.2.3, 3.3.1.

\subsubsection{The Conley-Zehnder index}
Let \(A_x\) denote the linearization of \({\cal V}_X\) at \(x\in {\cal
  P}(X)\). The critical point
\(x\) is said to be {\em nondegenerate} if \(A_x\) is surjective. In
this case, any lift of \(x\) in a covering of \(\lc\) is said
to be nondegenerate. 

Recall that the Conley-Zehnder index assigns continuously 
an integer \(\op{CZ}(A)\) to each surjective operator \(A\in \Sigma_C\), where
\(\Sigma_C\) is the space of operators \(A: L^2_1(S^1, \R^{2n})\to
L^2(S^1, \R^{2n})\) of the form \(A=J_0d/dt+\nu(t)\), where \(\R^{2n}\) is equipped with the standard symplectic
structure and complex structure \(J_0\), and \(\nu\) is a
\(C_\epsilon\)-function taking values in the space of self-adjoint matrices. Furthermore, the spectral
flow of a path of such operators \(A(s)\), \(s\in [0,1]\) is the
difference in the Conley-Zehnder indices of \(A(1)\) and
\(A(0)\), and \(\op{CZ}(J_0d/dt)=0.\)
(See e.g. \cite{salamon.park}. Our
\(\op{CZ}=\mu_H=n-\mu_{CZ}\) in \cite{salamon.park}). 

Let \(\hat{\lc}\) denote the universal covering of \(\lc\) and
\(\hat{\cal P}\) be the lift of \({\cal P}\) in \(\hat{\lc}\).
\begin{definition*} 
Fix a unitary trivialization of \(\gamma_0^*K\). 
 Together with a homotopy class of paths \(w\) from \(\gamma_0\) to \(x\),
this trivialization of \(\gamma_0^*K\) induces a (homotopy class of)
unitary trivialization \(\Phi_{x, [w]}: x^*K\to S^1\times \C^n\), for each 
\((x, [w])\in\hat{\cal P}\). The operator 
\({\mathbb A}_{(x, [w])}:=(\Phi_{x, [w]})_* A_x (\Phi_{x,
  [w]})_*^{-1}\) is contained in \(\Sigma_C\), and it is surjective when
\((x, [w])\) is nondegenerate. In this case, define
\[
\hat{\op{ind}}(x, [w])=\op{CZ}({\mathbb A}_{(x, [w])})\in \Z.
\]
\end{definition*}
A well-known index computation shows that if \((x, [w]), (x, [w'])\in
\hat{\cal P}\) are two different lifts of the same \(x\in {\cal P}\),
then
\[
\hat{\ind}(x,[w])-\hat{\ind}(x,[w'])=2c_1^f(\op{im}[w-w']),
\] 
which factors through \(\op{im}\). Thus, \(\hat{\ind}\) descends to
define a \(\Z\)-valued index \(\tilde{\ind}\) for nondegenerate
critical points in \(\tilde{\lc}\), with SF-homomorphism given by \(2\psi_c\).
This in turn defines an absolute \(\Z/2\Z\)-valued grading \(\ind\) for
nondegenerate critical points in \(\lc\), and a relative \(\Z\)-valued
grading on \(\tilde{\cal P}_{ndg}\), which satisfy the requirements in
\S2.1.

\subsubsection{Varying $f$ and Varying $(J, X)$.}
We now describe a well-known and very useful observation. 

We'll write \({\cal X}={\cal X}_f\), \(K=K_f\) to emphasize their dependence on
\(f\). 

First, notice that a diffeotopy from \(f\) to \(f'\) induces a
diffeomorphism from $\ls_fM$ to $\ls_{f'}M$. 
In particular, let \(\Phi=\{\phi_t\, \Big|\, t\in [0,1]\}\) be a symplectic isotopy
connecting the identity map to \(\phi_1\), where \(f'=\phi_1\circ f\),
and let \(\vartheta_t\) be the path
of closed 1-forms so that \(X_{\vartheta_t}\) generates \(\phi_t\). For \(\gamma\in
\Omega_{\gamma_0}(M;f)\), let
\[
\Phi\cdot \gamma(t):=\phi_t(\gamma(t)).
\]
This defines an isomorphism 
\[\Phi: \Omega_{\gamma_0}(M;f)\to \Omega_{\Phi\cdot\gamma_0}(M;f').\]
It also induces an action of \(\pi_1(\op{Symp}_0(M))\) on
\(\pi_0(\ls_fM)\).

Suppose further that the path \(t\mapsto X_{\vartheta_t}\) is in
\(X_f\). (This can always be arranged by e.g. a suitable
reparametrization of the path so that
\(X_{\vartheta_0}=0=X_{\vartheta_1}\)). 
It is easy to see that for any \(X\in {\cal X}_f\), there is an
\(X'\in {\cal X}_{f'}\) such that
\begin{gather*}
\Phi^*{\cal Y}_{X'}={\cal Y}_X\quad \text{and}\\
\phi_t^*\theta_{X',t}:=\theta_{X,t}+\vartheta_t.
\end{gather*}
Furthermore, if \(u(s, t)\) solves (\ref{pcr}), then
\[
w(s,\cdot):=\Phi\cdot u(s,\cdot)
\]
also satisfies (\ref{pcr}), but with \(f\) there replaced by
\(f'\), \(X\) replaced by \(X'\), and \(J_t\) replaced by
\(J'_t\), where
\[
J'_t(\phi_t(x)):=(D\phi_t)J_t(x)(D\phi_t)^{-1}.
\]
(Notice that \(J'\in {\cal J}_{K_f'}\) if \(J\in {\cal
  J}_{K_f}\)). 

This not only defines isomorphisms between the moduli spaces
\({\cal P}^{\gamma_0}(M, f; X)\), \(\cm_P^{\gamma_0}(M, f; J, X)\),
\(\cm_O^{\gamma_0}(M, f; J, X)\) and
\({\cal P}^{\Phi\cdot\gamma_0}(M, f'; X')\),
\(\cm_P^{\Phi\cdot\gamma_0}(M, f'; J', X')\),
\(\cm_O^{\Phi\cdot\gamma_0}(M, f'; J', X')\) respectively, but also equivalences of the
relevant deformation operators (\(A_x\), and \(E_u, \tilde{D}_u\) in
\S3.2, 3.3) by
similarity transformations. Thus the Floer theories associated to
\((M, f, \gamma_0;J, X)\) and \((M, f', \Phi\cdot\gamma_0;
J', X')\) are completely equivalent.

Because of this equivalence, in this paper we shall fix the
symplectomorphism \(f\) and vary the almost complex structure \(J\) and the symplectic vector field \(X\).

\subsection{Structure of the Moduli Spaces: ${\cal P}, \cm_P$.}
The goal of the rest of this section is to show that for regular \((J, X)\)
(to be defined in \S3.3), the associated formal flow is a
Floer system.

Since the results in this subsection all follow from simple adaptations of
the literature, we shall omit most details except the notions needed
for later sections. 

For the rest of this section, we fix a \(J\in {\cal J}_K^{reg}\) and
perturb \(X\in {\cal X}\) by hamiltonian vector fields to achieve
transversality.

\subsubsection{The space \({\cal P}\).}

Recall the definition of \(M_c(J_t), M_c(J)\) from the end of \S3.1.2.
\begin{definition*}
Given \(J\in {\cal J}_K^{reg}\), we say that \(X\in {\cal X}\) is
{\em \(J\)-nondegenerate} if the following hold:

(a) ${\cal P}(X)$ consists of finitely many nondegenerate points. 

(b) For any $\gamma \in{\cal P}(X)$, $\gamma(t)\not\in M_1(J_t)\,
\forall t$.
\end{definition*}
\begin{proposition*}
Given \(J\in {\cal J}_K^{reg}\) and \(X\in {\cal X}\), there is a Baire set
\({\cal H}^{ndg}(J, X)\subset{\cal H}\) such that 
\(X+\chi_H\) is \(J\)-nondegenerate for all \(H\in {\cal H}^{ndg}(J, X)\).
\end{proposition*}
\begin{proof}
These follow from simple adaptations of \cite{HS} Theorem
3.1, using Lemma \S3.1.2(a) above. For part (b), 
we replace the evaluation map in \cite{HS} by
\[
M_1(J)\times{\cal P}^{\cal H}(X)\to T_f\times_{S^1} T_f: (q, (x, H))\mapsto (q, x(\pi_S(q))),
\]
where \({\cal P}^{\cal H}(X)=\bigcup_{H\in {\cal H}}{\cal
  P}(X+\chi_H)\) is the universal moduli space of critical points that
fibers over \({\cal H}\), and \((x, H)\) denotes an element in this
universal moduli space: \(x\in {\cal P}(X+\chi_H)\). 
\end{proof}

In particular, this means that given a fixed pair \((J, X)\in {\cal J}_K^{reg}\times {\cal
  X}\), the pair \((J, X+\chi_H)\) satisfies (FS1) for
all \(H\in {\cal H}^{ndg}(J, X)\). 
\subsubsection{The Space $\cm_{P}$.}\label{moduli}

We now show that (FS2) holds for ``generic''  \((J,
X)\), in a sense we clarify next.

Consider Hamiltonian perturbations in the following space:
\begin{definition*}
Let \(J\in {\cal J}_K^{reg}\) and let \(X\in {\cal X}\) be \(J\)-nondegenerate.
Let \(\delta\in \R^+\) and \(k\in \Z^+\). Then 
\begin{equation}\label{pathPert}
V_\delta^k(J, X):=\Big\{H \, \Big| \, H\in {\cal H}, \, \|H\|_{C_{\epsilon}}\leq \delta; \nabla^i H_t(x(t))=0 \,
\forall x\in {\cal P} (X), \, i=0, 1, \ldots, k\Big\}.
\end{equation}
\end{definition*}
Notice that \({\cal P}(X+\chi_H)={\cal P}(X)\) \(\forall H\in
V_\delta^k(J, X)\).

Given \(J\in {\cal J}_K^{reg}\), an element \(u\) of \(\cm_P\) or \(\cm_O\) is said to be
{\em $J$-regular} if in addition to being nondegenerate, it
also satisfies: 
\[
u(\cdot, t)\, \cap M_0(J_t)=\emptyset \quad \forall\, t.
\]
A moduli space (any variant of \(\cm_P\) or \(\cm_O\)) 
is said to be {\em nondegenerate} if it consists of
nondegenerate elements; it is {\em $J$-regular} if it consists
of $J$-regular elements.

\begin{proposition*}
Let \(J, X\) be as in the previous Definition. Then for any \(k\geq2\)
and any small positive number \(\delta\), there 
is a Baire set \(V_{\delta}^{k, Preg}(J, X)\subset V_\delta^k(J, X)\)  
such that all \(H\in V_{\delta}^{k, Preg}(J, X)\) satisfies: 

(a) (FS2) holds for \((J, X+\chi_H)\);

(b) \(\cm_P^i(x, y;J,X+\chi_H)\) is \(J\)-regular for any \(i\leq2\),
\(x, y\in {\cal P}\).
\end{proposition*}
\begin{proof} 
These again follow from simple adaptations of \cite{HS}. 
Transversality follows the arguments of \cite{HS} Theorem
3.2, using again Lemma 3.1.2, and replacing 
the evaluation map by
\[
M_0(J)\times\cm_P^{i, V^k_\delta}(x, y);J, X)\to T_f\times_{S^1} T_f: (q, (u, H))\mapsto (q, u(0,\pi_S(q))),
\]
where \(\cm_P^{i, V^k_\delta}(x, y;J, X)=\bigcup_{H\in
  V^k_\delta(J,X)}\cm_P^i(x, y;J, X+\chi_H)\) is the
universal moduli space which fibers over \(V^k_\delta(J,X)\), the term
\((u, H)\) denotes an element over \(H\in V^k_\delta(J, X)\) in
\(\cm_P^{i, V^k_\delta}(x, y;J, X)\), and \(u\in \cm_P(x, y;J, X+\chi_H)\).

The compactness follows the argument of \cite{HS} Theorem 3.3. 
The main ingredients are: 
Gromov compactness, \(J\)-nondegeneracy of \(X+\chi_H\), and statement (b) in the
Proposition proven by the transversality argument above. 
The well-known fact that \(u\) decays exponentially near the
nondegenerate critical points \(x\), \(y\) (see e.g. \cite{floer.cmp} p.607, 
\cite{floer.cpam} pp-801-803) is also used. 
\end{proof}

\subsubsection{The Configuration Space \({\cal B}_P\), and the
  Deformation Operator \(E_u\).}

We now introduce the {\em configuration space} \({\cal B}_P\), which
endows \(\cm_P\) with the \({\cal B}\)-topology, and the relevant
deformation operator. These notions are used 
in the omitted details of the proof of Proposition
\ref{moduli}, and shall also be needed in later sections.

\label{def:beta}
Let \(\Theta:=\R\times S^1_1\); and recall that \(p_2: \Theta\to S^1_1\) is the
projection to the second factor. 

Let \(\beta:\R\to [0,1]\) be a
smooth cutoff function supported on \(\R^+\) such that \(\beta(s)=1\)
as \(s\geq 1\).

Given \(J\in {\cal J}_K\), we define a
{\it $t$-dependent} exponential map \(\exp: T\lc\to \lc\),
\[
\exp (x, \xi) (t):=
(\exp_{x(t)}^{g_t} \xi (t),t),\] where \(\exp_{x(t)}^{g_t}\) is the
  exponential map in the fiber $M$ with respect to $g_t$, the metric
  on \(M\) corresponding to \(J_t\). \(\xi\in L^p_1(x^*K)=T_x\lc\).

\begin{definition*}\label{B-P}
Suppose \(V\) is an euclidean or hermitian vector bundle over \(\Theta \);
\(\sigma_-,\sigma_+\in \R\). Let \(L^p_{k:(\sigma_-,\sigma_+)}(V)\) be the (exponentially)
weighted Sobolev space:
\[
\{\xi\, |\, \varsigma^{\sigma_-\sigma_+}\xi\in L^p_k\}, \,
\text{where \(\varsigma^{\sigma_-\sigma_+}(s):=e^{-(s\beta(s)\sigma_++s\beta(-s)\sigma_-)}\).}
\]
Let \(p_2^*T_f\to \Theta\) be the pullback bundle of $T_f\to S^1$. Given two 
critical points \(x, y\in {\cal P}(X)\), the {\em configuration space} 
is
\begin{eqnarray*}
\lefteqn{{\cal
  B}^{(\sigma_1,\sigma_2)}_P(x,y):=}\\
&&\Big\{ u|\, u\in L^p_{1,loc}(\Theta, p_2^*T_f), \\
&& \quad u(s,\cdot)=\exp (y, \xi_+(s, \cdot))\, \mbox{for some} \, \xi_+\in
L^p_{1:(0,\sigma_2)}(p_2^*(y^*K))\, \mbox{when} \, s>\rho_+(u);\\
&& \quad u(s,\cdot)=\exp (x, \xi_-(s, \cdot)) \, \mbox{for some} \, \xi_-\in
L^p_{1:(\sigma_1,0)}(p_2^*(x^*K)) \, \mbox{when} \, s<\rho_-(u)\Big\},
\end{eqnarray*}
where \(\rho_\pm\in \R\) are numbers depending on \(u\). The integer \(p>2\).
\end{definition*}
Understood as the space of paths in \(\lc\) between the two points
\(x\), \(y\), \({\cal B}_P(x, y)\) decomposes into many components
according to the homology classes of the paths, like \(\cm_P(x, y)\).
Thus, we may define \({\cal B}_P((x, [w]), (y, [v]))\), \({\cal
  B}_P^k(x, y)\) in complete analogy to \(\cm_P((x, [w]), (y, [v]))\),
\(\cm_P^k(x, y)\).

The space ${\cal  B}^{(\sigma_1,\sigma_2)}_P(x,y)$ is a Banach manifold (see
e.g. the argument of \cite{ floer.cpam} Theorem 3).
The local
model of a neighborhood of \(u\) is $L^p_{1:(\sigma_1,\sigma_2)}(u^*K)$ via
the $t$-dependent exponential map 
\[
\exp (u, \xi):=
(\exp_{u(s,t)}^{g_t} (\xi (s,t)),(s,t)).\]
The precise definition of weighted moduli space is:
\[
\cm_P^{(\sigma_1,\sigma_2)}(x,y; J,X):=\cm_P(J,X)\cap{\cal B}_P^{(\sigma_1,\sigma_2)}(x,y).
\]
Notice that we do not require any nondegeneracy condition on \(x, y\)
for the definitions of \({\cal B}^{(\sigma_1,\sigma_2)}_P(x, y)\) and
\(\cm^{(\sigma_1,\sigma_2)}_P(x, y; J,X)\).
However, with suitable nondegeneracy conditions on  
\(x, y\), \(\cm_P^{(\sigma_1,\sigma_2)}(x,
y; J, X)\) admits a description as the zero locus of the Fredholm section
\(\bar{\partial}_{JX}\) of a Banach
bundle. (Cf. e.g. \cite{floer.cpam}). See the next Lemma. 

Let \(E_u^{(\sigma_1,\sigma_2)}: L^p_{1:(\sigma_1, \sigma_2)}(u^*K)\to
L^p_{:(\sigma_1, \sigma_2)}(u^*K)\) denote the linearization of \(\bar{\partial}_{JX}\) at \(u\).

In this Floer theory, the weighted versions of gradings are defined
via \(\tilde{\ind}^\sigma\):
\begin{gather*}
\ind^\sigma(x)=\tilde{\ind}^\sigma(x, [w])\quad \mod 2;\\
\op{gr}^{(\sigma_1, \sigma_2)}((x_1, [w_1]), (x_2, [w_2]))=\tilde{\ind}^{\sigma_1}(x_1, [w_1])-\tilde{\ind}^{\sigma_2}(x_2, [w_2]),
\end{gather*}
where \(\tilde{\ind}^\sigma\) is the {\em generalized Conley-Zehnder
  index}, defined as follows.
Given \(\sigma\in \R\), we say that an \(x\in {\cal P}(X)\) is {\em \(\sigma\)-weighted
  nondegenerate} if \(A_x+\sigma\) is surjective. In this case, 
\[\tilde{\ind}^{\sigma}(x, [w]):=\tilde{\ind}^{\sigma}_{\gamma_0}(x, [w]):=\op{CZ}({\mathbb A}_{(x,[w])}+\sigma).\] 

A routine modification of the literature (e.g. \cite{floer.cmp},
\cite{salamon.park}, \cite{salamon-zehnder}) yields the next lemma.
\begin{lemma*}
Let $u\in \cm_P ((x,[z]), (y,[w]); J,X)$, and suppose that $x, y\in
{\cal P} (X)$ are \(\sigma_1\)-weighted nondegenerate and
\(\sigma_2\)-weighted nondegenerate respectively.
Then the operator
\[E_u^{(\sigma_1, \sigma_2)}: L^p_{1:(\sigma_1, \sigma_2)}(u^*K)\to
L^p_{:(\sigma_1, \sigma_2)}(u^*K)\] 
is Fredholm of index
\(\ind^{\sigma_1}(x,[w])-\ind^{\sigma_2}(y,[v])\). 
\end{lemma*}
In fact, the moduli space
\(\cm_P^{(\sigma_1, \sigma_2)}(x,
y)\) is independent of small perturbations to the weights \(\sigma_1, \sigma_2\)
When \(x\), \(y\) are \(\sigma_1\)-weighted nondegenerate and
\(\sigma_2\)-weighted nondegenerate respectively.
\begin{notation*}
By standard decay estimates, \(\cm_P(x, y)=\cm_P^{(0, 0)}(x,
y)\). Thus, we shall omit the superscripts \((\sigma_1, \sigma_2)\) when the
weights \((\sigma_1, \sigma_2)=(0, 0)\), and the critical points \(x, y\) are
nondegenerate. 
\end{notation*}

\subsection{Structure of the Moduli Spaces: $\cm_{O}$.}

We now verify the genericity of the condition (FS3); more precisely, we
prove:
\begin{proposition*}\label{thm:moduliO}
Given \(J\in {\cal J}_K^{reg}\), a \(J\)-nondegenerate \(X\in {\cal
  X}\), an integer \(k\geq2\), and a small \(\delta\in \R^+\), the set 
\[
V^{k, reg}_\delta(J, X):=\Big\{H\, \Big|\, H\in V^k_\delta(J, X), \,
\text{(FS2), (FS3) hold for \((J, X+\chi_H)\)}\Big\}
\]
is Baire in \(V^k_\delta (J, X)\). 
\end{proposition*}
Notice that by assumption, (FS1) already holds for
\((J, X+\chi_H)\), for all \(H\in V^k_\delta(J, X)\). We call \((J,
X)\) a {\em regular pair} if \(J\in {\cal J}_K^{reg}\),
\(X=X_0+\chi_H\), where \(X_0\) is \(J\)-nondegenerate and \(H\in
V^{k, reg}_\delta(J, X_0)\). Since the properties (FS1), (FS2), (FS3)
hold for regular pairs, \(\tau_F(J, X)\), \(\zeta_F(J, X)\), \(I_F(J,
X)\) are well-defined for regular pairs. Regular pairs are generic in
the sense described in this section.

This proposition follows from the combination of 
Lemmas \ref{cm-O-s}, \ref{lem:mo2}, \ref{lem:mo3} below: \(V^{k,
  reg}_\delta(J, X)\) is Baire since it contains \(\bigcap_{\Re\in \R^+}V^{k, \Re-reg}_\delta(J,
X)\), where \(V^{k, \Re-reg}_\delta(J,X)\) are the open dense sets in Lemma \ref{lem:mo3}. 
We first introduce some basic notions before stating these lemmas.  

\subsubsection{The Configuration Space \({\cal B}_O\) and the
  Deformation operator \(\tilde{D}_{(T,u)}\). }

Because we include closed orbits of all periods in \(\cm_O\), we shall
often write an element of \(\cm_O\) as a pair \((T, u)\), where
\(T\in \R^+\) is the period of the closed orbit, and \(u\) is a
section of the \(M\)-bundle \(p_2^*T_f\) over \(S_T^1\times S^1_1\), where \(p_2:S^1_T\times S^1_1\to S^1_1\) is again the projection to
the second factor.

The configuration space \({\cal B}_O\) for \(\cm_O\) is then
\[
{\cal B}_O:=\Big\{(T,u)\, \Big|\, T\in \R^+; u\in L^p_1(S^1_T\times S^1_1,
p_2^*T_f)\quad \text{s.t. \(u\Big|_{\{0\}\times S^1_1}\) is homotopic to
    \(\gamma_0\)}\Big\},
\]
It is easy to see that \({\cal B}_O\) is a Banach manifold modelled on \(\R\times
L^p_1(u^*K)\), and it fibers over \(\R^+\) by mapping each
\((T, u)\in{\cal B}_O\) to \(T\).
\(\cm_O(J, X)\) embeds in \({\cal B}_O\) as the zero locus of the
Fredholm section \(\bar{\partial}_{JX}\) of a Banach bundle. 

Let \(\tilde{D}_{(T, u)}\) be the linearization of \(\bar{\p}_{JX}\) at
\((T, u)\in \cm_O(J, X)\). It has the following expression:
\[
\tilde{D}_{(T, u)} (\varrho, \xi)=D_u\xi-\varrho /T \partial_su \quad
\text{for \((\varrho, \xi)\in \R\times L^p_1(u^*K)\)},
\]
where the operator \(D_u: L^p_1(u^*K)\to L^p(u^*K)\) is given by the same
formula as \(E_u\), but here \(s\) takes value in \(S^1_T\)
instead of \(\R\). 

\(\tilde{D}_{(T, u)}\) fits into the deformation complex:
\[
\begin{CD}
\R @>d_u>> \R\oplus L^p_1(u^*K) @>\tilde{D}_{(T, u)}>> L^p(u^*K)
\end{CD}
\]
where \(d_u\) is the linearization of the \(S^1\) action (by translation
in \(s\)) on \(\cm_O\),
\[
d_u(\lambda) = (0,\lambda\partial_su).
\]
It is useful to combine \(d_u\) and \(\tilde{D}_{(T,u)}\) into the
following operator \({\cal D}_{(T,u)}\), which is deformation
operator of the reduced moduli space
\(
\hat{\cm}_O=\cm_O/S^1:
\)
\[
\begin{split}
{\cal D}_{(T, u)}={\cal D}_{(T, u)}^{J, X}: \R\oplus L^p_1(u^*K)&\longrightarrow \R\oplus L^p(u^*K),\\
{\cal D}_{(T, u)}=\tilde{D}_{(T, u)}+d^*_u &=
\left(\begin{array}{cc}
0&\Pi_{\partial _s u} \\
-1/T\partial_s u & D_u\end{array}\right),
\end{split}
\]
where \(\Pi_{\partial_s u}\) denotes \(L^2\)-orthogonal projection to
\(\partial _s u\), and \(d_u^*\) is the formal \(L^2\)-adjoint of
\(d_u\).

The superscripts \(J, X\) are added to the operators when we wish to
emphasize their dependence on \(J, X\). 

\({\cal D}_{(T, u)}\) is obviously Fredholm, since \(D_u\) is
elliptic. 
By Riemann-Roch
\[
\ind({\cal D}_{(T, u)})=\ind(D_u)=2\psi_c(A)\quad \text{when $(T, u)\in \cm_O(A)$.}
\]

Notice that the nondegeneracy of a \((T, u)\in \cm_O\) is equivalent
to the nondegeneracy of its image in \(\hat{\cm}_O\): 
Since a closed orbit is by definition nonconstant in
\(s\), \(d_u\) is always injective, and thus the surjectivity of 
\(\tilde{D}_{(T, u)}\) and \({\cal D}_{(T, u)}\) are equivalent.

The standard transversality argument works (only) for the space of {\it simple}
closed orbits, which we denote by $\cm_{O}^{d,sim}\subset\cm_O^d$.
\begin{*lemma}\label{cm-O-s}
Let \((J, X)\) be as in the Proposition. Then the subset 
\[
\Big\{H\, \Big|\, H\in V^k_\delta(J, X), \, \text{\(\cm_O^{sim}(J, X+\chi_H)\) is nondegenerate}\Big\}
\]
is Baire in \(V^k_\delta(J, X)\). So is the subset
\begin{equation}\label{j-ndg-set}
\Big\{H\, \Big|\, H\in V^k_\delta(J, X), \, \text{\(\coprod_{l\leq 1}\cm_O^{l, sim}(J, X+\chi_H)\) is $J$-regular}\Big\}.
\end{equation}
\end{*lemma} 
\begin{proof}
Let \(\cm_O^{V^k_\delta, sim}(J, X)=\bigcup_{H\in V^k_\delta(J, X)}
\cm_O^{sim}(J, X+\chi_H)\) be the universal moduli space
which fibers over \(V^k_\delta(J, X)\). We want to show that it is
nondegenerate, namely, 
for any \(((T, u), H)\in \cm_O^{V^k_\delta, sim}(J, X)\), the operator
\[\begin{split}
\tilde{D}^{V^k_\delta}_{((T, u), H)}: \R\oplus L^p_1(u^*K) \oplus
V^k_\delta(J, X)\to L^p(u^*K),\\
 \tilde{D}^{V^k_\delta}_{((T, u),
  H)}((\varrho, \xi), h)=\tilde{D}_{(T, u)}^{J, X+\chi_H}(\varrho, \xi)+\nabla h
\end{split}\]
is surjective.
Equivalently, letting \(D^*_u\) denote the formal \(L^2\)-adjoint of
\(D_u\), we want to show that there is no \(\eta\) such that
\begin{equation}
\begin{split}
D_{u}^*\eta &=0,\\
\langle\partial_s u, \eta\rangle_2& =0,\label{g(t)}\\
\langle\nabla h, \eta\rangle_2&=0 \quad \forall  h.
\end{split}
\end{equation}
To apply the standard transversality argument
\footnote{There is vast literature on the transversality proof; 
unfortunately, many of the papers contain (minor) errors. \cite{trans} is
the one in which the author is unaware of errors, and it shall be
our main reference.}, 
we need the following
observations: first, by unique continuation (e.g. \cite{trans} Proposition
3.1) the set 
\[
\Omega(u):=\Big\{(s,t)\Big|\, u(s,t)=x(t)\, \mbox{for some}\, x\in
{\cal P}(X)\Big\}\subset S^1_T\times S^1_1
\] 
consists of discrete points for any \(u\in \cm_O(J, X+\chi_H)\); also, for simple \(u\) the
argument of \cite{trans} Proposition 4.1 shows that the set
\[
R(u)=\Big\{(s, t)\, \Big|\, (s, t)\not\in \Omega(u), \, u(s, t)\neq u(s',
t) \, \forall t, s\neq s'\Big\}\subset S^1_T\times S^1_1
\]
is open dense.
One may then 
follow the standard arguments (e.g. in the proof of
Theorem 5.1 (i) of \cite{trans}) 
to show that \(\eta=g(t)\partial_su\) for some real-valued function \(g(t)\). 
However, from (\ref{g(t)}), we have 
\[
\int g(t)|\partial _su(s,t)|^2\, ds\, dt\, =0.
\] 
This contradicts with the fact that $g(t)\neq 0 \,\, \forall t$ and
$\partial _su(s,t)\neq 0$ for all but discrete $s, t$
(cf. \cite{trans}).
To show that the set (\ref{j-ndg-set}) is also Baire, consider in
addition the evaluation map
\[
ev: M_0(J)\times\cm_O^{V^k_\delta, l,sim}(J, X)\to T_f\times_{S^1}T_f, \quad (y, ((
T, u), H))\mapsto (y, u(0, \pi_S(y))).
\]
This can be seen to be transverse to the
diagonal by the usual argument; so a generic fiber of $ev^{-1}(diag)$
is a manifold of dimension $2n-3+1
+l-2n=l-2<0$---hence empty. 
\end{proof}

Let $\cm_{O}^{1,sim, \Re}\subset \cm_O^{1,sim}$ be the subspace of
simple closed orbits with energy $\leq \Re$. 
\begin{*lemma}\label{lem:mo2}
Let \(J\in {\cal J}_K^{reg}\), \(X\) be $J$-nondegenerate, and \(\Re\in \R^+\) be
arbitrary. 
If \(\cm^{\leq1, \Re}(J, X):=\coprod_{l\leq1}(\cm_O^{l, \Re}(J, X)\,\cup \, \cm_P^{l,
  \Re}(J, X))\) is $J$-regular, then the reduced moduli space 
\(\hat{\cm}^{0, \Re}(J, X):=\hat{\cm}_O^{0, \Re}(J, X)\,\cup \, \hat{\cm}_P^{0, \Re}(J, X)\) is compact. 
\end{*lemma}
\begin{proof}
The compactness of \(\hat{\cm}_P^{0, \Re}(J, X)\) was a part of the
proof of Proposition \ref{moduli}, we therefore only need to consider
\(\hat{\cm}_O^{0, \Re}(J, X)\).

Note that there is a map $\hat{\cm}_{O}^{0, \Re}\to \hat{\cm}_{
O}^{0,sim, \Re}$ by mapping a multiple-cover to its underlying simple closed
orbit. The fiber of a \(u\in  \hat{\cm}_{
O}^{0,sim, \Re}\) consists of \(n(u)\) points, where
\begin{equation}\label{eq:n-u}
n(u):=[-\Re/{\cal Y}_X([u])]\in\Z^+
\end{equation}
is bounded above since \(-{\cal Y}_X([u])\) is bounded below by a
positive number. To see this, suppose there is a sequence of closed
orbits \(\{u_n\}\) with \({\cal E}(u_n)\to 0\) as \(n\to \infty\). By
Gromov compactness, there is a subsequence converging in \(C_0\)-norm
to a constant flow (i.e. a critical point). This means that the
homology class \([u_n]=0\) for large \(n\), contradicting the fact
that closed orbits have positive energy.

Thus, it suffices to show that \(\hat{\cm}^{0, sim,
  \Re}_O\) is compact. 
Gromov compactness assures us that there are only three ways
in which compactness can fail. One is bubbling off spheres; the second is
the existence of a sequence of elements in \(\hat{\cm}_O^{0,sim,\Re}(J,X)\) with
periods going to \(\infty\), which converge weakly to a broken orbit; the third is the existence of a
sequence \(\{(T_i, u_i)\}\subset \hat{\cm}_O^{0,sim,\Re}(J,X)\)
converging to a multiple-cover. 

The first possibility is eliminated by the assumption of \(J\)-regularity.
The second possibility is a codimension 1 phenomenon, and is
eliminated by the assumption that \(\cm^{\leq1, \Re}_P(J, X)\) is
nondegenerate. In the third case, the multiple cover that the sequence
converges to cannot be nondegenerate, contradicting the assumption
that \(\cm^{\leq1, \Re}_O(J, X)\) is nondegenerate. 
\end{proof}

It remains to show that multiple covers are also generically
nondegenerate. For this purpose, notice that to show that
\(\hat{\cm}_O^\Re\) is nondegenerate, it is equivalent to show that each
\(u\in \hat{\cm}_O^{sim, \Re}\) is ``\(n(u)\)-nondegenerate'', \(n(u)\)
being given by (\ref{eq:n-u}):

Let \(u\in \hat{\cm}^{sim}_O\); we denote the unique \(\Z/m\Z\)-cyclic
cover of \(u\) by \(u^m\). \(u\) is said to be {\em \(n\)-nondegenerate} if \(u^1,
u^2, \ldots , u^n\) are all nondegenerate.

\begin{*lemma}\label{lem:mo3}
Let \(J, X, k, \delta\) be as in the Proposition. Then
for each \(\Re\in \R^+\), the set
\[
V^{k, \Re-reg}_\delta(J, X):=\Big\{\, H\, \Big|\, H\in V^k_\delta(J, X), \,
\text{\(\cm^{\leq1, \Re}(J, X+\chi_H)\) is \(J\)-regular}\Big\}
\]
is open and dense in \(V^k_\delta(J, X)\).
\end{*lemma}
\begin{proof}
Because by the previous lemma \(\cm^{\leq1, \Re}\) consists of
finitely many isolated points, the proof is reduced to showing that
transversality can be achieved by Hamiltonian perturbations near each
\(u\in \cm^{\leq , \Re}\). 

The openness of \(V^{k, \Re-reg}_\delta(J, X)\) is obvious from the
compactness result of Lemma \ref{lem:mo2}. 
To show the denseness, we shall show that for any \(H\in V^k_\delta(J,
X)\), we may perturb \(H\) by some small \(h\) so that \(\cm^{\leq 1,
  \Re}(J, X+\chi_{H+h})\) is \(J\)-regular. By Proposition
\ref{moduli} and Lemma \ref{cm-O-s}, we may assume without loss of
generality that \(\cm_P^{\leq1, \Re}(J, X+\chi_H)\) and
\(\cm_O^{\leq1, sim, \Re}(J, X+\chi_H)\) are already
\(J\)-regular. Adding any sufficiently small \(h\) will not affect the
\(J\)-regularity of these moduli spaces. We want to show that \(h\)
may be chosen to make each $(T, u)\in
\hat{\cm}_{O}^{1,sim,\Re}$ \(n(u)\)-nondegenerate.

For this purpose, we examine how the
cokernel of the deformation operator changes under perturbation,
similar to Lemma 5.13 of  \cite{T2}.

Let \(m\leq n(u)\) be a positive integer. Recall that the domain of
\(u^m\), denoted \(\Theta_{u^m}\), is a \(\Z/m\Z\)-cyclic covering of
the domain of \(u\), \(\Theta_u\):
\[
c_m: \Theta_{u^m}\to \Theta_u.
\]
The group \(\Z/m\Z\) acts on the bundle \((u^m)^*K\) over
\(\Theta_{u^m}\) by deck transformation. On the other hand, let  
\(\R^m\) denote the Euclidean space spanned by the orthonormal basis
\(\{e_i|i\in\Z/m\Z\}\), endowed with a \(\Z/m\Z\) action by
cyclic permutation of the basis vectors: let \(L\) be a
generator of \(\Z/m\Z\), then \(Le_i=e_{i+1}\). 
Consider 
\[
V_m:=(u^m)^*K\otimes_{\Z/m\Z}\R^m. 
\]
This is an \(\R^{2nm}\)-bundle over \(\Theta_u\), and there is 
a standard isomorphism
\[
i_m: \Gamma((u^m)^*K) \to \Gamma(V_m),
\]
via the identification of the sheaf of sections of \(V_m\)
with the direct-image under \(c_{m*}\) of the sheaf of sections of \((u^m)^*K\).

Instead of \({\cal D}_{((mT), u^m)}\), we may equivalently consider the operator
\(^m{\cal D}_{(T, u)}\) induced by this isomorphism.
\[
\begin{split}
^m{\cal D}_{(T, u)}: \R\oplus L^p_1(V_m) &
\to \R\oplus L^p(V_m), \\
^m{\cal D}_{(T, u)} &=\left(\begin{array}{cc}
0& \Pi_{i_m(\partial_s u^m)} \\
-1/(mT) i_m(\partial_s u^m) &^mD_{(T, u)} \\
\end{array}\right).
\end{split}
\]
  
Suppose $^m{\cal D}_{(T, u)}$ has a $p_m$-dimensional cokernel, spanned by
$\underline{\eta}_1, \underline{\eta}_2, \ldots ,
\underline{\eta}_{p_m}$. Then \(\ker \,^m{\cal D}_{(T, u)}\) is also
$p_m$-dimensional, spanned by $\underline{\xi}_1, \underline{\xi}_2,
\ldots, \underline{\xi}_{p_m}$. Note that \(\Z/m\Z\) acts orthogonally on the kernel
and the cokernel; in fact, the action on the kernel and the cokernel
are the same, since the kernel and cokernel both vanish under some
deformation of the operator. Furthermore, since by assumption \((T, u)\) is 1-nondegenerate, there is no $\Z/m\Z$-invariant elements in \(\ker
\,^m{\cal D}_{(T, u)}, \, \cok \,^m{\cal D}_{(T, u)}\). Therefore
the bases $\{\underline{\xi}_i\}$, $\{\underline{\eta}_j\}$ may be chosen such that with respect to them, the
representations of \(\Z/m\Z\) on \(\ker \,^m{\cal D}_{(T, u)}, \, \cok \,^m{\cal D}_{(T, u)}\)
both have the form:
\begin{equation}\label{repr}
L^k:=\bigoplus_{l=1}^{p'} \left(\begin{array}{cc}\cos (2\pi q_l k/m) &-\sin (2\pi q_l
    k/m)\\
\sin (2\pi q_l k/m) &\cos (2\pi q_l k/m)\end{array}\right)
\bigoplus (-1)^{\oplus (p_m-2p')}.
\end{equation}
The \((-1)\) components appear only when \(m\) is even.
We choose the orientation of the bases such
that $q_l\in \{1,2, \ldots, [(m-1)/2]\}$.

Let \((\ker \,^mD_{(T, u)})^{\perp}, (\cok \,^mD_{(T, u)})^{\perp}\)
denote the \(L^2\)-orthogonal complement of the real line spanned by 
\(i_m(\partial _s u^m)\) in \(\ker \,^mD_{(T, u)}, \cok \,^mD_{(T, u)}\) respectively. Then
\[
(\ker \,^mD_{(T, u)})^{\perp}\subset \ker \,^m{\cal D}_{(T, u)}, \,(\cok \,^mD_{(T, u)})^{\perp}\subset\cok \,^m{\cal D}_{(T, u)}
\]
and their $L^2$-orthogonal complements are at most 1-dimensional and
are invariant under the $\Z/m\Z$ action. Say $\gamma$ is a nontrivial
element in either of these 1-dimensional spaces, then its
$\R$-component must be nontrivial. Observe however that the
$\R$-component is invariant under the $\Z/m\Z$ action, so in fact
$\Z/m\Z$ acts trivially on these orthogonal complements; they are
therefore empty by our 1-nondegeneracy assumption. Hence 
\[
(\ker \,^mD_{(T, u)})^{\perp}= \ker \,^m{\cal D}_{(T, u)}, \,(\cok \,^mD_{(T, u)})^{\perp}=\cok \,^m{\cal D}_{(T, u)},
\]
and we shall regard $\underline{\xi}_i, \underline{\eta}_j$ as
elements in the kernel/cokernel of $^mD_{(T, u)}$ respectively.

\begin{claim*}
Let \(J, X\) be as in the Lemma and let \(H\in V^k_\delta(J, X)\) be
such that \(\cm_P^{\leq1, \Re}(J, X+\chi_H)\) and
\(\cm_O^{\leq1, sim, \Re}(J, X+\chi_H)\) are \(J\)-regular. Let \((T, u)\in
\cm_O^{1,sim,\Re}(J,X+\chi_H)\), and \(m\in\{2,3, \ldots, n(u)\}\). 
Then there exists an arbitrarily small
\(h\in V^k_\delta(J, X)\) with the following properties:

(i) \(\nabla h_t(u(\cdot, t))=0\) \(\forall t\);

(ii) \(h\) is supported on a small tubular neighborhood of the image
of \(u\) in \(T_f\). 

(iii) \(\dim \cok
\,^mD^{J, X+\chi_{H+h}}_{(T, u)}<\dim\cok \,^mD^{J, X+\chi_H}_{(T, u)}\) if \(\dim\cok
\,^mD^{J, X+\chi_H}_{(T, u)}\neq0\). 
\end{claim*}
To finish the proof of the Lemma, apply the above Claim iteratively
for each \(u\) and \(m\) until the dimension of all
the relevant cokernels become 0, requiring that the perturbation
\(h_i\) at step \(i\) to be of \(C_{\epsilon}\)-norm smaller than
\(2^{-i}\varepsilon\), while also being small enough to preserve the
nondegeneracy conditions already obtained. 
This can be done in finite steps since \(\# (\hat{\cm}_O^{0, sim,
  \Re}(J, X+\chi_H))\), \(n(u)\) and \(p_m\) are all finite. Thus, 
we eventually obtain an \(h=\sum_i h_i\) with
\(\|h\|_{C_{\epsilon}}<\varepsilon\) satisfying the desired properties. 
\end{proof}

\medbreak
\noindent{\it Proof of the Claim.}
Suppose \(\dim\cok\,^mD^{J, X+\chi_H}_{(T, u)}\neq0\) and is nondecreasing under the
perturbation by $h$. That means that for each
$\underline{\eta}_i$, \(i=1, \ldots, p_m\), there is a small $\underline{\eta}_0
\in L^p_1(V_m)$ such that 
$$
(\, ^mD_{(T, u)}^{J, X+\chi_{H+h}})^* (\underline{\eta}_i +\underline{\eta}_0)
=(\, ^mD_{(T, u)}^{J, X+\chi_H})^{*}(\underline{\eta}_0)
+\nabla \nabla h (u) (\underline{\eta}_i+\underline{\eta}_0)=0. 
$$
Of course, this is possible only when 
$$
\Pi_{(\ker \, ^mD^{J, X+\chi_H}_{(T, u)})^{\perp}} (\nabla \nabla h (u) (\underline{\eta}_i+\underline{\eta}_0))=0.
$$
When $h$ is small, the set of such $h$'s is modeled on the solution
space of the equations:
\begin{equation}\label{contrad}
a_{ij}:=\langle\underline{\xi}_j, \nabla\nabla h (u)
\underline{\eta}_i\rangle_2=0\quad \forall i, j \in \{1, 2, \ldots, p_m\}.
\end{equation}
We shall show that there exists some choice
of $h$ such that the matrix \((a_{ij})\)
is nontrivial, contradicting (\ref{contrad}).

Since we shall choose \(h\) to be supported on a small neighborhood,
it suffices to consider restrictions of \(\underline{\eta}_j\),
\(\underline{\xi}_i\) to a contractible subset of \(\Theta_u\), over
which \(V_m\) splits as the direct sum of \(m\) copies of \(u^*K\).
Let \(\eta_j^1, \eta_j^2, \ldots, \eta_j^m\) be the
components of \(\underline{\eta}_j\) with respect to this splitting; 
similarly for \(\xi_i^k, \, k=1, 2, \ldots, m\). Locally, \(a_{ij}\)
may be written as
\begin{equation}\label{a-ij}
a_{ij}=\int \sum_k \Tr_{u^*K}(\nabla \nabla h(u)\eta^k_j \otimes
\xi_i^k)\, ds \, dt.
\end{equation}

A direct computation from (\ref{repr}) shows that the ($\Gamma(\End (u^*K))$-valued)
matrix 
\begin{eqnarray*}
\lefteqn{(A_{ij}):=\Big(\sum_k \eta_j^k\otimes \xi_i^k\Big)}\\
&&=\Big(\sum_k (L^k\eta_j^1)\otimes (L^k\xi_i^1)\Big)\\
&&=F\oplus G,
\end{eqnarray*}
where 
\[
\begin{split}
F& :=\left(\sum_{l,l'=1}^{p'}\frac{m \epsilon_{ll'}}{2}\iota_{ll'}\left(\begin{array}{cc}
\eta_{1}^1\otimes \xi^1_{1}+\eta^1_{2}\otimes \xi^1_{2}
&\eta^1_{1}\otimes \xi^1_{2}-\eta^1_{2}\otimes \xi^1_{1}\\
-\eta^1_{1}\otimes \xi^1_{2}+\eta^1_{2}\otimes \xi^1_{1} &
\eta_{1}^1\otimes \xi^1_{1}+\eta^1_{2}\otimes \xi^1_{2}
 \end{array}\right)\right);\\
G& :=m(\eta_j\otimes \xi_i)_{p_m\geq i, j>2p'}.
\end{split}
\]
$\epsilon_{ll'}=1$ when $q_l=q_{l'} \mod n$, and equals 0
otherwise; \(\iota_{ll'}\) is the embedding of the space of
\(2\times2\) matrices into the space of \((2p')\times(2p')\) matrices
\[
\iota_{ll'}((a_{pq}))=(b_{ij}), \quad b_{ij}=a_{i-2l+2,\, j-2l'+2},
\]
\(a_{pq}=0\) if \(p,q\not\in \{1,2\}\).

Note that by unique continuation (cf. \cite{trans}), any element in
\(\ker D_{((mT), u^m)}^{J, X+\chi_H}\) or \(\cok D_{((mT),
  u^m)}^{J, X+\chi_H}\) are non-vanishing 
for all but isolated $(s,t)$, so we can find an $(s_0, t_0)$ away from
\(\Omega(u)\cup R(u)\) such that
in a small neighborhood of which, $\xi_1^1, \eta^1_1$ are nonvanishing
and not colinear with $\partial _s u$ (this is possible because
\(i_m^{-1}\underline{\xi}_1, i_m^{-1}\underline{\eta}_1\) are
\(L^2\)-orthogonal to \(\partial _s u^m\). 

We shall choose \(h\) to be of the following form. Let \(\phi_t:
U(u(s_0,t_0))\to \R^{2n}=\{(\chi_0^t, \chi_1^t, \ldots ,
\chi_{2n-1}^t)\, |\, \chi_i^t\in \R\}\) be \(t\)-dependent local coordinate charts on a
neighborhood \(U(u(s_0,t_0))\) of \(u(s_0,t_0)\in M\) such that \(\phi_t(u(s, t))=(s-s_0, 0, \ldots,
0)\), and \(\phi_t\) is smooth in \(t\). Let 
\[
h_t(x)=\beta(t-t_0)\beta_{U(u(s_0, t_0))}(x) \sum_{i, j=1}^{2n-1}B^{ij}\chi_i^t(x)\chi_j^t(x),
\] 
where \(\beta: \R\to [0,1]\) is a smooth cutoff function supported in a
small neighborhood near 0, and \(\beta_{U(u(s_0, t_0))}\) is a smooth
  cutoff function on $M$ supported on $U(u(s_0, t_0))$. $(B_{ij})$ is
  a symmetric $(2n-1)\times (2n-1)$ matrix; we denote the
  corresponding bilinear form \(B: \mbox{Sym}^2 T_{u(s_0,t_0)}M\to \R\).

This choice of $h$ apparently satisfies $h_t(u(s,t))=0; \, \nabla     h_t(u(s,t))=0$.
Moreover, since \((a_{ij})=(\int\Tr_{u*K}(\nabla\nabla h
\,A_{ij})\, ds\, dt)\), if $A_{11}$ is in $G$, it is obviously possible to
choose a $(B_{ij})$ such that 
\[
a_{11}\sim CB(\xi_1^1(s_0, t_0),\eta_1^1(s_0,t_0))\neq0
\] 
for a constant \(C\).
If $A_{11}$ is in $F$ and \(\eta^1_2, \xi^1_2\) are nontrivial
multiples of \(\eta_1^1, \xi^1_1\) or \(\xi^1_1, \eta^1_1\)
respectively, then by elementary algebra, for generic \(B\) one of
\(a_{11}\) and \(a_{12}\) must be nonvanishing, since they are of the form
\begin{gather*}
a_{11}= C\Big(B(\xi_1^1(s_0,t_0),
\eta_1^1(s_0,t_0))+B(\xi_2^1(s_0,t_0), \eta_2^1(s_0,t_0))\Big)+o(1)\\
a_{12}= C\Big(B(\xi_2^1(s_0,t_0),
\eta_1^1(s_0,t_0))-B(\xi_1^1(s_0,t_0), \eta_2^1(s_0,t_0))\Big)+o(1).
\end{gather*}
In the remaining case, \(B(\xi_2^1(s_0,t_0), \eta_2^1(s_0,t_0))\) is
either zero or can be chosen independently of \(B(\xi_1^1(s_0,t_0),
\eta_1^1(s_0,t_0))\) so it is again easy to make \(a_{11}\neq 0\). 
\hfill\(\Box\)

\section{Main Theorems from a General Invariance Theorem.}

In this section, we derive the theorems stated in section 2 from 
a general invariance theorem, Theorem \ref{thm:gen-inv} below. A rough
outline of the proof for Theorem \ref{thm:gen-inv}
is provided in \S4.4, while the details are
subjects of the remaining sections in this paper.

\subsection{Proof of Theorem \ref{thm:A}, Assuming Theorem \ref{thm:gen-inv}.}

Let \(\Lambda\subset \R\) be an interval.
A path \(\{X_\lambda\}_{\lambda\in [1,2]}\) in \({\cal X}\) is said to
be {\em $H^1$-co-directional} if \(\exists\, \langle{\cal Y}\rangle\in \op{Hom}(\mathfrak{H};\R)\) so that 
\begin{equation}\label{Y-linear}
\langle{\cal Y}_{X_\lambda}\rangle\Big|_{\ker \psi}=\alpha_\lambda\langle{\cal
  Y}\rangle\Big|_{\ker \psi}\quad \text{for some \(\alpha_\lambda\geq0\)
  \(\forall\lambda\).}
\end{equation}

\begin{*theorem}\label{thm:gen-inv} {\bf (General Invariance Theorem)}
Let \(M\) be \(w^+\)-monotone, and \((J_1, X_1), (J_2, X_2)\) be two
regular pairs. Suppose there exists
an $H^1$-co-directional path \(\{X_\lambda\}_{\lambda\in [1,2]}\) in \({\cal X}\) connecting
\(X_1, X_2\).
Then \[\begin{cases}I_F(J_1, X_1)=I_F(J_2, X_2)&\text{if
    \(\alpha_1>0\), \(\alpha_2>0\), or \(\alpha_1=\alpha_2=0\)};\\
i_{-\langle{\cal Y}\rangle}I_F(J_1, X_1)=I_F(J_2, X_2) &\text{if \(\alpha_1=0\),
  \(\alpha_2>0\)};\\
i_{-\langle{\cal Y}\rangle}I_F(J_2, X_2)=I_F(J_1, X_1) &\text{if \(\alpha_2=0\),
  \(\alpha_1>0\)}.\end{cases}\] 
\end{*theorem}

\mysubsubsection
\noindent{\it Proof of Theorem \ref{thm:A}.}
{\em Part (a)}: We show that the conditions of Theorem \ref{thm:A}
imply the conditions of Theorem \ref{thm:gen-inv}.

Let \(X_1=X+\chi_{H_1}\) and \(X_2=X+\chi_{H_2}\). Then they are
connected by a path \(\{X_\lambda\}\) in \({\cal X}\) of the form
\(X+\chi_{H_\lambda}\), since \({\cal H}\) is path-connected. The
cohomology class \([{\cal Y}_{X_\lambda}]\) is independent of
\(\lambda\) in this case. Theorem \ref{thm:A} (a) now follows immediately
from Theorem \ref{thm:gen-inv}.

{\em Part (b)} is the consequence of the next, more general theorem. 
\qed

\begin{*theorem}\label{cor:monotone}
Let \((M, \omega)\) be \(w^+\)-monotone. If \(M, f, [\gamma_0]\)
satisfy either of the additional conditions (a) or (b) below, then 
\(I_F^{f, [\gamma_0]}(M)\) is invariant under symplectic isotopies in the sense of
(\ref{sym-inv}):

(a) There is an \(\alpha\in \R\) such that 
\begin{equation}\label{monotone-M}
\omega\Big|_{\pi_2(M)}=\alpha c_1\Big|_{\pi_2(M)} \quad \text{as cohomology classes of \(M\).}
\end{equation}

(b) \((f, [\gamma_0])\) is ``monotone'' in the sense that
\(\psi_\omega^{f, [\gamma_0]}=\alpha \psi_c^{f, [\gamma_0]}\) for some
\(\alpha\in \R\).
\end{*theorem}
In the last line, we added the superscripts \(f, [\gamma_0]\) to the
notation \(\psi_\omega\), \(\psi_c\) to emphasize their dependence on \(f, [\gamma_0]\).
\begin{examples*}
{\em (a)} Condition (a) holds when \((M, \omega)\)  is monotone.

{\em (b)} When \(M\) is a K\"{a}hler-Einstein manifold, any \((f,
[\gamma_0])\) is monotone. A similar notion of ``monotone
symplectomorphisms'' was introduced for the \(n=1\) case in \cite{Sei}.
\end{examples*}

\noindent{\em Proof of Theorem 4.1.3.}
Because \(-\langle{\cal Y}_X\rangle=\psi_\omega^{f,
  [\gamma_0]}-\partial_f^*[\theta_X]\) and 
\({\cal X}\) is path-connected, if for fixed \(f,
[\gamma_0]\) there exists an \(X_0\in{\cal X}\) such that \(\langle{\cal
Y}_{X_0}\rangle\Big|_{\ker\psi}=0\), any \(X\in {\cal X}\) can be connected to \(X_0\) via a
path satisfying (\ref{Y-linear}). Invariance under symplectic
isotopies would then follow from Theorem \ref{thm:gen-inv}.

We now verify that this is true for both cases (a) and (b). In case
(b), \(\psi_\omega\Big|_{\ker\psi}=\alpha\psi_c\Big|_{\ker\psi_c}=0\),
so \(X_0\) can be simply taken to be \(0\).

For case (a), note that an examination of the maps in (\ref{c-dgrm})
shows that there is a short exact sequence
\[
0\to \mathfrak{H}'\to \mathfrak{H}\to \mathfrak{H}''\to 0,
\] 
where \(\mathfrak{H}'=\op{Image}(\iota_{f*}\circ \op{Hurewicz})\), and
\(\mathfrak{H}''=\op{Image}(\partial_f\circ \op{im})\).
On the other hand, (\ref{monotone-M}) and 
the fact that \([c_1(TM)]=\iota_f^*[c^f]\),
\([\omega]=\iota_f^*[\omega_f]\) imply
\((\psi^{f,[\gamma_0]}_\omega-\alpha\psi^{f,[\gamma_0]}_c)\Big|_{\mathfrak{H}'}=0\).
Thus, there exists a (non-unique) \(\theta\in H^1(M)\) such that 
\[
\psi^{f,[\gamma_0]}_\omega-\alpha\psi^{f,[\gamma_0]}_c=\partial_f^*\theta\Big|_{\mathfrak{H}}.
\]
Now choose \(X_0\in
{\cal X}\) such that its Calabi invariant \([\int\theta_{X_{0,t}} \, dt]=\theta\).
\qed

\subsection{A Fundamental Example.}

We explain here how Corollary \ref{cor:C} follows from Theorem
\ref{thm:A} (b). It is restated in a slightly more general form as
follows. 

\begin{*corollary}
Let \((M,\omega)\) be \(w^+\)-monotone and supppose it satisfies (\ref{monotone-M}).
Let \(\Phi=\{\phi_\lambda\, |\lambda\in[0,1]\}\) be a symplectic
isotopy with \(\phi_0=\op{id}\), \(\phi_1=f\), and
\(\gamma_0=\Phi\cdot \gamma_p\), \(\gamma_p\) being the constant path
at \(p\in M\).
Then \(I_F\) is invariant under symplectic isotopies in the sense of
(\ref{sym-inv}); furthermore, \(I_F^{f, [\gamma_0]}(M)=\tau(M)\).
\end{*corollary}
\begin{proof}
In this case, the fibration (\ref{weak-fib}) has a standard
section, and therefore via the homotopy exact sequence of the
fibration, we have the decomposition
\(\pi_1(\lc)=\pi_1(M)\times\pi_2(M)\). On the other hand, since in
this case \(T_f\simeq M\times S^1\), \(H_2(T_f;\Z)=H_1(M;\Z)\oplus
H_2(M;\Z)\), and \(\op{im}\) is simply the abelianization map. Thus, 
\begin{equation}\label{H1C:dec}
\mathfrak{H}=H_1(\lc;\Z)=H_1(M;\Z)\oplus \pi_2(M).
\end{equation} 
In terms of this decomposition, 
\[
\begin{split}
\psi_c&=0\oplus [c_1(TM)]:H_1(M)\oplus\pi_2(M)\to\Z,\\
\psi_\omega&=\theta_\Phi\oplus[\omega]:H_1(M)\oplus\pi_2(M)\to\R, 
\end{split}
\]
\(\theta_\Phi\) being the Calabi invariant of \(\Phi\).
The coefficient ring for the Floer complex is
\[
\Lambda_F=\Nov(H_1(M)\oplus \Ker(\psi_c|_{\pi_2(M)}),
-\theta\oplus[\omega];\Z).
\]
By the monotonicity of \(M\), 
Theorem \ref{cor:monotone} (a) ensures 
the invariance under symplectic isotopies. Furthermore, by
\S3.1.4 we may assume without loss of generality that \(f=\op{id}\),
and \(\gamma_0=\gamma_p\).

Thus, \(I_F^{\op{id}, [\gamma_p]}(M)\) can be computed at a regular pair
\((J, X)\) which is \(t\)-independent, with \(X=\chi_H\), \(H\) being
a small Morse function on \(M\). By a
standard argument (see e.g. \cite{floer.cmp}, \cite{salamon-zehnder}
section 7), such
regular pairs exist, and \({\cal P}(X)\), \(\cm_P(J, X), \cm_O(J, X)\) in this case
consist of \(t\)-independent elements, namely, critical points,
gradient flow lines, and periodic orbits of \(H\) respectively. Thus, 
\[
\op{CF}_{k\in \Z/2\Z}(J, X)=\op{CM}_{k\in \Z/2\Z}(H)\otimes_{\Z[H_1(M;\Z)]} \Z[H_1(M;\Z)\oplus \Ker(\psi_c|_{\pi_2(M)})],
\]
where \(\op{CM}_*(H)\) is the twisted Morse complex of \(H\). Namely,
it is defined in exactly the same way as the construction of
\(\op{CF}\) in \S2.2.2, modelling on the lift of the gradient flow of
\(H\) to the universal abelian covering \(\tilde{M}\).

We also need to compare the orientations on the moduli spaces that
enter into the definition of the boundary maps \(\partial_F\) and \(\partial_M\). As
mentioned in Remark 2.2.4, while on the Floer theory side the moduli
spaces are endowed with coherent orientations (as in \cite{FH}, and
II.7.2 later), on the Morse theory side the moduli spaces are oriented
via their interpretation as intersection spaces (as in
e.g. \cite{hl1}). Linearized versions of the adiabetic analysis which
identified the Floer-theoretic moduli spaces with the Morse-theoretic
moduli spaces also identify the Floer-theoretic coherent orientations
with coherent orientations of solution spaces of the gradient flow
equation of \(H\). On the other hand,  the arguments in \cite{schwarz}
section 3 and Appendix B shows that this analytic version of
orientation agrees with the geometric version of orientation.
Thus, the complexes 
\(\op{CF}_*(J,X)\) and \(\op{CM}_*(H)\) have the same boundary map \(\partial_F=\partial_M\). 

Lastly, \(\cm_O(J, \chi_H)=\emptyset\) since the gradient
flow of the real valued function \(H\) has no periodic orbits. 
Thus, \[I_F^{\op{id},
  [\gamma_p]}(M)=\iota_*\tau(\op{CM}(H))=\iota_*\tau(M).\]
The second equality above is the well-known equivalence between 
torsions of Morse complexes (of real-valued Morse functions) and
torsions of cell chain complexes. See e.g. \cite{Mil}.
\end{proof}

\begin{remark*}\label{compare:LO}
In \cite{le-ono}, Le and Ono considered a Floer theory under 
assumptions similar to Theorem \ref{cor:monotone}. There 
they showed that in this case, a (twisted version of) Floer
homology is equivalent, up to a change of coefficients, to a 
Novikov homology (a version of twisted Morse homology).

Their definition of Floer complex is 
different from ours in that they used a smaller covering of
the loop space. The advantage of this choice is that the associated
{\it homology} retains more information (i.e. is less likely to be trivial),
though it has worse invariance properties. In this paper our main
concern is {\it torsion}, and the torsions associated to different coverings
are related by a simple change of coefficients induced by maps
between covering groups. For our purpose the larger 
covering is therefore the better choice.
\end{remark*}

\subsection{Regular Homotopy of Floer Systems (RHFS).}
The proof of Theorem
\ref{thm:gen-inv} will rely on the existence of a certain path of formal
vector fields \(\{{\cal V}_\lambda\}\) connecting the two generated by
\((J_1, X_1)\) and \((J_2, X_2)\). 
The purpose of this subsection is to describe such paths of 
vector fields, called ``regular homotopy of Floer systems''.  

We shall work on an abstract level for the rest of this section, as we
did in \S2.1. The eventual goal is to give a formal proof of
the invariance of \(I_F\) (cf. Proposition \ref{prop:abs-inv}), which applies to any Floer theory satisfying
the outline in \S2.1.

\mysubsubsection[CHFSs.]
A {\em co-directional homotopy of Floer systems} (CHFS for short) is 
a 1-parameter family of formal flows 
\(\{(\lc, \mathfrak{H}, \ind; {\cal Y}_\lambda, {\cal
  V}_\lambda)\}_{\lambda\in \Lambda=[1,2]}\) (abbreviated as
\(\{{\cal V}_\lambda\}_{\lambda\in \Lambda=[1,2]}\)) such that:
\begin{itemize}
\item \({\cal V}_1, {\cal V}_2\) generate Floer
systems with the same SF-homomorphism \(\psi\), and 
\item \(\exists\) \(\lambda\)-independent \(\langle{\cal Y}\rangle\in \op{Hom}(\mathfrak{H}, \R)\),
  so that \(\langle{\cal Y}_\lambda\rangle\Big|_{\ker
    \psi}=\alpha_\lambda\langle{\cal Y}\rangle\Big|_{\ker \psi}\) 
for some \(\alpha_\lambda\geq0\) \(\forall\lambda\). 
\end{itemize}
Such a \(\langle{\cal Y}\rangle\) will be called a {\em c-class} of
the CHFS.
Notice that for a CHFS, the associated
Novikov ring \(\Nov (\ker \psi, -\langle{\cal Y}_\lambda\rangle;\Z)=\Nov (\ker \psi,
-\langle{\cal Y}\rangle; \Z)\) remains the same through the homotopy.

Given a closed interval \(S\subset \Lambda\), we say that \(\{{\cal
  V}_\lambda\}_{\lambda\in S}\) is a {\em sub-homotopy} over \(S\) of
the CHFS \(\{{\cal V}_\lambda\}_{\lambda\in \Lambda}\).

Define the parametrized moduli spaces \({\cal
  P}^\Lambda=\bigcup_\lambda{\cal P}_\lambda\),
\(\cm_P^{\Lambda, k+1}=\bigcup_\lambda \cm_{P, \lambda}^k\),
\(\cm_O^{\Lambda, k+1}=\bigcup_\lambda \cm_{O, \lambda}^k\), where \({\cal
  P}_\lambda, \cm_{P, \lambda}, \cm_{O, \lambda}\) are moduli spaces
of the formal flow generated by \({\cal V}_\lambda\). We denote the
projection of these parametrized moduli spaces to \(\Lambda\) by
\(\Pi_\Lambda\) in general. In addition, define
\[
{\cal P}^{\Lambda, deg}={\cal P}^\Lambda\backslash \bigcup_{\lambda\in\Lambda}{\cal P}_{\lambda, ndg}.
\]

A CHFS is {\em regular} (i.e. is an RHFS)
if the parametrized moduli spaces satisfy all the properties (RHFS*) below.
These properties are written modelling on the generic behavior of a 1-parameter
family of flows in finite-dimensional Morse-Novikov theory.
\subsubsection{Properties of \({\cal P}^\Lambda\), \({\cal
  P}^{\Lambda, deg}\).}
\medbreak
\noindent{\bf (RHFS1)} {\it 
\({\cal P}^\Lambda\) is a compact Zariski smooth 1-manifold with
boundary \(\p{\cal P}^\Lambda={\cal P}_2-{\cal P}_1\). The space \({\cal P}^{\Lambda,
  deg}\) coincides with the set of
critical points of \(\Pi_\Lambda\); in addition, it 
consists of finitely many points in the interior of
\({\cal P}^\Lambda\), which are extrema of
\(\Pi_\Lambda\).
}
\medbreak
In fact, the identifcation of \({\cal P}^{\Lambda, deg}\) with
critical points of \(\Pi_\Lambda\) is a direct consequence of the
Zariski smoothness of \({\cal P}^\Lambda\). To see
this, note that the deformation operator for \({\cal P}^\Lambda\) at
\(x_\lambda\), denoted \(\hat{A}_{x_\lambda}\), is a rank 1 extension
of \(A_{x_\lambda}\) given by
\[
\hat{A}_{x_\lambda}(\alpha, \xi)=\alpha \p_\lambda{\cal
  V}_\lambda+A_{x_\lambda}\xi \quad \text{for \(\alpha\in \R\),
\(\xi\in T_{x_\lambda}C\)}.
\]
Apparently, \(\ker A_{x_\lambda}\subset \ker
\hat{A}_{x_\lambda}\). Zariski smoothess of \({\cal P}^\Lambda\)
implies that the latter is 1-dimensional, and therefore the former 
is either 0-dimensional or
1-dimensional. In the first case, \(x_\lambda\) is a nondegenerate
element of \({\cal P}_\lambda\); in the second case, it is said to be {\it
  minimally degenerate}. Furthermore, notice that
\(d\Pi_\Lambda:T_{x_\lambda}{\cal P}^\Lambda\to T_\lambda\Lambda\) is given by projecting \(\ker
\hat{A}_{x_\lambda}\subset\R\oplus T_{x_\lambda}\lc\) to the
\(\R\)-summand. This projection is trivial precisely when
\(x_\lambda\) is minimally degenerate, and hence \({\cal P}^{\Lambda, deg}\) consists precisely of
critical points of the map $\Pi_{\Lambda}$. 
\medbreak

\noindent{\bf (RHFS1i)} {\sl (Injectivity)}
{\it The restriction \(\Pi_\Lambda\Big|_{{\cal P}^{\Lambda,
    deg}}\) is injective.}
\medbreak

Elements of \({\cal P}^{\Lambda, deg}\) are
called {\em death-births}: a local maximum of \(\Pi_\Lambda\) is
called a {\em death}, and a local minimum is called a {\em birth}.
We denote \[\Lambda_{db}:=\Pi_\Lambda({\cal P}^{\Lambda,
    deg})\subset\Lambda.\]
A small neighborhood \(S\) of a subset \(P\subset
\Lambda_{db}\) is called a {\em death-birth neighborhood}. The half
\(\Pi^\Lambda ({\cal P}^S \backslash {\cal P}^{S, deg} )\subset S\) is called a {\em
  birth-neighborhood}, and the other half \(S\backslash\Pi^\Lambda
({\cal P}^S)\) is called a {\em death-neighborhood}.

Let \(\tilde{\cal
  P}^\Lambda\), \(\tilde{\cal
  P}^{\Lambda, deg}\) be respectively the lifts of \({\cal
  P}^\Lambda\), \({\cal P}^{\Lambda, deg} \subset \Lambda\times \lc\) to \(\Lambda\times\tilde{\lc}\).

Note that \(\ind, \op{gr}\) define locally constant functions on \({\cal
  P}^\Lambda\backslash{\cal P}^{\Lambda, deg}\) and \((\tilde{\cal
  P}^\Lambda\backslash\tilde{\cal P}^{\Lambda, deg})^{\times 2}\). Let
\(\op{gr}_+\), \(\op{gr}_-: (\tilde{\cal P}^\Lambda)^{\times 2}\to \Z\) denote the
upper-semicontinuous and the lower-semicontinuous
extension of \(\op{gr}\) respectively. Similarly for \(\ind_+, \ind_-:
{\cal P}^\Lambda\to \Z/2\Z\).
It follows from the definition of \(\ind, \op{gr}\) via spectral flow
and the fact that \(\ker A_{x_\lambda}=1\) for any \(x_\lambda\in
{\cal P}^{\Lambda, deg}\) that:
\begin{equation}\label{RHFS-1g}\begin{split}
\op{gr}_-((x_\lambda,[w_\lambda]),
\cdot) &=\op{gr}_+(\cdot,(x_\lambda, [w_\lambda]))-1;\\
\op{gr}_-(\cdot, (x_\lambda,[w_\lambda])) &=\op{gr}_+((x_\lambda,
[w_\lambda]), \cdot)-1;\\
\ind_-(x_\lambda) &=\ind_+(x_\lambda)-1\quad \text{for any \((x_\lambda, [w_\lambda])\in \tilde{\cal P}^{\Lambda, deg}\).}
\end{split}\end{equation}

By (RHFS1, 1i), \({\cal P}^\Lambda\backslash {\cal P}^{\Lambda, deg}\)
consists of finitely many path components, and \(\Pi_\Lambda\)
restricts to a diffeomorphism from each such component \({\bf x}\) to its 
image in \(\Lambda\), which we denote by \(\Lambda_{\bf x}\). Denote the set
of such path components by \(\aleph_\Lambda\). Similarly, let
\(\tilde{\aleph}_\Lambda\) be the set of path components of 
\(\tilde{\cal P}^\Lambda\backslash \tilde{\cal P}^{\Lambda, deg}\). Given
\(\lambda\in \Lambda\), we denote by \(x_\lambda\)
\[x_\lambda=\Pi_\Lambda^{-1}(\lambda)\cap {\bf x}.\] 
Since \(\ind\) is constant over \({\bf x}\), we denote \(\ind ({\bf x})=\ind (x_\lambda)\) for an
arbitrary \(\lambda\in \Lambda_{\bf x}\). Similarly for \(\op{gr}\).

\mysubsubsection[Weight truncation associated to \(\langle{\cal
  Y}\rangle\) and \(e_{\cal P}\).]
Before proceeding to the next property of RHFS, we shall describe a truncation
of the moduli spaces by certain weight determined by \(\langle{\cal
  Y}\rangle\). This is essentially equivalent to truncating by energy, but
is better suited for a uniform description of the moduli spaces in a
CHFS. 

Given a formal flow \((\lc, \mathfrak{H}, \ind; {\cal Y}_\chi, {\cal
  V}_\chi)\), choose a \(Y\in
\op{Hom}(\mathfrak{H}, \R)\) and an basis \(e_{\cal P}\) for
\(\op{CF}\) as in \S\ref{def:torsion} (3). Namely, 
a lifting \(\op{Lf}: {\cal P}\to \tilde{\cal P}\), and a choice of
\(A_0\in \mathfrak{H}\) such that \(\psi(A_0)=2\mathbb{N}_\psi\) when
\(\mathbb{N}_\psi\neq0\). In this case, this choice decides a
splitting
\[
\mathfrak{H}=\ker \psi\oplus \Z A_0.
\]
A connecting flow line \(u\) belongs to certain
\(\cm_{P}(\op{Lf}(x), A\cdot \op{Lf}(y))\) for a unique
\(A\in \mathfrak{H}\), and \([u]_{\mathfrak{H}}:=A\) is said
to be the {\em \(\mathfrak{H}\)-class} of \(u\). On the other hand,
when \(u\) is a closed orbit, then
\([u]_{\mathfrak{H}}:=A\) when \(u\in \cm_{O}(A)\).
The {\em \(\ker\psi\)-class} of \(u\), denoted \([u]_{\ker\psi}\), is
the projection of \([u]_{\mathfrak{H}}\) to \(\ker\psi\) according to
the above splitting determined by \(A_0\); in the case when
\(\mathbb{N}_\psi=0\), \([u]_{\ker\psi}:=[u]_{\mathfrak{H}}\).
Similarly, for \(A\in \mathfrak{H}\), \([A]_{\ker\psi}\) denotes its
projection to \(\ker \psi\), in accordance with the above splitting when
\(\mathbb{N}_\psi\neq0\) (otherwise \(\ker \psi=\mathfrak{H}\), and \([A]_{\ker\psi}=A\)).
The {\em weight} of a \(u\in \cm_P\) or
\(\cm_O\), denoted \(\op{wt}_{-Y,
  e_{\cal P}}(u)\), is defined to be \(-
  Y[u]_{\ker\psi}\). 
 
For a CHFS, we define the weight for elements in \(\cm_P^\Lambda\) or
\(\cm_O^\Lambda\) by fixing a lifting \({\cal P}^\Lambda \to \tilde{\cal
  P}^\Lambda\) (which in turn gives a consistent system of liftings
\(\op{Lf}: {\cal P}_\lambda \to \tilde{\cal
  P}_\lambda\) for each \(\lambda\in \Lambda\)), fixing a
\(\lambda\)-independent \(A_0\), and letting
\(Y=\langle{\cal Y}\rangle\) \(\forall \lambda\) for a c-class
\(\langle{\cal Y}\rangle\). Such a weight is said to be {\em adapted} to the CHFS
\(\{{\cal V}_\lambda\}\).

Given \(\Re\in \R\), let \[\cm_P^{S, k}({\bf x}, {\bf y}; \op{wt}_{-\langle{\cal Y}\rangle,
  e_{\cal P}}\leq\Re),\,  \hat{\cm}_{P,\lambda}^k(x_\lambda, y_\lambda;\op{wt}_{-\langle{\cal Y}\rangle,
  e_{\cal P}}\leq\Re)\] be the moduli spaces consisting of flows with
weight \(\leq \Re\) and appropriate additional constraints. Similarly
for other variants of moduli spaces and \(\cm_O\), \(\cm_O^\Lambda\).

Let \(\tilde{\cal A}_\lambda\) be the primitive of \({\cal
  Y}_\lambda\) with \(\tilde{\cal A}_\lambda(\gamma_0)=0\). Then 
\[
{\cal E}(u_\lambda)=\begin{cases}
-\langle{\cal Y}_\lambda\rangle[u_\lambda]_{\mathfrak{H}}-\tilde{\cal
  A}_\lambda(\op{Lf}(y_\lambda))+\tilde{\cal A}_\lambda(\op{Lf}(x_\lambda))
&\text{when \(u_\lambda\in \cm_{P, \lambda}(x_\lambda,y_\lambda)\)};\\
-\langle{\cal Y}_\lambda\rangle[u_\lambda]_{\mathfrak{H}} &\text{when
  \(u_\lambda\in \cm_{\lambda, O}\)}.\end{cases}
\]
Thus, we have the following relation between energy and weight:
\begin{equation}\label{E-bdd}
{\cal E}(u_\lambda)=\begin{cases}\alpha_\lambda\op{wt}_{-\langle{\cal Y}\rangle,
  e_{\cal P}}(u_\lambda)+\beta_\lambda\psi([u_\lambda]_{\mathfrak{H}})&+\tilde{\cal A}_\lambda(\op{Lf}(x_\lambda))-\tilde{\cal
  A}_\lambda(\op{Lf}(y_\lambda))\\
&\text{when \(u_\lambda\in \cm_{P, \lambda}(x_\lambda,y_\lambda)\)};\\
\alpha_\lambda\op{wt}_{-\langle{\cal Y}\rangle,
  e_{\cal P}}(u_\lambda)+\beta_\lambda\psi([u_\lambda]_{\mathfrak{H}})
&\text{when
  \(u_\lambda\in \cm_{\lambda, O}\)},
\end{cases}
\end{equation}
where \(\beta_\lambda=-\langle{\cal
  Y}_\lambda\rangle(A_0)/(2\mathbb{N}_\psi)\) is continuous in \(\lambda\).

For elements in \(\cm_{P,\lambda}^k\) or \(\cm_O^k\), all terms above
are fixed except for the first term involving the weight. Thus, the
filtration by weights on these moduli spaces are equivalent to
filtration by energy when \(\alpha_\lambda>0\), while when
\(\alpha_\lambda=0\), the energy filtration is the trivial
filtration. 
Furthermore, as the constant terms above vary continuously in
\(\lambda\), and \({\cal P}^\Lambda\) is compact, the above formula
gives a uniform bound on energy by weight.

\subsubsection{Structure of \(\cm_{P, \lambda}\) for \(\lambda\in \Lambda_{db}\).}

In the case when \(x, y\in {\cal P}\) are either nondegenerate or
minimally-degenerate, let \(\cm_P^k(x, y)\) be the moduli space of
connecting flow lines which lift to be paths from \((x, [w])\) to
\((y, [v])\) in \(\tilde{\lc}\), with \(\op{gr}_+((x, [w]), (y,
[v])=k\).
\medbreak
\noindent{\bf (RHFS2d)} {\sl (moduli at Death-births)} {\it 
For any \(\lambda\in \Lambda_{db}\), any integer \(k<2\), and any pair \(x_\lambda,
y_\lambda\in {\cal P}_\lambda\), 
\(\cm_{P, \lambda}((x_\lambda, [w_\lambda]), (y_\lambda,[v_\lambda]))\) is
a Zariski smooth manifold of dimension \(k\).
Furthermore, for any real constant \(\Re\), 
\(\hat{\cm}_P^{0}(x_\lambda, y_\lambda; \op{wt}_{-\langle{\cal Y}\rangle,
  e_{\cal P}}\leq\Re)\) consists of finitely many Zariski smooth points.
}
\medbreak

When \(x_\lambda, y_\lambda\) are both nondegenerate, the above
statement is part of property (FS2). 

\medbreak

\subsubsection{The Structure of \(\cm_P^\Lambda\).}
Consider the parameterized moduli spaces over an interval \(S\subset \Lambda\): 
\[\hat{\cm}_P^S(({\bf x, [w]}), ({\bf y,
  [v]})):=\bigcup_{\lambda\in S\cap \Lambda_{\bf x}\cap \Lambda_{\bf y}}\hat{\cm}_{P, \lambda}((x_\lambda,
[w_\lambda]), (y_\lambda, [v_\lambda]));\quad \hat{\cm}_O^{S}:=\bigcup_{\lambda \in S}\hat{\cm}_{O, \lambda}.\]
These are included in parameterized versions of moduli spaces of broken
trajectories/orbits, 
\[
\hat{\cm}_P^{S, +}((\underline{\bf x, [w]}), (\underline{\bf y,
  [v]})):=\bigcup_{\lambda \in \bar{S}\cap \bar{\Lambda}_{\bf x}\cap
  \bar{\Lambda}_{\bf y}}\hat{\cm}^+_{P, \lambda}((\underline{\bf x, [w]}), (\underline{\bf y,
  [v]})),
\quad \hat{\cm}_O^{S, +}=\bigcup_{\lambda \in \bar{S}}\hat{\cm}_{O, \lambda}^+,
\]
where \(\hat{\cm}^+_{P, \lambda}((\underline{\bf x, [w]}), (\underline{\bf y,
  [v]}))\) is the reduced, broken-trajectories variant of:
\[
\cm_{P, \lambda} ((\underline{\bf x, [w]}), (\underline{\bf y,
  [v]})):=\cm_{P, \lambda}^{(\sigma_1, \sigma_2)}((x_\lambda,
[w_\lambda]), (y_\lambda, [v_\lambda])), \]
with weights \(\sigma_1\), \(\sigma_2\) chosen such that \[
\op{gr}^{(\sigma_1, 0)}((x_\lambda, [w_\lambda]), \cdot)=\op{gr} (({\bf x, [w]}), {\bf \cdot});\quad 
\op{gr}^{(\sigma_2, 0)}((y_\lambda, [v_\lambda]), \cdot)=\op{gr}(({\bf y, [v]}), {\bf\cdot}).\] 
Other variants of the moduli spaces such as
\(\hat{\cm}_P^{S, k, +}(\underline{\bf x}, \underline{\bf y}; \op{wt}_{-\langle{\cal Y}\rangle,
  e_{\cal P}}\leq\Re)\) can be defined in a similar way. 
The definition of chain topology extends to these parameterized
versions in an obvious manner.
\medbreak

\noindent {\bf (RHFS2)} {\it 
Given \({\bf x}, {\bf y}\in \aleph_\Lambda\), an integer \(k\leq 1\),
and an interval \(S\subset\Lambda\), 
the parameterized moduli space \(\hat{\cm}_P^{S, k}({\bf x}, {\bf y
})\) is a Zariski smooth manifold of
dimension \(\op{gr}({\bf x},{\bf y})\); furthermore, 
the parameterized moduli space of broken trajectories \(\hat{\cm}_P^{\Lambda, k, +}(\underline{\bf x}, \underline{\bf y}; \op{wt}_{-\langle{\cal Y}\rangle,
  e_{\cal P}}\leq \Re)\) is compact for any \(\Re\).
}
\medbreak
Notice that since \(\hat{\cm}_P^{S, 0,+}(\underline{\bf x}, \underline{\bf y}; \op{wt}_{-\langle{\cal Y}\rangle,
  e_{\cal P}}\leq\Re)=\hat{\cm}_P^{S, 0}({\bf x}, {\bf y}; \op{wt}_{-\langle{\cal Y}\rangle,
  e_{\cal P}}\leq\Re)\), (RHFS2) implies that this moduli space consists
of finitely many nondegenerate
points.
\medbreak
\noindent {\bf (RHFS2i)} {\sl (Injectivity)} {\it 
The projection \(\Pi_\Lambda\) is injective on the space
\[\hat{\cm}_P^{S,0}:=\coprod_{{\bf x}, {\bf y}\in \aleph_\Lambda}\hat{\cm}_P^{S, 0}({\bf x}, {\bf y}).\]
}
\medbreak
Elements in this space are called {\em handle-slides}. They are said
to be of {\em type II} if \({\bf x}={\bf y}\). Otherwise they are of
{\em type I}. 
Let 
\[\begin{split}
\Lambda_{hs}^\Re &:=\Pi_\Lambda\Big(\coprod _{{\bf x}, {\bf y}\in
  \aleph_\Lambda}\hat{\cm}_P^{\Lambda, 0}({\bf x}, {\bf y};
\op{wt}_{-\langle{\cal Y}\rangle, e_{\cal P}}\leq \Re)\Big);\\
\Lambda_{hs:ii}^\Re &:=\Pi_\Lambda\Big(\coprod _{{\bf x}\in \aleph_\Lambda}\hat{\cm}_P^{\Lambda, 0}({\bf x}, {\bf x}; \op{wt}_{-\langle{\cal Y}\rangle,
  e_{\cal P}}\leq \Re)\Big);\\
\Lambda_{hs} & := \bigcup_{\Re}\Lambda_{hs}^\Re; \quad \Lambda_{hs:ii}:=\bigcup_{\Re}\Lambda_{hs:ii}^\Re.
\end{split}\]
These subsets of \(\Lambda\) are disjoint from \(\Lambda_{db}\).
\medbreak

To state the next property, 
introduce the following 0-dimensional stratum in \(\hat{\cm}_P^{\Lambda, 1, +}(\underline{\bf x}, \underline{\bf y}; \op{wt}_{-\langle {\cal Y}\rangle,
  e_{\cal P}}\leq \Re)\): let \({\mathbb T}_{P,hs-m}({\bf x}, {\bf
  y};\Re)\) be the subset consisting of \(k\)-th broken trajectories
connected at {\em nondegenerate} critical points, with \(k>1\).
(RHFS2i) implies that such broken trajectories occur (only) with Type II
handleslides (namely, \(\Pi_\Lambda{\mathbb T}_{P,hs-m}({\bf x}, {\bf
  y};\Re)\) surjects to \(\Lambda_{hs: ii}^\Re\)).

Let \({\mathbb I}_P({\bf x}, {\bf y} ;\Re) :=\hat{\cm}_P^{\Lambda, 1, +}(\underline{\bf x}, \underline{\bf y}; \op{wt}_{-\langle {\cal Y}\rangle,
  e_{\cal P}}\leq \Re)\backslash \hat{\cm}_P^{\Lambda, 1}({\bf x}, {\bf y}; \op{wt}_{-\langle {\cal Y}\rangle,
  e_{\cal P}}\leq \Re)\).

\medbreak
\noindent{\bf (RHFS2c)} {\sl (Corner structure)}
{\it The moduli space of broken trajectories 
\(\hat{\cm}_P^{\Lambda, 1, +}(\underline{\bf x}, \underline{\bf y}; \op{wt}_{-\langle {\cal Y}\rangle,
  e_{\cal P}}\leq \Re)\) is lmb along \({\mathbb I}_P({\bf x}, {\bf
  y} ;\Re)\backslash {\mathbb T}_{P, hs-m}({\bf x}, {\bf y} ;\Re)\).}

\subsubsection{The Structure of \(\cm_O^\Lambda\).}
\noindent{\bf (RHFS3)} {\it
Given \(k\leq1\), 
the parameterized moduli space \(\cm_O^{S, k+1}(A)\)
is a Zariski smooth manifold of dimension \(k+1\) with a semi-free
\(S^1\) action. Furthermore, The space \(\hat{\cm}_O^{S, k, +}
(\op{wt}_{-\langle {\cal Y}\rangle}\leq\Re)\) is compact.}
\medbreak
\noindent{\bf (RHFS3c)} {\sl (Corner structure)} {\em
 The moduli space of broken orbits 
\(\hat{\cm}_O^{\Lambda, 1, +}(\op{wt}_{-\langle {\cal Y}\rangle,
  e_{\cal P}}\leq \Re))\) is lmb along \({\mathbb I}_O(\Re)\backslash
{\mathbb T}_{O, hs-m}(\Re)\), where \({\mathbb I}_O(\Re)\), 
\({\mathbb T}_{O, hs-m}(\Re)\) are defined similarly to \({\mathbb I}_P({\bf x}, {\bf
  y} ;\Re)\), \({\mathbb T}_{P, hs-m}({\bf x}, {\bf y} ;\Re)\).}

\subsubsection{Orientations.}

Let \(\sign (p:P)\) denote the sign of \(p\in
\partial P\) oriented as a boundary point of the oriented 1-manifold 
\(P\),
and the sign of an \(x_\lambda\in {\cal P}^{\Lambda, deg}\) be defined by
\[
\sign (x_\lambda)=\begin{cases}+1&\text{when \(x_\lambda\) is a
    death;}\\-1&\text{when \(x_\lambda\) is a birth.}\end{cases}
\]
Notice that that \(\sign (x_\lambda)=\sign (\lambda:\Lambda_{\bf x})\) when
\(x_\lambda\in {\cal P}^{\Lambda, deg}\) are on the boundary of the
path component \({\bf x}\subset{\cal P}\backslash {\cal P}^{\Lambda,
  deg}\). 

The notions of coherent orientation and 
grading-compatible orientation can be extended to parameterized moduli
spaces, as will be explained in II.7.2.
\medbreak
\noindent{\bf (RHFS4)}{\em 
  The parameterized moduli spaces of connecting flow lines
\(\{\hat{\cm}_P^{S, k}({\bf x}, {\bf y})\}\) can be endowed with a
coherent orientation, and the parameterized moduli space of closed
orbits \(\cm_O^{S, 2}\) can be endowed with a grading-compatible
orientation. With respect to these orientations, the strata
\[\text{ \({\mathbb I}_P({\bf x}, {\bf
  y} ;\Re)\backslash {\mathbb T}_{P, hs-m}({\bf x}, {\bf y} ;\Re)\) and \({\mathbb I}_O(\Re)\backslash
{\mathbb T}_{O, hs-m}(\Re)\),}\] 
oriented respectively as boundary
components of the parameterized moduli spaces
\[\text{\(\hat{\cm}_P^{\Lambda, 1, +}(\underline{\bf x}, \underline{\bf y}; \op{wt}_{-\langle {\cal Y}\rangle,
  e_{\cal P}}\leq \Re)\) and \(\hat{\cm}_O^{\Lambda, 1, +}(\op{wt}_{-\langle {\cal Y}\rangle,
  e_{\cal P}}\leq \Re)\), }\]
are expressed in terms of lower dimensional moduli spaces as follows:
\begin{gather*}
{\mathbb I}_P({\bf x}, {\bf
  y} ;\Re)\backslash {\mathbb T}_{P, hs-m}({\bf x}, {\bf y} ;\Re)
= {\mathbb J}_P({\bf x}, {\bf y} ;\Re)\sqcup {\mathbb T}_{P,db}({\bf
x}, {\bf y};\Re)\sqcup {\mathbb T}_{P,hs-s}({\bf
x}, {\bf y};\Re);\\ 
{\mathbb I}_O(\Re)\backslash
{\mathbb T}_{O, hs-m}(\Re)={\mathbb J}_O(\Re)\sqcup {\mathbb
  T}_{O,db}(\Re)
\sqcup {\mathbb T}_{O,hs-s}(\Re),
\end{gather*}
where 
\begin{equation}\label{sign-J-P}
{\mathbb J}_P(S,{\bf x}, {\bf y} ;\Re) =
\coprod_{\lambda\in \partial (S\cap \Lambda_{\bf x}\cap
  \Lambda_{\bf y})}\sign (\lambda: S\cap \Lambda_{\bf x}\cap
  \Lambda_{\bf y})\hat{\cm}_{P,\lambda}^{0}(\underline{\bf x},
  \underline{\bf y};\op{wt}_{-\langle {\cal Y}\rangle,
  e_{\cal P}}\leq\Re),
\end{equation}
\begin{eqnarray}\label{sign:TPd}
&{\mathbb T}_{P,db} ({\bf
x}, {\bf y};\Re) =& \hat{\cm}_P^{\Lambda, 1, +}(\underline{\bf x}, \underline{\bf y}; \op{wt}_{-\langle {\cal Y}\rangle,
  e_{\cal P}}\leq \Re)\cap \nonumber\\&
 \coprod_{z_\lambda\in {\cal P}^{\Lambda, deg}}  \coprod_{k=0}^{\infty}&
\sign (z_\lambda)(-1)^{k}
\hat{\cm}_{P, \lambda}^0(x_\lambda, z_\lambda)\times \hat{\cm}_{P,
  \lambda}^0(z_\lambda, z_\lambda)^{\times k}\times \hat{\cm}_{P, \lambda}^0(z_\lambda, y_\lambda),\qquad
\end{eqnarray}
\begin{eqnarray}\label{sign:TPs}
&{\mathbb T}_{P,hs-s}({\bf x}, {\bf y};\Re)= & \hat{\cm}_P^{\Lambda, 1, +}(\underline{\bf x}, \underline{\bf y}; \op{wt}_{-\langle {\cal Y}\rangle,
  e_{\cal P}}\leq \Re)\cap \nonumber\\
&&\coprod_{i=0,1}
\coprod _{{\bf z}\in \aleph_\Lambda}(-1)^{i}
\hat{\cm}^{S, i}_P({\bf x}, {\bf z})\times_{S}
\hat{\cm}^{S, 1-i}_P({\bf z},{\bf y}),
\end{eqnarray}
\begin{equation}
{\mathbb J}_O(S, \Re)=\coprod_{\lambda\in \partial S}\sign (\lambda:
S)\hat{\cm}_{O,\lambda}^0(\op{wt}_{-\langle {\cal Y}\rangle, e_{\cal P}}\leq\Re),
\end{equation}
\begin{equation}\begin{split}
{\mathbb T}_{O, db}(\Re)=&\hat{\cm}_O^{\Lambda, 1, +}(\op{wt}_{-\langle {\cal Y}\rangle,
  e_{\cal P}}\leq \Re)\, \cap \\
&\coprod_{y_\lambda\in {\cal P}^{\Lambda,  deg}}\coprod_{k\in \Z^+}
\sign (y_\lambda) (-1)^{k+\ind_+ (y_\lambda)} \hat{\cm}_{P,
  \lambda}^0(y_\lambda, y_\lambda)^{\times k},\\
\end{split}\end{equation}
\begin{equation}
{\mathbb T}_{O, hs-s}(\Re)=
\coprod_{{\bf x}\in \aleph_\lambda}(-1)^{\ind ({\bf
    x})+1}\hat{\cm}_P^{S,0}({\bf x},{\bf x}; \op{wt}_{-\langle {\cal Y}\rangle, e_{\cal P}}\leq\Re).
\end{equation}}

\begin{remarks*}
\({\mathbb T}_{P, db}\), \({\mathbb T}_{O, db}\)
consist of broken trajectories/orbits connected at a
degenerate critical point (a death-birth); \({\mathbb T}_{P, hs-s}\), \({\mathbb T}_{O,
  hs-s}\) consist of once-broken trajectories/orbits connected at a
non-degenerate critical point, which occur with handle-slides.

The local structure of the parameterized moduli spaces of broken
trajectories/orbits near \({\mathbb T}_{P, hs-m}\), \({\mathbb T}_{O,
  hs-m}\) is not understood. In fact, it seems that this 
is very sensitive to perturbations,
and thus is not described by a universal formula. (See II.1.2.5
for more discussion of this difficulty from the perspective of gluing
theory). 

\end{remarks*}

We shall show later (section 6 and Part II)
that under the assumptions of Theorem
\ref{thm:gen-inv}, there is an RHFS connecting the Floer
systems associated with the pairs \((J_1, X_1)\), \((J_2, X_2)\).

\subsection{Invariance from Regular Homotopy of Floer Systems.}
Let \(\{(\lc, \mathfrak{H}, \ind; {\cal Y}_\lambda, {\cal
  V}_\lambda)\, |\, \lambda\in [1,2]\}\) be an RHFS. If \((\lc, \mathfrak{H}, \ind; {\cal Y}_\lambda, {\cal
  V}_\lambda)\) is a Floer system, let \((\op{CF}({\cal V}_\lambda), \p_F({\cal
  V}_\lambda))\) denote the associated Floer
complex. Let \(\eta_F({\cal V}_\lambda)\), \(\zeta_F({\cal V}_\lambda)\), \(I_F({\cal V}_\lambda)\)
etc be similarly defined. We are interested in how they vary with \(\lambda\).
By understanding this, we will see in this subsection that
the existence of an RHFS satisfying some 
extra condition ((NEP) in \ref{neq-pert} below), imply that
\(I_F({\cal V}_1)=I_F({\cal V}_2)\). 

The condition (NEP) roughly says that the induced CHFSs
on finite-cyclic covers of \(\lc\) may be perturbed into RHFSs.
It is needed because (RHFS*) says nothing about the structure of the
parameterized moduli spaces near \(\mathbb{T}_{P, hs-m}\) and
\(\mathbb{T}_{O, hs-m}\); thus from an RHFS one may only conclude 
the invariance of \(I_F\) ``up to first order''. By considering
induced CHFSs on finite-cyclic covers of \(\lc\), the higher order
contributions are recovered from the relation between torsion
invariants of a space and its finite-cyclic covers.

The contents of this subsection are included only for the sake of
completeness: they are straightforward
consequences of (RHFS*), and are by no means new. 
We shall therefore
leave out some details of the arguments which the reader may easily
fill in, or by consulting the literature. A similar
account with complete details may be found in \cite{Hut}.

\subsubsection{Preparations.}
Recall that throughout an RHFS with c-class \(\langle{\cal Y}\rangle\), 
the associated Novikov ring is \(\Lambda_F=\Nov (\ker \psi, -\langle{\cal Y}\rangle;
\Z)\). Its ring of fractions \(Q(\Lambda_F)\) embeds via (\ref{Q-emb})
into a Novikov ring, which we denote by \(\Omega_F\). Both
\(Q(\Lambda_F)\) and \(\Omega_F\) decompose as finite direct sums of
fields, 
\[Q(\Lambda_F)=\bigoplus_\kappa  F_\kappa, \quad \Omega_F=\bigoplus _\kappa \mathbb{K}_\kappa, \] 
where each \(\kappa\) corresponds to an equivalence of
characters \(\kappa: \op{Tors}(\ker \psi)\to \C^\times \), each
\(\mathbb{K}_\kappa\) is again a Novikov ring, and the embedding
(\ref{Q-emb}) is compatible with the decomposition.

For the invariance proof, it is convenient to work with the Novikov
ring \(\Omega_F\) instead of \(Q(\Lambda)\), namely, we identify
\(\mathfrak{Q}\) as an element in \(\Omega_F\) via the embedding
(\ref{Q-emb}), for \(\mathfrak{Q}=\partial _F\), \(\tau_F\),
\(\eta_F\), \(\zeta_F\), or \(I_F\). It is also sometimes 
convenient to work with a field component of \(\Omega_F\) at a time. 
In such case, we denote by
\(\mathfrak{Q}^\kappa=\mathfrak{Q}\otimes_{\Lambda_F}\mathbb{K}_\kappa\),
namely, the image of \(\mathfrak{Q}\) under the composition of (\ref{Q-emb}) and
the projection to the field component \(\mathbb{K}_\kappa\).
Of course, invariance of \(I_F^\kappa\) for each \(\kappa\) will
imply invariance of the total \(I_F\).

Given \(a=\sum_g a_g g\in\Nov (G, N; R)\), let
\[\op{Tc}_{N; \Re}[a]:=\sum_{g\in G, N(g)\leq\Re}a_g g.\]

Choose an ordered basis \(e_{\cal P}\) for each \(\op{CF}({\cal
  V}_\lambda)\) as in \S\ref{def:torsion} and 4.3, thus
identifying \(\partial_F({\cal V}_\lambda)=(\partial_F({\cal
  V}_\lambda)_{ij})\) with a matrix with entries 
in \(\Omega_F\). \(\op{Tc}_{-\langle {\cal Y}\rangle; \Re}[\partial_F({\cal V}_\lambda)]\) will denote the
matrix with entries \(\op{Tc}_{-\langle {\cal Y}\rangle; \Re}[\partial_F({\cal V}_\lambda)_{ij}]\).
With the basis so fixed, \(\tau_F({\cal V}_\lambda)\) and \(I_F({\cal
  V}_\lambda)\) will now be viewed as elements in the Novikov ring
\(\Omega_F\). 

In contrast, let \(\partial_F({\cal V}_\lambda; \Re)\),
\(\eta_F({\cal V}_\lambda; \Re)\) be defined in the same way as
\(\partial_F\), \(\eta_F\) in \S2.2 using weight-truncated versions of
moduli:
\[\begin{split}
\partial_F({\cal V}_\lambda; \Re)_{x_\lambda, y_\lambda}&:=\sum_{A\in\mathfrak{H}}\chi\Big(\hat{\cm}_{P, \lambda}^0(\op{Lf}(x_\lambda),
A\cdot\op{Lf}(y_\lambda); \op{wt}_{-\langle{\cal Y}\rangle, e_{\cal P}}\leq \Re)\Big)\, [A]_{\ker\psi};\\
\eta_F({\cal V}_\lambda; \Re)&:=\sum_{A\in\ker\psi} \chi\Big(\hat{\cm}^0_{O, \lambda}(A); \op{wt}_{-\langle{\cal Y}\rangle, e_{\cal
    P}}\leq \Re)\Big)\, A.
\end{split}\]
\subsubsection{Invariance in regular range and left/right limits.}\label{cobord}
Below are some simple consequences of the smoothness and compactness
properties of the parameterized moduli spaces in (RHFS2, 3).
Let
\[\begin{split}\Lambda_{reg}&:=\Lambda\backslash (\Lambda_{hs}\,
  \cup\,\Lambda_{db});\\
\Lambda_{floer}&:=\{\lambda\, |\, \lambda\in \Lambda, \, \text{\({\cal
    V}_\lambda\) generates a Floer system}\}.\end{split}\]
\begin{lemma*}\begin{description}\itemsep -1pt
\item[(a)] \(\Lambda_{floer}\subset \Lambda_{reg}\subset \Lambda\) is Baire. 
\item[(b)] \(\partial _F({\cal V}_\lambda; \Re)\) is well defined and locally
constant for \(\lambda \) in \(\Lambda\backslash (\Lambda_{hs}^\Re\,
\cup\,\Lambda_{db})\). 
\(\eta_F({\cal V}_\lambda; \Re)\) is well defined and locally
constant for \(\lambda\) in \(\Lambda\backslash (\Lambda_{hs:ii}^\Re\,
\cup\,\Lambda_{db})\). 
\item[(c)] For \({\frak Q}=\partial_F, \eta_F, \tau_F, \zeta_F\) or
\(I_F\), \({\frak Q}({\cal V}_\lambda)\) is well
defined for all \(\lambda\in\Lambda_{reg}\).
In fact, for \({\frak Q}=\partial_F\) or \(\eta_F\), \({\frak Q}({\cal V}_\lambda; \Re)=\op{Tc}_{-\langle{\cal Y}\rangle; \Re}[{\frak Q}({\cal
  V}_\lambda)]\) for any
\(\lambda\in \Lambda_{reg}\). 
\item[(d)] The statements (i)--(iii) below hold for \({\frak Q}=\partial_F,
\eta_F, \zeta_F\). They also hold for \({\frak Q}^\kappa=\tau_F^\kappa, I_F^\kappa\) under the additional condition that
\begin{equation}\label{acyclic}
H_*(\op{CF}({\cal V}_\lambda)^\kappa)=0 \quad \forall\, \lambda\in \Lambda_{reg},
\end{equation}

\begin{description}\itemsep -1pt
\item[ (i)] For all \(\lambda\in \Lambda\), the following left and right
limits are well-defined.
\[{\frak Q}_{\lambda-}:=\lim_{\lambda'\in
  \Lambda_{reg}, \, \lambda'\nearrow\lambda}{\frak Q}({\cal
  V}_{\lambda'}), \quad {\frak Q}_{\lambda+}:=\lim_{\lambda'\in
  \Lambda_{reg}, \, \lambda'\searrow\lambda}{\frak Q}({\cal
  V}_{\lambda'}).\] 
\item[ (ii)] \({\frak Q}_{\lambda-}={\frak Q}_{\lambda+}={\frak Q}({\cal
  V}_\lambda)\) when \(\lambda\in \Lambda_{reg}\).
\item[ (iii)] For all sufficiently large \(\Re\), \(\op{Tc}_{-\langle{\cal Y}\rangle; \Re}[{\frak
  Q}_{\lambda-}]=\op{Tc}_{-\langle{\cal Y}\rangle  ; \Re}[{\frak Q}_{\lambda+}]\) is
locally constant on \(\Lambda\backslash (\Lambda_{db}\cup
\Lambda_{hs}^{2\Re})\).
\end{description}\end{description}
\end{lemma*}
\noindent {\it Sketch of proof.}
Given \({\bf x}, {\bf y}\in \aleph_\Lambda\), and \([a, b]\in
\Lambda\backslash (\Lambda_{db}\cup \Lambda_{hs}^\Re)\), 
(RHFS2) implies that 
\(\hat{\cm}_P^{[a,b], 1}({\bf x}, {\bf y}; \op{wt}_{-\langle{\cal
    Y}\rangle, e_{\cal P}}\leq \Re)\) is a
compact 1-manifold with boundary \(\hat{\cm}_{P,b}^{0}(x_b,
y_b; \op{wt}_{-\langle{\cal
    Y}\rangle, e_{\cal P}}\leq \Re)-\hat{\cm}_{P,a}^{0}(x_a,y_a; \op{wt}_{-\langle{\cal
    Y}\rangle, e_{\cal P}}\leq \Re)\). 
Namely, the parameterized moduli space forms a cobordism between
\(\hat{\cm}_{P, b}^0\) and \(\hat{\cm}_{P, a}^0\).

The fiber over \(\lambda\in [a, b]\),
\(\hat{\cm}_{P,\lambda}^{0}(x_\lambda,y_\lambda; \op{wt}_{-\langle{\cal
    Y}\rangle, e_{\cal P}}\leq \Re)=\Pi_\Lambda|_{\hat{\cm}_P^{[a,b],
    1}({\bf x}, {\bf y}; \op{wt}_{-\langle{\cal
    Y}\rangle, e_{\cal P}}\leq \Re)}^{-1}(\lambda)\), is
Zariski smooth and compact unless \(\lambda\) is a critical value of
\(\Pi_\Lambda\). (This is because \(\ker E_{u_\lambda}\subset\ker D\Pi_\lambda\),
where \(E_{u_\lambda}\) is the deformation operator of the fiber at
\(u_\lambda\)). By Sard's theorem, non-critical values form a Baire
set.
Claim (a) follows from this observation and its analog for
\(\hat{\cm}_O\).
 
On the other hand, notice that one does not need \(\lambda\in\Lambda_{floer}\) for \(\partial_F({\cal V}_\lambda; \Re)\),
\(\eta_F({\cal V}_\lambda; \Re)\) to be well-defined. Their definitions
involve only Euler characteristics of various moduli spaces
\(\hat{\cm}^0_\lambda\), and as long as there is a smooth compact
parameterized moduli space \(\Pi_\Lambda: \hat{\cm}^{[a,b],1}\to
[a,b]\) containing \(\hat{\cm}^0_\lambda\) as a fiber,
\(\chi(\hat{\cm}^0_\lambda)\) is well-defined and constant in \(\lambda\in[a,b]\).
Claims (b) and (c) then follow. 

To justify (d), one just needs to show that the terms
in 
\(\op{Tc}_{-\langle {\cal Y}\rangle ; \Re}[{\frak
  Q}_{\lambda}]\) (\(\lambda\in \Lambda_{reg}\)) 
only depend on Euler characteristics of 
moduli spaces with weight-truncation \(\op{wt}_{-\langle{\cal Y}\rangle,
  e_{\cal P}}\leq2\Re\). The local constancy of \(\op{Tc}_{-\langle {\cal Y}\rangle ; \Re}[{\frak
  Q}_{\lambda}]\) in \(\lambda\), and hence the existence of
left/right limits, would then follow from the same
cobordism argument used for (b).

For \({\frak Q}=\partial_F, \eta_F, \zeta_F\), this is clear from the
definitions; in fact, only flows with \(\op{wt}_{-\langle {\cal Y}\rangle, e_{\cal P}}\leq \Re\)
contribute. 

The story with \({\frak Q}=\tau_F^\kappa, I_F^\kappa\) is less straightforward and
requires the assumption that \(\Re\) is sufficiently large. We shall
concentrate on \({\frak Q}=\tau_F^\kappa\) since the case of \({\frak
  Q}=I_F^\kappa\) follows from the cases of \({\frak Q}=\tau_F^\kappa\) and \({\frak
  Q}=\zeta_F\).

First, note that by the compactness of \(\hat{\cm}_P^{\Lambda,
  1}(\op{wt}_{-\langle {\cal Y}\rangle, e_{\cal P}}\leq \Re)\), there
is a constant \(C>0\) such that 
\begin{equation}\label{wt-lbdd}
\op{wt}_{-\langle {\cal Y}\rangle,
  e_{\cal P}}(u_\lambda)\geq -C\quad\forall u_\lambda\in
\hat{\cm}_P^{\Lambda, 1}.
\end{equation}

Write \(\Lambda\backslash \Lambda_{db}=\coprod _k S_k\), where each
\(S_k\) is a path component of \(\Lambda\backslash \Lambda_{db}\).
The Floer groups \(\op{CF}({\cal V}_\lambda)\) are
independent of \(\lambda\) on each \(S_k\). By the acyclicity of
\((\op{CF}({\cal V}_\lambda)^\kappa, \partial _F({\cal V}_\lambda)^\kappa)\), there is a splitting
\(\op{CF}_{i}^\kappa=A_i^\kappa\oplus B_i^\kappa\), \(A_i^\kappa, B_i^\kappa\)  being each spanned by
elements in \({\cal P}\), such that the submatrices
\(\underline{\partial }_i^\kappa:=\partial_F^\kappa: A_i^\kappa\to B_{i-1}^\kappa\) are isomorphisms. (See \cite{hl2}
Lemma 2.7). By the local constancy of \(\partial_F({\cal V}_\lambda;
\Re)\) and (RHFS2*) (see also Lemmas 4.4.4, 4.4.5 below), \(A_i^\kappa\) and
\(B_i^\kappa\) may be chosen independently of \(\lambda\) in each
\(S_k\). Furthermore, \(\op{deg}(\underline{\partial }_i^\kappa)\) is well
defined and constant for each \(\lambda\) in \(S_k\); we denote this
value by \({\mathbb M}_{i, k}\). Let \({\mathbb
  M}:=\sup_k\sum_i|{\mathbb M}_{i,k}|\).

By (RHFS2), we may require \(\Re\) to be large enough so that the square
matrices \(\partial_F({\cal V}_\lambda;\Re)^\kappa: A_i^\kappa\to B_{i-1}^\kappa\) are
isomorphisms for all \(i\) and \(\lambda\in \Lambda\backslash
(\Lambda_{db}\cup \Lambda_{hs}^\Re)\). In addition, suppose 
\begin{equation}\label{Re-large}
\Re>2{\mathbb M}+C\cdot \#(\aleph_\Lambda).
\end{equation}
According to \cite{hl2} Lemma 2.7,
\(\tau_F^\kappa=\prod_i(\det\underline{\partial }_i^\kappa)^{(-1)^i}\).
Expanding, we see that for \(\lambda\in S_k\), 
the degree of any term in \(\tau_F({\cal V_\lambda})^\kappa\) can be
written as \(\sum_i(-1)^i{\mathbb M}_{i, k}\) plus a sum of {\em
  non-negative}
terms of the form \[\sum_{j=1}^{\dim A_i}\op{wt}_{-{\cal Y},
  e_{\cal P}}(u_j)-{\mathbb M}_{i, k},\] 
where \(u_j\) are flows from an
index i critical point to an index \(i-1\) critical point. 
The constraint that \(\deg \leq\Re\) for any term in \(\op{Tc}_{-\langle {\cal
    Y}\rangle ; \Re}[\tau_F({\cal V}_\lambda)^\kappa]\), together with (\ref{Re-large}) and
(\ref{wt-lbdd}) then imply that \(\op{wt}_{-{\cal Y},
  e_{\cal P}}(u_j)\leq 2\Re\) for any flow \(u_j\) contributing to \(\op{Tc}_{-\langle {\cal
    Y}\rangle ; \Re}[\tau_F({\cal V}_\lambda)^\kappa]\).
\hfill\(\Box\)

\subsubsection{At a Death-birth bifurcation.}
We shall focus on the case of a death, since the case of a birth
may be obtained from the case of a death by reversing the orientation
of \(\Lambda\). 

Suppose \(z_\lambda\in {\cal P}^{\Lambda, deg}\) is a death with
\(\ind_+(z_\lambda)=i\). Let \(\lambda_-, \lambda_+\in \Lambda_{floer}\) be
close to \(\lambda\) with \(\lambda_-<\lambda<\lambda_+\). By (RHFS1*), \({\cal P}_{\lambda_-}={\cal
  P}_{\lambda_+}\sqcup \{z^+_{\lambda_-}\}\sqcup
\{z^-_{\lambda_-}\}\), where \(\ind(z^+_{\lambda_-})=i\), \(\ind(z^-_{\lambda_-})=i-1\).
Thus, 
\[
\op{CF}_{i}({\cal V}_{\lambda_-})=\R\{z^+_{\lambda_-}\}\oplus\op{CF}_{
  i}({\cal V}_{\lambda_+});\quad  
\op{CF}_{i-1}({\cal V}_{\lambda_-})=\R\{z^-_{\lambda_-}\}\oplus\op{CF}_{
  i-1}({\cal V}_{\lambda_+}).
\]
Since the left and right limits \(\op{CF}_{\lambda\pm,
  i}=\op{CF}_{i}({\cal V}_{\lambda_\pm})\), they have similar
decompositions. 

The properties (RHFS2*, 3*, 4), especially those describing the
structures of the moduli spaces near \({\mathbb T}_{P, db}\) and 
\({\mathbb T}_{O, db}\), lead to:
\begin{lemma*}\begin{description}\itemsep -1pt
\item[(a)] In terms of the above decomposition of \(\op{CF}_{\lambda-, i}\), 
\(\op{CF}_{\lambda-, i-1}\),  the left limit \((\partial_{F,
  \lambda-})_i: \op{CF}_{\lambda-, i}\to \op{CF}_{\lambda-, i-1}\) is:
\[
(\partial_{F, \lambda-})_i=\left(\begin{array}{cc}1+b& v\\ w & N\end{array}\right),
\]
where \[\begin{split}
b &=\sum_{A\in \ker\psi} \chi(\hat{\cm}_{P,\lambda}^0(\op{Lf}(z_\lambda), \op{Lf}(z_\lambda))\, A=\sum_j b_j\\
b_j & =\sign (u_j)[u_j]_{\ker\psi} \quad 
\text{where \(u_j\in\hat{\cm}^0_P(z_\lambda, z_\lambda)\), \(u_j\neq z_\lambda\)},
\end{split}\]
$v=(v_{y_\lambda})$, $w= (w_{x_\lambda})$, \(N=(N_{x_\lambda,
  y_\lambda})\) are respectively a row vector, a column vector, and a matrix, with entries given by 
\[\begin{split}
v_{y_\lambda}&=\sum_{A\in\mathfrak{H}}\chi(\hat{\cm}_{P,
    \lambda}^0(\op{Lf}(z_\lambda), A\cdot \op{Lf}(y_\lambda))\, [A]_{\ker\psi},\\
w_{x_\lambda}&=\sum_{A\in \mathfrak{H}}\chi(\hat{\cm}_{P,
    \lambda}^0(\op{Lf}(x_\lambda), A \cdot\op{Lf}(z_\lambda))\, [A]_{\ker\psi},\\
N_{x_\lambda, y_\lambda}&=\sum_{A\in \mathfrak{H}}\chi(\hat{\cm}_{P,
    \lambda}^0(\op{Lf}(x_\lambda), A \cdot\op{Lf}(y_\lambda))\, [A]_{\ker\psi},
\end{split}\]
\(x_\lambda, y_\lambda\in {\cal P}_\lambda\)
being nondegenerate critical points with indices \(i\) and \(i-1\) respectively.
On the other hand, 
\begin{equation}\label{db-p}
(\partial_{F, \lambda+})_i=N+\sum_{k=0}^\infty
(-1)^{k+1}wb^kv,\end{equation}

\item[(b)] \((\partial_{F, \lambda-})_k=(\partial_{F, \lambda+})_k\) when
\(k\neq i, i\pm1\). When \(k=i+1\) or \(i-1\), \((\partial_{F,
  \lambda+})_k\) is obtained from \((\partial_{F, \lambda-})_k\) by
deleting the row corresponding to \(z_{\lambda -}\), or the column
corresponding to \(z_{\lambda +}\) respectively.

\item[(c)] We have
\begin{equation}
\label{db-eta}
\eta_{F, \lambda+}-\eta_{F, \lambda-}=\sum_{k=1}^\infty\sum_{i_1,
  \ldots, i_k=1}^\infty \frac{(-1)^{k+i+1}}{k}b_{i_1}\cdots b_{i_k}.
\end{equation}

\item[(d)] Hence, by elementary algebra
\[\begin{split}
\tau_{F, \lambda+}/\tau_{F, \lambda-}&=(1+b)^{(-1)^{i+1}}\\
\zeta_{F, \lambda+}/\zeta_{F, \lambda-}&=(1+b)^{(-1)^i}\\
I_{F, \lambda+}&=I_{F, \lambda-}.
\end{split}\]
\end{description}\end{lemma*}
\begin{remark*}
(a) By (RHFS2), the entries of \(b, v, w, N\) all take values in \(\Lambda_F\).

\noindent (b) Let \({\bf e}_{ij}\) denote a square matrix with the \(ij\)-th entry 1,
and 0 for all other entries. We refer to the conjugation action by an
elementary matrix of the form \(\op{Id}+\mathbf{e}_{ij}\) as an {\em
  elementary transformation}. By an {\em expansion}, we mean the direct
sum with the rank 1 identity matrix, and  a {\em collapse} refers to
the converse operation.
When \(\langle{\cal Y}\rangle=0\), \(b=0\), and the above lemma shows that in
this case the left and right limits \(\partial_{F,\lambda-}\),
\(\partial_{F, \lambda+}\) are related by expansion/collapse modulo 
elementary transformations.
\end{remark*}

\subsubsection{At a Type I handle-slide.}
Let \(u_\lambda\in\hat{\cm}_P^{\Lambda, 0}({\bf x},{\bf y})\), where \({\bf x}\neq {\bf y}\).

Let $T=\oplus_i T_i$,
where \(T_i: \op{CF}_i({\cal V}_\lambda)\to \op{CF}_i({\cal
  V}_\lambda)\) is
\[T_i=\begin{cases}Id &\text{when $i\neq \ind x_\lambda$}, \\
Id+\sign(u_\lambda) [u_\lambda]_{\ker\psi}{\bf e}_{x_\lambda, y_\lambda} &\text{when
  $i=\ind x_\lambda=\ind y_\lambda$}. \end{cases}\]
\begin{lemma*} In this case, 
\begin{description}\itemsep -1pt
\item[(a)] $\p_{F, \lambda+}=T\p_{F, \lambda-}T^{-1}$.
\item[(b)] \(\eta_{F, \lambda+}=\eta_{F, \lambda-}\). 
\item[(c)] Hence \({\frak Q}_{\lambda+}={\frak Q}_{\lambda-}\) for \({\frak Q}=\zeta_F,
\tau_F, I_F\).
\end{description}
\end{lemma*}
This is again a consequence of (RHFS2*, 3*, 4): (a) follows from the 
description of the structure of the moduli spaces near \({\mathbb
  T}_{P, hs-s}\), and (b) follows from the fact that 
in this case, \(\lambda\in \Lambda_{hs}\backslash
\Lambda_{hs:ii}\), and hence \(\Pi_\Lambda^{-1}(\lambda)\, \cap \partial \hat{\cm}_O^{\Lambda,
  1,+}=\emptyset\). 

\subsubsection{At a Type II handle-slide.} \label{neq-pert}
Up to ``first order'', this can be described similarly to
\S4.4.4 above. To understand the higher order contributions, we need to
further assume (NEP) below to apply Hutchings's argument.
Some preparation is required to state the assumption (NEP).

For convenience, we introduce the following weaker version of RHFS:
a homotopy of formal flows is {\em \(R\)-regular} with respect to a
weight truncation \(\op{wt}_{-\langle {\cal Y}\rangle, e_{\cal P}}\), 
if in addition to (RHFS1*), the statements in 
(RHFS2*, 3*, 4) are only required to hold for all weight-truncated versions \(\hat{\cm}_P^{0}(x_\lambda, y_\lambda; \op{wt}_{-\langle {\cal Y}\rangle ,
  e_{\cal P}}\leq\Re)\), \(\hat{\cm}_P^{S, k}({\bf x}, {\bf y}; \op{wt}_{-\langle {\cal Y}\rangle ,
  e_{\cal P}}\leq\Re)\), \(\hat{\cm}_O^{S, k}(\op{wt}_{-{\cal
    Y}, e_{\cal P}}\leq\Re)\) with \(\Re\leq R\).
Notice that we do {\em not} require an \(R\)-regular homotopy to be a CHFS.

A finite-cyclic cover of \(\lc\) is said to be {\em \({\frak
    H}\)-adapted} if it is a subcover of \(\tilde{\lc}\), and its
monodromy factors as 
\[\pi_m\circ \nu\circ \op{im}: \pi_1(\lc)\to \Z/m\Z,\]
where \(\nu \) is a homomorphism \(\nu: {\frak H}\to
\Z\), and \(\pi_m\) is the epimorphism \(\Z\to \Z/m\Z\).
We shall denote the \(m\)-cyclic cover of such monodromy by
\(\lc^{\nu, m}\). 
 
Let \(u_\lambda\in \hat{\cm}_P^{\Lambda, 0}(\op{Lf}({\bf x}),
A\cdot\op{Lf}({\bf x}))\) be a Type II handle-slide. 
Let \(\op{div}(A)\) denote the {\em divisibility} of \(A\in \ker\psi\subset{\frak H}\).

An \({\frak H}\)-adapted \(m\)-cyclic cover, \(\lc^{\nu, m}\), 
is said to be {\it \(u_\lambda\)-breaking} if  
\(\nu(A)=\op{div}(A)\), and 
\(m>1\) does not devide \(\op{div}(A)\). 
Indeed, in this case the end points of \(u_\lambda\) lift to different
points in \(\lc^{\nu, m}\).

We are now ready to state the assumption. 
Let \(\{(\lc, \mathfrak{H}, \ind, {\cal Y}_\lambda, {\cal V}_\lambda)\}_{\lambda\in \Lambda}\) be an RHFS,
and \(\op{wt}_{-\langle {\cal Y}\rangle, e_{\cal P}}\) be a weight filtration adapted to
this RHFS.
Let \(u_\lambda\in \hat{\cm}_P^{\Lambda, 0}(\op{Lf}({\bf x}),
A\cdot\op{Lf}({\bf x}))\) be a Type II handle-slide. 

\medbreak
\noindent {\rm\bf (NEP)} {\sc Existence of nonequivariant
  perturbations.} {\it We say that (NEP) holds for \(u_\lambda\) if,
given: 
\begin{itemize}\itemsep -1pt
\item an arbitrarily large 
\(\Re\in \R^+\), 
\item a small neighborhoods \(S_\lambda', S_\lambda\subset \Lambda\) 
about \(\lambda\), such that \(S_\lambda'\) is a proper sub-interval
of \(S_\lambda\), and \(\hat{\cm}_P^{S_\lambda, 0}(\op{wt}_{-\langle
  {\cal Y}\rangle, e_{\cal P}}\leq \Re)\) contains a unique element, \(u_\lambda\), 
\item an arbitrary \(u_\lambda\)-breaking, \({\frak
  H}\)-adapted finite-cyclic cover
\(\lc^{\nu,m}\), 
\end{itemize}
there is an \(\Re\)-regular homotopy of formal flows, \(\{{\cal V}_\lambda^+\,
|\, \lambda\in S_\lambda\}\) over \(\lc^{\nu, m}\), such that \({\cal V}_\lambda^+\) agrees with the
lift of \({\cal V}_\lambda\) to  \(\lc^{\nu,
  m}\) \(\forall \lambda\in \Lambda\backslash S_\lambda'\).}
\medbreak

\begin{lemma*}
Assuming that (NEP) holds for 
\(u_\lambda\in \hat{\cm}_P^{\Lambda, 0}(\op{Lf}({\bf x}),
A\cdot\op{Lf}({\bf x}))\). Then
\begin{description}\itemsep -1pt
\item[(a)] $\p_{F, \lambda+}=T \p_{F, \lambda-}T^{-1}$, where 
$T=\oplus_i T_i$, \(T_i: \op{CF}_i({\cal V}_\lambda)\to \op{CF}_i({\cal V}_\lambda)\)
being
\[T_i=\begin{cases}Id &\text{when $i\neq \ind x_\lambda$}, \\
Id+\chi{\bf e}_{x_\lambda, x_\lambda} &\text{when $i=\ind x_\lambda$}, \end{cases}\]
and 
\(\chi=\sign(u)[u]+\sum_{j=2}^\infty\chi_j[u]^j\in \Nov^+(\ker\psi,
-\langle{\cal Y}\rangle;\Z)\).
\item[(b)] \(\eta_{F, \lambda+}-\eta_{F, \lambda-}=(-1)^{\ind x}
\sign(u)[u]+\sum_{j=2}^{\infty}\chi'_j[u]^j=:\chi'\in \Nov^+(\ker\psi,
-\langle{\cal Y}\rangle;\Z)\).
\item[(c)] \(\tau_{F, \lambda+}=(1+\chi)^{(-1)^{\ind x+1}}\tau_{F, \lambda-}\).
\item[(d)] \(I_{F, \lambda+}=I_{F, \lambda-}\).
\end{description}
\end{lemma*}
\noindent{\it Sketch  of proof.} (a) follows from the
description of the structure of the parameterized moduli space of
broken trajectories near \({\mathbb T}_{P, hs-s}\) in (RHFS2c, 4), 
the property (NEP), and the proof of \cite{Hut} Lemma 3.7.

(b) is a consequence of the description of the structure of
parameterized moduli space of broken orbits near
\({\mathbb T}_{O, hs-s}\) in (RHFS3c, 4).

(c) follows from (a) by elementary algebra, and (b), (c) together imply 
\(I_{F, \lambda-}=cI_{F, \lambda+}\), where
\(c=1+\sum_{i=2}^{\infty}c_i[u]^i\).

Now by (NEP), the proof of 
\cite{Hut} Lemma 3.10 can be applied to verify (d).
\hfill\(\Box\)
\subsubsection{Conclusion.} Summarizing the above, we state:
\begin{proposition*}\label{prop:abs-inv}
Let \(\{{\cal V}_\lambda\}_{\lambda\in[1,2]}\) be an RHFS such that
(NEP) holds for all Type II handle-slide bifurcations. 
Then \(I_F({\cal V}_1)=I_F({\cal V}_2)\).
\end{proposition*}
More precisely, \(I_F({\cal V}_i)\) needs to be replaced by
\(i_{-\langle {\cal Y}\rangle}I_F({\cal V}_i)\) when \(\langle{\cal
  Y}_i\rangle=0\) for the above equality to make sense. 
\medbreak

\begin{proof}
We shall show that \(I_F({\cal V}_1)^\kappa=I_F({\cal V}_2)^\kappa\)
for each \(\kappa\).

It follows from the invariance of Floer homology that
either \(H_*(\op{CF}({\cal V}_\lambda)^\kappa)=0\) \(\forall
\lambda\in \Lambda_{reg}\), or \(H_*(\op{CF}({\cal V}_\lambda)^\kappa)\neq0\) \(\forall
\lambda\in \Lambda_{reg}\). (Cf. \cite{le-ono}. Alternatively, one may
derive it from the description of bifurcation behaviors in
\S4.4.3--5). In the latter case, \(I_F({\cal V}_1)^\kappa=I_F({\cal
  V}_2)^\kappa=0\) and we are done. 
Thus, we may assume \(H_*(\op{CF}({\cal V}_\lambda)^\kappa)=0\) \(\forall
\lambda\in \Lambda_{reg}\). The invariance of \(I_F^\kappa\) then follows
from Lemma 4.4.2 (d) and the comparison of left and right limits at
bifurcation points in Lemmas 4.4.3--5.
\end{proof}

With this Proposition in place, the proof of Theorem \ref{thm:gen-inv}
then consists in showing that the conditions of
Theorem \ref{thm:gen-inv} imply the existence of an RHFS with property
(NEP), that connects the Floer systems associated to \((J_1,X_1)\),
\((J_2, X_2)\). This is the goal
of the rest of this paper. 
First, to obtain (RHFS2d) we need to
generalize the structural results for \(\cm_P(x, y)\) in \S3.2 to case
where one or both of the critical points \(x, y\) are degenerate. This
is done in the next section. In section 6, we show that the conditions
of Theorem \ref{thm:gen-inv} imply the existence of certain CHFSs
satisfying all the properties (RHFS*) except for (RHFS2c, 3c, 4),
called ``admissible \((J, X)\)-homotopies''. It will be shown in part
II that the remaining properties (RHFS2c, 3c, 4), (NEP) 
can be derived from the
definition of admissible \((J, X)\)-homotopies. This will then conclude the
proof of Theorem \ref{thm:gen-inv}.

\section{Near a Minimally-Degenerate Critical Point $y$.}
We now return to the context of the specific version of Floer theory
described in \S3. 
This section contains detailed description of 
the behaviors of various flows in
the loop space near a minimally degenerate critical point \(y\): In
\S5.1 we give decay estimates of a flow in \(\cm_P\) ending at \(y\);
in \S5.2, these estimates are applied to obtain the structure theory
of moduli spaces such as \(\cm_P(x, y)\); in \S5.3 we assume that
\(y\) belongs to a family \(y\in {\cal P}^\Lambda\) for a homotopy of
Floer systems generated by a 1-parameter family \(\{(J_\lambda,
X_\lambda)\}_{\lambda\in\Lambda}\), and describe
the neighborhood of \(y\) in \({\cal P}^\Lambda\).
These results will be useful for verifying the property
(RHFS2d), and the estimates will be important for
the proofs of the gluing theorems in Part II.

\subsection{A Decay Estimate for the Flows.}

Following \cite{MMR}, we use center manifold theory on Hilbert manifolds
to describe the asymptotics of flows ending at \(y\) (Proposition
\ref{prop:decay} below). 

\mysubsubsection[Preparations.] Some preliminaries are required to state
the assumptions of the Proposition. 

First, notice that for any (possibly degenerate) \((y, [v])\in
\tilde{\cal P}\), the generalized Conley-Zehnder index 
\(\tilde{\ind}^\sigma(y, [v])\) is well-defined for any
non-zero \(\sigma\) with sufficiently small \(|\sigma|\). Furthermore,
it only depends on the sign of \(\sigma\).
\begin{definition*}
For any \((y, [v])\in\tilde{\cal P}\), 
\[\begin{cases}
\tilde{\ind}_+(y, [v]):=\tilde{\ind}^{\sigma}(y, [v]) &\text{ for a \(\sigma\in \R^-\)
with \(|\sigma|\ll1\);}\\
\tilde{\ind}_-(y, [v]):=\tilde{\ind}^{\sigma}(y, [v]) &\text{ for a \(\sigma\in \R^+\)
with \(|\sigma|\ll1\);}\end{cases}
\]
\end{definition*}
It is immediate from the definition that \[\tilde{\ind}_+(y, [v])-\tilde{\ind}_-(y,
[v])=\dim\cok A_y, \]
in particular, it is 0 if \(y\) is nondegenerate; 1 if \(y\) is
minimally degenerate.

It is also immediate that \(\tilde{\ind}_+, \tilde{\ind}_-\) are respectively the
upper semicontinuous and lower semicontinuous extension of the function
\(\tilde{\ind}\) on \(\tilde{\cal P}\backslash \tilde{\cal P}^{deg}\). Thus,
with the definitions 
\[\begin{split}
\ind_{\pm} &:=\tilde{\ind} \mod 2, \\
\op{gr}_+((x, [w]), (y,[v])) &:=\tilde{\ind}_+(x, [w])-\tilde{\ind}_-(y,
[v]),\\
\op{gr}_-((x, [w]), (y,[v])) &:=\tilde{\ind}_-(x, [w])-\tilde{\ind}_+(y,
[v]),
\end{split}\] 
(\ref{RHFS-1g}) indeed follows from (RHFS1). 

Next, given a small positive \(\epsilon\), let 
\[U_{\epsilon}\subset L^2_2(y^*K)=T_y\lc\] be an \(\epsilon\)-neighborhood of
\(0\). It corresponds to a small neighborhood of \(y\) in \(\lc\) via 
the exponential map \(\exp (y, \cdot)\). The vector field \({\cal
  V}_X\) on \(\lc\) pulls back to a vector field \(\Xi\) on \(U_{\epsilon}\):
\begin{equation}\label{Xi}
\Xi (\xi)=A_y(\xi)+n_y (\xi),
\end{equation}
\(n_y\) being the nonlinear term. We leave the reader the easy
task of finding the explicit formula for \(n_y\), and limit
ourselves to the following remark:
\(n_y\) is in general a function of both \(\xi\) and
\(\dot{\xi}\), but when \(J_{t}\) is integrable near \(y(t)\), it
depends on \(\xi\) only.
(An formula for \(n_y\) in an analogous situation may be found in \cite{floer.cpam}
Lemma 3.2. Note that in the second line of Floer's formula, \(X\) and \(Y\) should be interchanged. Moreover,
since we consider perturbation by symplectic vector fields, there will
be an additional term \(\Theta_N\), \(\Theta(\theta, y,
\xi):=(D_3e)^{-1}\check{\theta}_X(\theta, \exp(y, \xi))\) in the
notation of \cite{floer.cpam}).

We now state the conditions of the Proposition.
\mysubsubsection[The conditions.] 
Let \((x,[w]), (y, [v])\in \tilde{\cal P}(X)\) be either nondegenerate or minimally
degenerate, with 
\(
\op{gr}_+((x, [w]), (y, [v]))=i\leq1.\)
We want to describe the asymptotic behavior of any \(u\in
\cm_P((x,[w]), (y, [v]))\) when one of \(x\), \(y\) is minimally degenerate. (The
nondegenerate case is well-known). Without loss of generality, suppose
 \(y\) is minimally degenerate; the case when \(x\) is minimally
 degenerate is completely analogous. Assume additionally the following:
\begin{description}
\item[(1)] \(J\in {\cal J}\) is such that \(\forall t\in S^1\), \(J_t\) is
integrable in a small neighborhood of \(y(t)\in M\).

\item[(2)] Let \({\bf e}_y\) be a unit vector in \(\ker A_y\). Then 
\(\Pi_{\ker A_y}\nabla_{{\bf e}_y}\nabla_{{\bf e}_y}n_y(0)\neq 0\). 
In fact, we shall choose \({\bf e}_y\) to be the unique unit vector 
such that 
\begin{equation}\label{2ndDerv}
\Pi_{\ker A_y} \nabla _{{\bf e}_y}\nabla _{{\bf e}_y}n_y(0)=C_y{\bf
  e}_y \quad \text{for a real \(C_y>0\). }
\end{equation}

\item[(3)] For any sufficiently small \(\sigma>0\), \(\cm_P^{(\sigma, -\sigma)}((x,[w]),  (y, [v]))=\emptyset\).
\end{description}
Note that conditions (2) and (3) above are in fact generic: Condition
(2) holds for Hamiltonian perturbations in an open dense set of
\(V^2_\delta(J, X)\).
Condition (3) holds for a Baire set in \(V^k_\delta(J, X)\), \(k>1\), by the
weighted version of Proposition \ref{moduli}, since according to Lemma 3.2.3,
the expected dimension of \(\hat{\cm}_P^{(\sigma, -\sigma)}((x, [w]),
(y, [v]))\) is
\[\op{gr}_+((x, [w]), (y, [v]))-1<0.\]

We now state the main result in this subsection:
\begin{*proposition}\label{prop:decay}
Assuming \S5.1.2, let 
$u(s,t)$ be any element in $\cm_P((x,[w], (y, [v]))$, and write 
\[u(s)=\exp(y, \mu(s))\] 
for large \(s\). Then there exist 
constants \(C_i\), \(C_i'\), \(i=0, 1, 2\), such that: 
\[\begin{split}
C'_0/s\leq \|\mu(s)\|_{2,2,t}\leq C_0/s;\\
C'_1 s^{-2}\leq \|\partial_s\mu(s)\|_{2,1,t}\leq C_1s^{-2};\\
C'_2 s^{-3}\leq \|\partial_s^2\mu(s)\|_{2,t}\leq C_2 s^{-3}.
\end{split}\]
\end{*proposition}
In fact the argument below also yields estimates for higher derivatives, but we
shall only need them up to order 2.

The rest of this subsection is devoted to the proof of this Proposition.
\medbreak

First, recall that the usual Gromov compactness (such as in the proof of
\cite{floer.cpam} Theorem 1) implies:
\begin{*lemma}\label{gromov-close}
For any sufficiently small $\epsilon$, there is an $s_0=s_0(\epsilon)\gg1$
such that 
\begin{equation}\label{mu:nbhd}
\|\mu(s)\|_{2,2,t}+\|\partial_s\mu(s)\|_{2,1,t}+\|\partial_s^2\mu(s)\|_{2,t}\leq
\epsilon       \quad \text{for all $s\geq s_0$.}
\end{equation}
\end{*lemma}
In fact, \(\|\mu\|_{C_\epsilon([s_0, \infty)\times S^1)}\) is small. 

This enables us to focus on the neighborhood $U_\epsilon$. 

\mysubsubsection[Existence of center manifold.]
Let \(\ker A_y^{\perp}\) be the \(L^2\)-orthogonal complement of
\(\ker A_y\subset L^2(y^*K)\).
\begin{lemma*}
There is a smooth map 
\[
\zeta: \ker A_y \cap U_{\epsilon} \to \ker A_y^{\perp}\cap U_{\epsilon},
\]
such that its graph is a 1-submanifold of \(U_{\epsilon}\) tangent to
\(\ker A_y\) at \(0\), and tangent to \(\Xi\) elsewhere.
\end{lemma*}
This lemma follows from \cite{MMR} Theorem 5.1.1 and the smoothness of
\(J, X\), since \(A_y\) is
\(L^2\)-self adjoint with 1-dimensional 0-eigenspace. The graph of
\(\zeta\) above is a {\it ``center manifold''} of the flow generated
by \({\cal V}_X\); 
It can be made unique by extending \(\Xi\) linearly outside
\(U_{\epsilon}\subset L^2_2(y^*K)\). (See \cite{MMR}.) 

\paragraph{A simple example.} Let $M=\R^2=\{(x,y)| x,y\in \R\}$,
\(\omega=dx\wedge dy\). Suppose \(J\)  and \(H\) are both
\(t\)-independent; $J=J_0$ is the standard complex structure, and
$H_{\lambda}(t, x, y)=x^2+\lambda y^2+y^3$. For the Floer system
generated by \((J, \chi_{H_0})\), 
the center manifold is simply the $y$-axis, and the flow
on the center manifold (to 0) is easily seen to decay as $C/s$ by a straightforward
computation.

\mysubsubsection[Estimates for flows on the center manifold.]
We now show that in general, the flow on the center manifold 
has the same decay behavior exhibited by the example.

To solve for the flow on the center
manifold, it suffices to solve for the flow on its projection to $\ker
A_y$.

Let \(b: \R\to \ker A_y\).
If $b(s)+\zeta(b(s))$ is a flow of \(\Xi\), then $b(s)$ satisfies: 
\begin{equation}\label{Nflow}
\frac{db}{ds}+\Pi_{\ker A_y} n_y(b+\zeta(b))=0,
\end{equation}
where \(\Pi_{\ker A_y}\) denotes the \(L^2\)-orthogonal projection
from \(L^2(y^*K)\) to \(\ker A_{y}\).
\begin{lemma*}
Assuming (\ref{2ndDerv}) and (\ref{mu:nbhd}), then 
any non-constant flow on the center manifold \(b(s)+\zeta(b(s))\)
satisfies the following estimate:
\[
\bar{C}'_k/s^{k+1}\leq \|\partial_s^k(b(s)+\zeta(b(s))\|_{2,t} \leq \bar{C}_k/s^{k+1}
\]
for some  positive constants \(\bar{C}_k,\bar{C}'_k\),  \(k=0,1,2\).
\end{lemma*}
\begin{proof}
First observe the useful fact that on \(\ker A_y\), all the Sobolev
norms are commensurate.
Let $b(s)=\beta (s){\bf e}_y$ for an \(\R\)-valued function
\(\beta(s)\). 
We have the following key inequalities: for all sufficiently small \(b\),
\begin{equation}\label{keyIneq}
C'\beta^2\geq|\langle {\bf e}_y,  n_y(b+\zeta(b))\rangle_{2,t}|\geq C\beta^2.
\end{equation}
The left inequality is straightforward from the form of \(n_y\); the
right inequality is also straightforward from the definitions of
\(n_y\) and \(\zeta\) by (\ref{2ndDerv}).

By the right inequality in (\ref{keyIneq}) and (\ref{Nflow}) we have
\[\frac{d\beta}{ds}\leq -C\beta^2\quad \text{for a positive constant $C$}.\] 
Suppose \(\beta\) is a nontrivial solution. Then \(\beta\) is nowhere
zero since it satisfies an ODE. We may thus divide
both sides of the above expression by $\beta^2$ and integrate over
$s$, which gives 
$\beta(s)\leq 1/(Cs+B)$ if $\beta(s_0)\geq 0$, and if $\beta(s_0)<0$, $\beta(s)$ goes to
$-\infty$ as $s$ increases, contradicting our assumption.
Similarly, using the left inequality of (\ref{keyIneq}), we obtain
$\beta(s)\geq 1/(C's+B')$. (\(B, B'\) are some constants.) These estimates together with the fact that
$\zeta$ vanishes up to first order establish the \(k=0\) case 
of the lemma.

The \(k=1\) case of the lemma follows from 
(\ref{Nflow}) and (\ref{keyIneq}).
Taking the $s$-derivative of (\ref{Nflow}) and using the \(k=0,1\)
cases of the lemma and (\ref{2ndDerv}), we get the 
\(k=2\) case of the lemma.
\end{proof}

\mysubsubsection[Approximation of \(u(s)\) to the center manifold.]
Next, we show that other flows to
$y$ approach the center manifold exponentially.
\begin{lemma*}\label{c:decay}
Let $u(s)=\exp (y,\mu (s))$ be a flow (i.e. a solution of (\ref{pcr}))
to $y$ satisfying (\ref{mu:nbhd}). Write
\[
\mu(s)=a(s)+c(s),
\] 
where \(a(s):=b(s)+\zeta(b(s))\) is a path on the center manifold,
\(b(s)\in \ker A_y\), and 
$c(s)\in \ker A_y^{\perp}$.
Then 
there exists positive constants \(C\), \(\gamma\) such that 
\[
\|c(s)\|_{2,2,t}+\|\partial_s c(s)\|_{2,1,t}+\|\partial _s^2
c(s)\|_{2,t}\leq Ce^{-\gamma s} \quad \forall s\geq s_0.
\] 
\end{lemma*} 

\begin{proof}
Observe that since \(a(s)\) is in the center manifold, the
vector field \(\Xi\) is parallel to \(\partial_s a\).  Thus, 
the flow equation (\ref{pcr}) becomes:
\begin{gather}\label{b:flow}
-\frac{da}{ds}=(1+\nabla_{b}\zeta)\Pi_{\ker A_y} n_y(a+c);\\
\label{c:flow}
-\frac{dc}{ds}=A_yc+(1-\Pi_{\ker A_y}-\nabla_{b}\zeta\Pi_{\ker A_y})((n_y(a+c)-n_y(a)).
\end{gather}
Let's further decompose $c$ into 
\[
c=c_+ +c_-,
\] 
where $c_+$ is in the completion of the direct sum of eigenspaces of \(A_y\) with
positive eigenvalues; similarly for \(c_-\). 

Note that the estimate 
\begin{equation}\label{0th}
\|c(s)\|_{2,t}\leq Ce^{-\gamma s}
\end{equation} 
follows from the routine
arguments in \cite{MMR} pp.91-92 and the next lemma. 
The constant \(\gamma\) is a positive number slightly
smaller than the number \(\nu_-\)  below.
\begin{sublemma*}\label{c-pm:est}
Assume \S5.1.2.
Let $\nu_+$ be the minimal positive eigenvalue of $-A_y$, and let
\(-\nu_-\) be the maximal negative eigenvalue. Then there exist positive constants \(\epsilon_{\pm}\), with
$\epsilon_{\pm}\ll\nu_{\pm}$, such that the following hold:
\begin{gather*}
\frac{d}{ds}\|c_+\|_{2,t}\geq \nu_+ \|c_+\|_{2,t}-\epsilon_+ \|c\|_{2,t};\\
\frac{d}{ds}\|c_-\|_{2,t}\leq -\nu_- \|c_-\|_{2,t}+\epsilon_- \|c\|_{2,t},
\end{gather*}
\end{sublemma*}
\begin{proof}
The first inequality follows by 
taking the inner product of (\ref{c:flow}) with $c_+$,
using the facts that $C^1\subset L^2_2$ in 1-dimension, that
\(\|a\|_{2,2,t}+\|c\|_{2,2,t}\) is small, and that \(n_y\) is a
function of \(\xi\) only under the assumption \S5.1.2 (1).

The second inequality is similar.
\end{proof}

To get estimates on higher derivatives of \(c\), we perform elliptic
bootstrapping in the following way. 
Re-arranging (\ref{c:flow}) and writing 
\[
A_y=J_t\partial_t+\Upsilon_y,
\]
we have
\begin{equation}\label{c0}
(\partial_s+J_t\partial_t)c=-\Upsilon_y c
+(1-\Pi_{\ker A_y}-\nabla_{b}\zeta\Pi_{\ker A_y})(n_y(a+c)-n_y(a)).
\end{equation}
Let \(\beta_{\underline{s}, j}: \R\to\R\) be a smooth cutoff function that is
supported on \((\underline{s}-2j, \underline{s}+2j)\), and is 1 over \([\underline{s}-j,
\underline{s}+j]\). Let \(\Theta_{\underline{s}, j}:=[\underline{s}-j, \underline{s}+j]\times S^1\). Multiplying the above expression by
\(\beta_{\underline{s}, 1}(s)\) and
squaring both sides, we obtain
\[\begin{split}
\|\beta_{\underline{s},1}\partial_s c\|^{2}_{L^2(\Theta_{\underline{s}, 2})}+\|\beta_{\underline{s},1}\partial_t
c\|^{2}_{L^2(\Theta_{\underline{s}, 2})}
&\leq C\|\beta_{\underline{s},1}c\|^2_{L^2(\Theta_{\underline{s}, 2})}+C'\|\beta'_{\underline{s},1}c\|_{L^2(\Theta_{\underline{s}, 2})}\|\beta_{\underline{s},1}\partial_t c\|_{L^2(\Theta_{\underline{s}, 2})}\\
&\leq C''\|c\|^2_{L^2(\Theta_{\underline{s}, 2})}+1/2\|\beta_{\underline{s}}\partial_t c\|^2_{L^2(\Theta_{\underline{s}, 2})}
\end{split}\]
by the triangle inequality. Rearranging, we get:
\begin{equation}\label{c:Dbound}
\|\partial_s c\|^{2}_{L^2(\Theta_{\underline{s}, 1})}+\|\partial_t
c\|^{2}_{L^2(\Theta_{\underline{s}, 1})}
\leq C_3 \|c\|^2_{L^2(\Theta_{\underline{s}, 2})}\leq C_4 e^{-2\gamma\underline{s}}
\end{equation}
by integrating (\ref{0th}).

Now, we may iterate this argument, replacing (\ref{c0}) by s- or
t-derivatives of the equation, and using (\ref{mu:nbhd}) repeatedly.
The generalizaton of (\ref{c:Dbound}) boundes
\(\|\nabla^kc\|^2_{L^2(\Theta_{\underline{s}, 2^{1-k}})}\) in terms of
\(\|\nabla^jc\|^2_{L^2(\Theta_{\underline{s}, 2^{1-j}})}\) for \(j=0,
1, \ldots, k-1\), and hence by \(C'_k e^{-2\gamma \underline{s}}\).

On the other hand, by Sobolev embedding 
\[
\|\eta(\underline{s})\|_{2,t}\leq C_\varepsilon \|\eta\|_{L^2_1(\Theta_{\underline{s}, \varepsilon})};
\]
so the above estimates of \(\|\nabla^kc\|_{L^2(\Theta_{\underline{s},
    \varepsilon})}\) gives the other estimates in Lemma \ref{c:decay}.
\end{proof}

\mysubsubsection[Approximation of \(a(s)\) to a flow on the center
manifold. ]\label{518}
In Lemma \ref{c:decay}, we showed that the flow \(u(s)\) approximates a
path \(a(s)\) in the center manifold. However, \(a(s)\) might not be a
flow itself. Nevertheless, 
the next lemma shows that \(a(s)\) approximates
the flow on the center manifold described in Lemma 5.1.6 exponentially fast.

\begin{lemma*}
Let $a(s)$ be a path in the center manifold in \(U_{\epsilon}\)
satisfying 
\begin{equation}\label{smallE}
\|\partial_sa(s)+A_y(a(s))+n_y(a(s))\|_{2,1,t}+\|\partial_s[\partial_sa(s)+A_y(a(s))+n_y(a(s))]\|_{2,t}\leq
Ce^{-\gamma s}
\end{equation}
for all \(s \geq s_0\).
Then there is a unique
flow $z(s)$ on the center manifold  such that
\[
\|z(s)-a(s)\|_{2,2,t}+\|\partial_s(z(s)-a(s))\|_{2,1,t}+\|\partial_s^2(z(s)-a(s))\|_{2,t}\leq
C'e^{-\gamma s} \quad \forall s \geq s_0.
\]
\end{lemma*}
\begin{proof}
The estimate for \(\|z(s)-a(s)\|_{2,t}\) is proved in \cite{MMR}
Lemma 5.3.1 by a simple contraction mapping theorem argument. In fact,
it is proved there that \(\|\Pi_{\ker A_y}(z(s)-a(s))\|_{2,t}\leq
C_1e^{-\gamma s}\). However, as we noted before, since on the
1-dimensional space \(\ker A_y\)
the Sobolev norms are commensurate, this implies \(\|\Pi_{\ker A_y}(z(s)-a(s))\|_{2,2,t}\leq
C_2e^{-\gamma s}\), which in turn implies the exponential decay of
\(\|z(s)-a(s)\|_{2,2,t}\) by the uniform bounds on 
\(\nabla\zeta\) in \(U_{\epsilon}\). Now the estimates for
\(\|\partial_s(z(s)-a(s))\|_{2,1,t},\|\partial_s^2(z(s)-a(s))\|_{2,t}\)
can be derived from this iteratively via the condition (\ref{smallE})
and the fact that \(\partial_s z(s)+A_y z(s)+n_y(z(s))=0\).
\end{proof}

\mysubsubsection[Concluding the proof of Proposition \ref{prop:decay}.]
Following the notation of Lemma \ref{c:decay}, 
if \(\exp (y, a(s)+c(s))\) is a flow to \(y\), then $a(s)$
satisfies (\ref{smallE}) by the estimates on \(c(s)\) in Lemma
\ref{c:decay}. So we may apply Lemma \ref{518}. On the other hand, the flow
\(z(s)\) obtained in Lemma \ref{518} cannot be constant, because
otherwise by Lemmas \ref{c:decay} and \ref{518}, \(u(s)\)
approaches \(y\) exponentially. This would then contradict the
assumption \S5.1.2 (3). Proposition \ref{prop:decay} now follows from
the combination of Lemmas 5.1.6, \ref{c:decay}, and \ref{518}.

\subsection{Fredholm Theory and Structure Theorem for $\cm_P(x, y)$.}

In this subsection we prove the structure theorem for moduli spaces
\(\cm_P(x, y)\) when one or both of \(x, y\) are minimally degenerate,
analogous to the standard Proposition \ref{moduli}. 
(See Proposition \ref{y-moduli}
below). 

Without loss of generality we again assume that \(y\) is minimally degenerate.

Because of the degeneracy of \(y\), the usual \(L^p\) Fredholm theory
fails in this situaion. One {\em
  does} have a good Fredholm theory by working with the exponentially
weighted Sobolev spaces introduced in \S3.2.3: in particular, Lemma
3.2.3 works in this situation for a non-zero weight \(\sigma_2\neq 0\).
However, these Fredholm theories are unsuitable for the purpose of 
describing \(\cm_P(x, y)\) for the following reasons: 
Let \(\sigma\in \R\) with \(|\sigma|\ll1\). If \(\sigma>0\), one may work with the \(L^p_{k:(\sigma,
  -\sigma)}\) norms to describe the moduli space \(\cm_P^{(\sigma,
  -\sigma)}(x, y)=\cm_P\cap {\cal B}_P^{(\sigma, -\sigma)}(x,
y)\). However, \(\cm_P^{(\sigma,
  -\sigma)}(x, y)\) does not contain the whole \(\cm_P(x, y)\), since
the former contains only exponentially decaying flows, while the
latter contains flows that are polynomially decaying, as seen in
\S5.1.  On the other hand, one cannot take \(\sigma<0\) either,
because \(L^p_{k:(\sigma, -\sigma)}\) would be too large for an
appropriate local description of the moduli space: for example, 
the Kuranishi map would not be suitably bounded with
respect to these norms (i.e. estimates for nonlinear terms of the type of 
(II.3) fail). 

The decay estimates in \S5.1 suggest that polynomially-weighted
spaces would be the natural alternative. Fredholm theories of such
spaces is relatively uncommon in the literature, see however 
\cite{floer.cpam, L} for similar theories.  
We shall follow \cite{floer.cpam} section 4 closely.

\mysubsubsection[LT-decomposition and polynomial weights.]
Let $u\in \cm_P^i(x,y)$, 
and write \(u(s)=\exp(y, \mu(s))\) as before. 
Assume without loss of generality that $x$ is nondegenerate and 
$y$ is minimally degenerate.

We saw in last subsection that as $s\to\infty$, $\mu '(s)/\|\mu'(s)\|_{2,t}$ converges
to the degenerate direction $\ker A_{y}$.
We shall call it the ``longitudinal direction''. Given an section of $u^*K$, 
$\xi$, we denote its longitudinal and transversal component by $\xi_L$ and
$\xi_T$ respectively. Namely,
\begin{equation}\label{decom}
\xi=\xi_T+\xi_L,
\end{equation} 
where
$$
\xi_L(s):=\|u'(s)\|_{2,t}^{-2}\langle u'(s), \xi(s)\rangle_{2,t}u'(s).
$$
Roughly speaking, in the transversal direction, things go the same way
as in the nondegenerate case, while in the longitudinal direction novelties
arise. Let \(\Pi_L, \Pi_T\) denote the \(L^2\)-orthogonal projection
to the longitudinal/transversal directions respectively.

Define the weight 
$$
\sigma_u(s)=\left\{ \begin{array}{cc}
\|u'(0)\|_{2,t}^{-1} & \mbox{for $s\leq 0$.}\\
\|u'(s)\|_{2,t}^{-1} & \mbox{for $s\geq 0$.}\end{array}\right.
$$
From the previous subsection we know that $C'/s^2\geq\sigma_u\geq C/s^2$ for large $s$. 
\begin{definition*}
Let \(u\in \cm_P^i(x, y)\) as above. We define the following norms on
$C^{\infty}_0(u^*K)$: 
\begin{gather*}
\|\xi\|_{W_u}:=\|\sigma_u^{1/2}\xi\|_{p,1}+\|\sigma_u\xi '_L\|_p;\\
\|\xi\|_{L_u}:=\|\sigma_u^{1/2}\xi\|_p+\|\sigma_u \xi_L\|_p.
\end{gather*}
Let $W_u=W_u(u^*K)$, $L_u=L_u(u^*K)$ denote the completion of $C^{\infty}_0$ with respect to 
the above two norms respectively.
\end{definition*}
\begin{*proposition}\label{prop:Fred} {\rm\bf (Fredholmness)}
Let \(u\in \cm_P^i(x,y)\), where \(i\leq 1\), \(y\) is minimally degenerate, and
\(x\) is nondegenerate.
Then the operator $E_u: W_u\to L_u$ is bounded and
Fredholm of index \(i\).
\end{*proposition}

With the obvious modification of the defnition of \(W_u\) and \(L_u\),
the statement of this Proposition also holds when \(x\) is minimally
degenerate, and \(y\) is nondegenerate or minimally degenerate. 
\begin{notation*}
In the case when at least one of the two end critical points is
minimally degenerate, we shall always reserve the notation \(E_u\) for
the deformation operator between \(W_u\) and \(L_u\). The deformation
operator between exponentially weighted Sobolev spaces will come with
a superscript, of the form \(E_u^{(\sigma_1, \sigma_2)}\)
(cf. \S3.2.3).
\end{notation*}
The proof of this Proposition occupies \S5.2.3--5.2.5 below.
Since it is not very different from the proof of
\cite{floer.cpam} Theorem 4a, we shall only provide details where
modification or clarification is needed.

\mysubsubsection[Boundedness of \(E_u\).]
Via (\ref{decom}), \(E_u\) decomposes into four components:
\[E_{LT}:=\Pi_TE_u\Pi_L,\] and similarly \(E_{TT}\), \(E_{TL}\), \(E_{LL}\). 
The boundedness is straightforward to check by looking at each of these components:
The boundedness of \(E_{TT}\) is the same as the usual (nondegenerate)
case; the boundedness of \(E_{LT}, E_{TL}, E_{LL}\) follows respectively
from the observations (\ref{ob:a}), (\ref{ob:b}), (\ref{ob:c}) below.

\mysubsubsection[Fredholmness of \(E_u\).]
Simple linear algebra shows that if the diagonal components
\(E_{TT}, E_{LL}\) are Fredholm, and one of the off-diagonal
components vanishes while the other is bounded, then \(E_u\) is
Fredholm and 
\[
\ind E_u=\ind E_{TT}+\ind E_{LL}.
\]
Our task is thus to verify each of the above.

\begin{equation}\label{ob:a}
\hbox to 10pt{\hss\hbox to \hsize{\sc $\bullet \, E_{LT}$ is zero:\hfil}\hss}
\end{equation} 
This follows immediately from the fact that $E_u u'=0$.
\begin{equation}\label{ob:b}
\hbox to 10pt{\hss\hbox to \hsize{\sc $\bullet \, E_{TL}$ is bounded:\hfil}\hss}
\end{equation}
 Following Floer's computation (4.18) of \cite{floer.cpam} (but replacing \(J\partial_t\)
there by \(A_{u(s)}\)), when \(s\) is large
\[
(E_u \xi_T)_L (s)=-2\sigma_u^2(s)\langle u''_T(s), \xi_T(s)\rangle_{2,t}u'(s).
\]
$\|u''_T(s)\|_{2,t}$ is estimated as follows.
Write \(u(s)=\exp (y, \mu(s))\), where \[\mu(s)=b(s)+\zeta(b(s))+c(s)\] as in the previous subsection, we see
that 
\begin{gather*}
\mu '=(1+\nabla\zeta)b'+c'\\
\mu ''=b''+\nabla_{b''}\zeta+\nabla_{b'}\nabla_{b'}\zeta+c''.
\end{gather*}
We know from Lemma \ref{c:decay} that $\|c'(s)\|_{2,t}, \|c''(s)\|_{2,t}$ decays exponentially. 
Thus, modulo exponentially decaying terms, $\|u_T''(s)\|_{2,t}$ is bounded by
\[
C\|b'(s)\|^2_{2,1,t}+\tilde{C}\|b(s)\|_{2,1,t}\|b''(s)\|_{2,t}\leq
C's^{-4}.\] 
Thus, for large \(s\)
\begin{equation}\label{est:E-TL}
\begin{split}
\|(E_u \xi_T)_L\|_{L_u}& =\|\sigma_u (E_u \xi_T)_L\|_p\\
&\leq 2\|\sigma_u^2\langle u'',\xi_T\rangle_{2,t}(s)\|_p\\
&\leq C\|\,\|\xi_T\|_{2,t}(s)\|_{p,s}\\
&\leq C'\|\xi_T\|_p,
\end{split}
\end{equation}
where the subscripts $s$ and $t$ in the third line above indicate which variable
the norm is respect to; and the final $L^p$ norm is with respect to both
variables $s$ and $t$. In the last step we used the fact that on $S^1$ the
$L^2$-norm is bounded by the $L^p$ norm.
\medbreak

\noindent {\sc \(\bullet\) \(E_{TT}\) and \(E_{LL}\) are Fredholm:}
By a typical excision argument (see e.g. \cite{floer.cpam}), it suffices
to consider only \(\xi\) supported on \(\{(s,t): s\in (R, \infty)\}\)
for a large \(R\), and to replace \(u\) by $v$,
which is a small loop in $\lc$ from $y$ 
to itself, such that $v$ agrees with $u$ for large $s$.
We assume further that \(\|v'(s)\|_{2,1,t}\neq0\) and  
that $v$ is symmetric about $s=0$ for simplicity.

The Fredholmness of \(E_{TT}\) is shown in \cite{floer.cpam} in a way
similar to the nondegenerate case. To see the Fredholmness of
\(E_{LL}\), write $\xi_L=: f\sigma_v v'$, where \(f\) is an
\(\R\)-valued function on \(s\). Then $E_{LL}$ is
equivalent to the operator $\frac{d}{ds}+\sigma_v'\sigma_v^{-1}$ with
domain and range being the completion of \(C_0^{\infty}(\R)\) with
respect to the norms \(N_1\), \(N_2\) respectively:
\[
\|f\|_{N_1}:=\|\sigma_v^{1/2}f\|_p+\|\sigma_v f'\|_p; \quad
\|f\|_{N_2}:=\|\sigma_v f\|_p.
\]
Conjugating this operator by $\sigma_v$ makes it the
operator $\frac{d}{ds}$ between the completion with respect to the
\(N\) norm and \(L^p\), where
\[
\|g\|_{N}:=\|\sigma_v^{-1/2}g\|_p+\|g'\|_p.
\]
\begin{equation}\label{ob:c}
\hbox to 10pt{\hss\hbox to \hsize{This operator $\frac{d}{ds}: N\to
    L^p$ is obviously bounded,\hfil}\hss}
\end{equation}
and has a 1-dimensional kernel, namely, the
space of constant functions. 
Let \(N_0\subset N\) be the codimension 1 subspace
\[
N_0:=\{g(s)\, |\, g(s)\in N, \, g(0)=0\}.
\]
We claim that the restriction \(d/ds|_{N_0}\) has a bounded
inverse. This would then imply the Fredholmness of \(d/ds: N\to L^p\),
and hence \(E_{LL}\).

A natural candidate for its right-inverse is the integral:
\begin{equation}\label{int-inverse}
(\frac{d}{ds})^{-1}q (s):=\int_{0}^s q(\tilde{s})\, d\tilde{s}.
\end{equation}
To see that this is well-defined and bounded, notice the following
estimate for compactly supported \(L^p_{1, loc}\) functions on \((0, \infty)\):
\begin{equation}\label{eq:int-est}
\begin{split}
\|g'\|_p+\|\sigma^{-1/2}_vg\|_p &\leq \|g'\|_p+C\|s^{-1}g\|_p\\
&\leq \|s\frac{d}{ds} (s^{-1}g)\|_p+(C+1)\|s^{-1}g\|_p\\
&\leq C'\|(s\frac{d}{ds} +1)(s^{-1}g)\|_p\\
&=C'\|\frac{d}{ds}g\|_p.
\end{split}
\end{equation}
The third step above follows from the fact that the integral kernel of
the operator \(s\frac{d}{ds}+1\) decays as \(s^{-1}\) for large \(s\).

Let \(g(s)=\chi_R(s)\tilde{g}(s)\), where \(\tilde{g}(s)=\int_{0}^s q(\tilde{s})\, d\tilde{s}\), and
$\chi_R(s)=:\chi(s/R)$, $\chi$ being a smooth cutoff function supported on
$(-2,2)$, with \(\chi(s)=1\) for \(s\in [-1, 1]\). Then the above
estimate (\ref{eq:int-est}) yields
\[\begin{split}
\|\chi_R\tilde{g}'\|_{L^p(\R^+)}+\|\sigma^{-1/2}_v\chi_R\tilde{g}\|_{L^p(\R^+)}
& \leq C'\|\chi_Rq\|_{L^p(\R^+)}+C''R^{-1}\|\chi'(s/R)\tilde{g}(s)\|_{L^p(\R^+)}\\& \leq C'\|\chi_Rq\|_{L^p(\R^+)}+C_2\|q\|_{L^p([0, 2R])}.
\end{split}\]
(Cf. Lemma II.3.3.3(c) for the second step). 
Taking \(R\to \infty\) and combining with the analogous estimates on the
other half of the real line, we see that \(\tilde{g}\in N_0\)
when \(q\in L^p\); in fact, 
\[
\|\tilde{g}\|_N\leq C_n\|q\|_p.
\]
Thus, the expression (\ref{int-inverse}) does give
a well-defined, bounded inverse of \(d/ds: N_0\to L^p\).

We may now conclude that \(E_{LL}\) is Fredholm, and hence so is
$E_u$. 

\mysubsubsection[Computation of the index.]
Let \(\sigma\) be a small positive number. We claim that
\[\ind E_u=\ind E_u^{(-\sigma, \sigma)}.\]
The index computation of \(E_u\) stated in Proposition \ref{prop:Fred} would then
follow from Lemma 3.2.3. 

Since \(L^p_{1:(-\sigma, \sigma)}\) and \(L^p_{(-\sigma, \sigma)}\) contain $W_u$ and $L_u$ respectively,
it suffices to show that $\ker (E^{(-\sigma, \sigma)}_u)\subset W_u$, and 
$\cok E_u\subset (L_{(-\sigma, \sigma)}^p)^*$. We shall only show the first, since
the second is similar. This boils down to estimating the solution 
$\xi$ of $E_u^{(-\sigma, \sigma)}\xi=0$ where $s$ is large.
If $\xi$ solves $E_u^{(-\sigma, \sigma)} \xi=0$, then 
\begin{gather}\label{eq:1}
(E_u^{(-\sigma, \sigma)}\xi_T)_T=0;\\ 
(E_u^{(-\sigma, \sigma)}\xi_L)_L+(E_u^{(-\sigma, \sigma)}\xi_T)_L=0, \label{eq:2}
\end{gather}
since as noted before, $(E_u^{(-\sigma, \sigma)}\xi_L)_T=0$. 

Using these, we estimate $\|\xi_T(s)\|_{2,t}$ like (4.23) of \cite{floer.cpam} as follows.
\[
\begin{split}
&\frac{d^2}{ds^2}\langle\xi_T, \xi_T\rangle_{2,t} \\
&\quad =2\frac{d}{ds}\langle \xi'_T,\xi_T\rangle_{2,t}\\
&\quad=-2\frac{d}{ds}\langle A_{u(s)}\xi_T,\xi_T\rangle_{2,t}\\
&\quad=-4\langle A_{u(s)}\xi_T,\xi_T'\rangle_{2,t}-2\langle
A'_{u(s)}\xi_T,\xi_T\rangle_{2,t}\\
&\quad=4\|A_{u(s)}\xi_T\|_{2,t}^2+4\langle(A_{u(s)}\xi_T)_L, (E_u^{(-\sigma, \sigma)}\xi_L)_L\rangle_{2,t}-2\langle A'_{u(s)}\xi_T,\xi_T\rangle_{2,t}\\
\end{split}
\]
When $s$ is large (so that $u(s)$ is close to $y$), the terms in the 
last line above can be estimated as follows. For the first term, we
have the standard estimate:
\[
4\|A_{u(s)}\xi_T\|_{2,t}^2\geq 4\nu^2\|\xi_T(s)\|_{2,t}^2,
\]
where $\nu$ is a positive number slightly smaller than $\min (|\mu_i|)$, and
$\mu_i$ are nonzero eigenvalues of $A_{y}$. 

For the second term, use the fact that \((\xi_T')_L=-\sigma_u\langle
u''_T, \xi_T\rangle_{2,t}\) (take the
\(s\)-derivative of the equation \(\langle u', \xi_T\rangle_{2,t}=0\)
to see this) and the estimates for $E_{TL}$ and \(\|u''_T\|_{2,t}\) obtained earlier. We get:
\[
\begin{split}
4\langle(A_{u(s)}\xi_T)_L, &(E_u^{(-\sigma, \sigma)}\xi_L)_L\rangle_{2,t}\\
&=-4\langle(E_u^{(-\sigma, \sigma)}\xi_T)_L,
(E_u^{(-\sigma, \sigma)}\xi_T)_L\rangle_{2,t}+4\langle(\xi'_T)_L, (E_u^{(-\sigma, \sigma)}\xi_L)_L\rangle_{2,t}\\
&\geq -C\sigma_u^{-2}\|\xi_T(s)\|_{2,t}^2.
\end{split}
\] 
Under the assumption \S5.1.2 (1), the last term can both be bounded below by
\[
\begin{split}
& -C\|u'(s)\|_{\infty,1,t}\|\xi_T(s)\|_{2,t}^2\\
&\quad \geq -C'\|u'(s)\|_{2,2,t}\|\xi_T(s)\|_{2,t}^2\\
&\quad \geq -C''s^{-2}\|\xi_T(s)\|_{2,t}^2\quad \text{for large \(s\)}.
\end{split}
\]
In summary, letting $q(s):=\|\xi_T(s)\|_{2,t}^2$, we have
\begin{equation}\label{deltaDecay}
q''(s)\geq 4\nu^2 q(s)\quad \text{for all large enough $s$. }
\end{equation} 
Now if $q'>-2\nu q$ at
some $s_0$, then $q$ must grow faster than $e^2\nu s$ for large $s$.
This is because that by comparison principle, $q$ is larger than the
solution of $q_0''=4\nu^2 q_0$ with initial conditions
$q_0^{(k)}(s_0)=q^{(k)}(s_0)$ for \(k=0,1\). 
And when $q'(0)>-2\nu q(0)$ such a solution grows 
exponentially with exponent $2\nu$ at infinity. However, such a $q$ cannot
come from an $\xi\in L^p_{1:(-\sigma, \sigma)}$, since the weight $\sigma$ is 
chosen to be much smaller than $\min |\mu_i|$, and thus smaller than $\nu$.
On the other hand, if $q'\leq -2\nu q$ for all sufficiently large $s$, by
integration one easily sees that $q\leq Ce^{-2\nu s}$; namely 
$\|\xi_T(s)\|_{2,t}\leq C'e^{-\nu s}$. Since
$\xi_T'=-A_u\xi_T-(E_u^{(-\sigma, \sigma)}\xi_T)_L$ by (\ref{eq:1}), 
we have a similar estimate for $\xi_T'(s)$, and therefore 
$\|\xi_T\|_{W_u}<\infty$. 

Next we estimate the $\xi_L$ component. Writing again 
$\xi_L(s)=f(s)\sigma_u(s)u'(s)$, it satisfies the equation:
\[
f'(s)+\sigma_u'(s)\sigma_u^{-1}f+\langle u',E_u^{(-\sigma, \sigma)}\xi_T\rangle_{2,t}\sigma_u=0.
\]
We have seen that $\|(E_u^{(-\sigma, \sigma)}\xi_T(s))_L\|_{2,t}\leq C\sigma_u(s)^{-1}\|\xi_T(s)\|_{2,t}$
decays exponentially. Thus 
\[\begin{split}
f(s)&=C\sigma_u^{-1}(s)+\mbox{exponentially decaying terms}, \\
|\sigma_u f'(s)|&\leq |\sigma_u'(s)f(s)|+\mbox{exponentially decaying
  terms}.
\end{split}\] 
\[
|\sigma_u'(s)|=\sigma_u^3(s)|\langle u'(s), u''(s)\rangle_{2,t}|
\leq \sigma_u^2(s)\|u''(s)\|_{2,t}\leq C\sigma_u(s) s^{-1}.\]
Thus $\|\xi_L\|_{W_u}\leq C\|\beta s^{-1}\|_p+\varepsilon$ is bounded.
The proof of Proposition \ref{prop:Fred} is now complete. 

\begin{*proposition}\label{y-moduli} {\rm\bf (Structure of \(\cm_P(x, y)\))}
Suppose \({\cal P}(X)\) consists of either nondegenerate or minimally
degenerate critical points, and that \(J\) satisfies \S5.1.2 (1) near
any minimally degenerate critical point. Let \(\delta\) be an
arbitrary small positive number, and either \(k>2\), \(X\) satisfies
(\ref{2ndDerv}), or \(k=2\), \(X\) arbitrary. Then there is a Baire
set \(V^{k,Preg}_\delta(J, X)\subset V^k_\delta(J, X)\), such that
for any \(H\in V^{k,Preg}_\delta(J, X)\), (RHFS2d) holds for the Floer
system generated by \((J, X+\chi_H)\) (regarded as a constant homotopy
of Floer systems).
\end{*proposition}
With the Fredholm theory established, this basically 
follows from the routine arguments in section 3. Thus, we shall only
very briefly comment on the proof. 

Since we restrict the Hamiltonian perturbation to be within
\(V^2_\delta(J, X)\), the set of critical points \({\cal
  P}(X+\chi_H)\) is independent of \(H\), so are the deformation operators
\(A_x\) for any \(x\in {\cal P}(X+\chi_H)\).
The definition of the norms \(W_u\), \(L_u\) given in \S5.2.1 involves
\(u\), and hence depends on \(H\); however, it is useful to observe
that the {\em spaces} \(W_u\) and \(L_u\) have an alternative
description which is independent of \(H\), as follows. 
\begin{notation*}\label{notation:T-pq}
Let \(p,q\in\lc\), and \(q=\exp(p,
\tilde{q})\), then \(T_{p,q}:=D_2\exp(p, \tilde{q})\) denotes the
isomorphism from \(T_p\lc\) to \(T_q\lc\), where
\(D_i\) denotes the differential with respect to the \(i\)-th variable.
\end{notation*}
The following Lemma follows from Proposition \ref{prop:decay} by 
straightforward estimates. 
\begin{lemma*}
Let \(x, y, u\) be as in Proposition \ref{prop:Fred}. Then on
\(C_0^\infty(u^*K)\), the norms \(W_u\) and \(L_u\)
are commensurate with the the norms \(W_y\), \(L_y\) below respectively:
\begin{gather*}
\|\xi\|_{W_y}:=\|\sigma_y^{1/2}\xi\|_{p,1}+\|\sigma_y\Pi_{T_{y,
    u(s)}{\bf e}_y}\xi '(s)\|_p;\\
\|\xi\|_{L_y}:=\|\sigma_y^{1/2}\xi\|_p+\|\sigma_y \Pi_{T_{y,
    u(s)}{\bf e}_y}\xi\|_p, 
\end{gather*}
where 
$\sigma_y(s)=s^{-2}$ for $s\geq1$, and \(\sigma_y=1\) otherwise. 
\end{lemma*}
As usual, \(W_y(u^*K), L_y(u^*K)\) will denote the completion of
\(C_0^\infty(u^*K)\) with respect to these norms. The necessary
modification on the definition in the case when \(x\) is minimally degenerate,
\(y\) is nondegenerate or minimally degenerate will also be implied. 
We usually prefer to work with the original defintions of \(W_u, L_u\)
in \S5.2.1, because they give simpler estimates due to the fact that
\(E_u u'=0\). 
However, the alternative definition just given is more convenient for
describing the configuration space: Let
\begin{eqnarray*}
\lefteqn{{\cal B}^W_P(x,y)=}\\
&&\Big\{ u|\, u\in L^p_{1,loc}(\Theta, p_2^*T_f), \\
&& \quad u(s,\cdot)=\exp (y, \xi_+(s, \cdot))\, \mbox{for some} \,
\xi_+\in W_y(p_2^*(y^*K))\, \mbox{when} \, s>\rho_+(u);\\
&& \quad u(s,\cdot)=\exp (x, \xi_-(s, \cdot)) \, \mbox{for some} \, \xi_-\in
L^p_{1}(p_2^*(x^*K)) \, \mbox{when} \, s<\rho_-(u)\Big\},
\end{eqnarray*}
where \(\rho_+, \rho_-, p_2\) are as in Definition \ref{B-P}. 
The usual arguments show that this is a Banach manifold, with open
neighborhoods modelling on \(T_u{\cal B}^W_P(x,y)=W_u\), and there is a Banach
bundle over \({\cal B}^W_P(x,y)\) with fibers \(L_u\). Thus,
\(\cm_P^{i,W}(x,y):=\cm_P^i\cap {\cal B}^W_P(x,y)\) is again described as the 
zero locus of a Fredholm section.
Since we showed in \S5.2.5 that \(\cok E_u=\cok E_u^{(-\sigma,
  \sigma)}\), an exponentially weighted version of the transversality
arguments in section 3 implies that \(\cm_P^{i,W}(x,y; J, X+\chi_H)\) is 
nondegenerate
for \(H\) in a Baire set \(U_1\subset V^2_\delta(J, X)\). 
On the other hand, as remarked before, the conditions \S5.1.2 (2) and
(3) hold for a Baire set, say \(U_2\), in \(V^2_\delta(J, X)\).
The decay estimate of Proposition \ref{prop:decay} implies that for
\(H\in U_2\), \(\cm_P^{i,W}(x, y; J, X+\chi_H)\) contains the whole
\(\cm^i_P(x, y; J, X+\chi_H)\). (Thus we are
justified in dropping the superscript \(W\) henceforth). 
Let \(V^{2, reg}_\delta(J, X; x, y):=U_1\cap U_2\). \(V^{2, Preg}(J,
X)\) is the intersection of such \(V^{2, reg}_\delta(J, X; x,
y)\) for all pairs of \(x\), \(y\).

The compactness of the moduli space is the consequence of Gromov
compactness plus the decay estimate of Proposition \ref{prop:decay}.

\subsection{Estimates for the New Critical Points $y_+, y_-$.}
Let \(\{(J_\lambda, X_\lambda)\}_{\lambda\in \Lambda}\) be a path in
\({\cal J}_K\times {\cal X}\) generating a CHFS, where \(\Lambda=[-1,
1]\). Let \(y\in {\cal P}^{\Lambda, deg}\), with
\(\Pi_\Lambda(y)=\lambda_0\).  We saw that if \(y\) is a (Zariski) smooth point of
\({\cal P}^\Lambda\), then it is a minimally degenerate element of
\({\cal P}(X_{\lambda_0})\).  
For explicit estimates near \(y\), it is convenient to further assume
that \(y\) is standard in the following sense. 
\begin{*definition}\label{def:std}
A minimally degenerate \(y\) in \({\cal P}^{\Lambda, deg}\) is said to be in a
{\em standard d-b neighborhood} if the following hold:
\begin{description}
\item   [(1a)] \S5.1.2 (2) holds, namely $\Pi_{{\bf e}_y}(\nabla_{{\bf
    e}_y}\nabla_{{\bf e}_y}n_y(0))=C_y {\bf e}_y$ for \(C_y>0\),
\item[(1b)] $\Pi_{{\bf e}_y}(\nabla_{{\bf e}_y}\nabla_{c}n_y(0))= 0$ for all
$c\perp {\bf e}_y$,
\item[(1c)] $\Pi_{{\bf e}_y}(\nabla_{{\bf e}_y}\nabla_{{\bf e}_y}\nabla_{{\bf e}_y}n_y(0))=0$.
\item[(2a)] \(J_{\lambda}\) is constant in \(\lambda\) in a small
neighborhood of \(\Lambda\) containg \(\lambda_0\),
and \(J_{\lambda, t}\) is integrable in a small
neighborhood of \(y(t)\in M\)  for all \(t\in S^1\). (Namely \S5.1.2 (1)).\item[(2b)] At \(\lambda=\lambda_0\), $\p_{\lambda}\check{\theta}_{X_{\lambda}}(y)=C'_y{\bf e}_y$ for a
constant \(C'_y\neq0\),
\item[(2c)] At \(\lambda=\lambda_0\), $\Pi_{{\bf e}_y}\p_{\lambda}\nabla_{v}\check{\theta}_{X_{\lambda}}(y)=0 \quad \forall v\in L^2_t(y^*K)$,
\item[(2d)] $\p_{\lambda}\check{\theta}_{X_{\lambda}}$ is supported in a neighborhood of $y$ away from all the other critical points.
\end{description}
\end{*definition}
This assumption will be particular useful in the proofs of the gluing
theorems discussed in Part II.

Assuming this, in this subsection we 
show that \(y\) bifurcates into two nondegenerate critical
points \(y_{\lambda\pm}\), and give some
essential estimates of these new critical points. 
These estimates will be useful for the gluing theorems in Part II.
\begin{*lemma}\label{ev-est}
Let \(y\in {\cal P}^{\Lambda, deg}\) be in a standard d-b
neighborhood. Then \(y\) is a local extremum of \(\Pi_\Lambda\): it is
a local maximum when the constant \(C_y'\) in \S\ref{def:std} (2b) is
positive, and a local minimum otherwise. 

Furthermore, let \(\lambda\in \Lambda\) be close to
\(\lambda_0:=\Pi_\Lambda(y)\), and \(\lambda<\lambda_0\) when
\(C_y'>0\); \(\lambda>\lambda_0\) otherwise. Let \[y_{\lambda
  +}=\exp(y, \eta_{\lambda+}),\, y_{\lambda -}=\exp (y, \eta_{\lambda
  -})\in {\cal P}_\lambda\] 
be the two points near \(y\in \lc\). 
Then they are both nondegenerate, of index $\ind_+(y)$ and
$\ind_-(y)$ respectively. Moreover, 
$$
C_{\pm}|\lambda-\lambda_0|^{-1/2}\leq \pm \langle{\mathbf{e}_y}, \eta_{\lambda \mp}\rangle_{2,t}\leq C'_{\pm}|\lambda-\lambda_0|^{-1/2},
$$
where $C_{\pm}, C'_{\pm}$ are positive constants. Furthermore, the
eigenvalue of $A_{y_{\lambda\mp}}$ with minimal absolute value is bounded above and below by multiples of $\pm
|\lambda-\lambda_0|^{1/2}$ respectively.
\end{*lemma} 
\begin{proof} 
Without loss of generality, let \(C_y'>0\).

For each \(\eta\in L^2_t(y^*K)\), write \[\eta_L:=\Pi_{\ker
  A_y}\eta=:\underline{\eta}_L{\bf e}_y, \quad \eta_T:=\eta-\eta_L.\] 
Using the assumption that \(y\) is in a
standard d-b neighborhood, 
the defining equation of critical points,
\({\cal V}_X=0\), takes the following form
in the local coordinates about \(y\):
\begin{gather}\begin{split}
&A_y\eta_T+(1-\Pi_{\ker
  A_y})n_y(\eta_L+\eta_T,\eta_L+\eta_T)+O(|\lambda-\lambda_0|^2)\\
&\quad\qquad\qquad+O\Big(|\lambda-\lambda_0|(\|\eta_L\|^2_{2,1,t}+\|\eta_T\|^2_{2,1,t})^{1/2}\Big)=0;\end{split}\nonumber\\
\begin{split}
&C_y\underline{\eta}_L^2+C_y'(\lambda-\lambda_0) +\Pi_{\ker A_y}
n_y(\eta_T,\eta_T)+O\Big(\|\eta\|_{2,1,t}(\|\eta_L\|_{2,1,t}+\|\eta_T\|_{2,1,t}^2)\Big)\\
&\quad\qquad\qquad +O(|\lambda-\lambda_0|^2)+O\Big(|\lambda-\lambda_0|(\|\eta_L\|^2_{2,1,t}+\|\eta_T\|^2_{2,1,t})\Big)=0.
\end{split}
\label{P}
\end{gather}
\begin{notation*} In this paper, \(O(\cdot)\) or \(o(\cdot)\) can be either a
number, or a function of \(t\) whose
\(L^2_{1,t}\) norm is of the order indicated.
\end{notation*}

When $|\lambda-\lambda_0|$ is small, we can assume that $(\|\eta_L\|^2_{2,1,t}+\|\eta_T\|^2_{2,1,t})^{1/2}$ is
small, and thus the first equation above implies:
\begin{equation}\label{c-est}
\|\eta_T\|_{2,1,t}\leq C_1(\|\eta_L\|^2_{2,1,t}+|\lambda-\lambda_0|^2);
\end{equation}
Substituting this into the second equation in (\ref{P}), we see that there are two
solutions for \(\eta_L\), and for both
\begin{equation}\label{b-est}
C_2|\lambda-\lambda_0|^{1/2}\leq \|\eta_L\|_{2,1,t}\leq C_2'|\lambda-\lambda_0|^{1/2}.
\end{equation}
This then implies via (\ref{c-est}) that 
$$
\|\eta_T\|_{2,1,t}\leq C_3|\lambda-\lambda_0|.
$$
The solution with positive/negative value of $\eta_L$ is denoted by
\(\eta_{\lambda -}\), \(\eta_{\lambda  +}\) respectively. 

We now estimate the small eigenvalue of
$A_{y_{\lambda \pm}}$. 
Let $\xi_{\lambda}\in L^2_{1,t}(y_{\lambda -}^*K)$, and let $(\xi_\lambda)_0\in L^2_{1,t}(y^*K)$ be
defined by \[\exp (y, \eta_{\lambda -}+(\xi_\lambda)_0)=\exp (y_{\lambda-},
\xi_{\lambda}).\] 
Using the relation between $\xi_{\lambda}$ and
$(\xi_\lambda)_0$, we see that the lowest order term of $A_{y_{\lambda -}}$ is
conjugate to the linearization of the left hand side of
(\ref{P}) at $\eta_{\lambda-}$. The latter has the form
\begin{equation}\label{A-y:lambda}
\left(\begin{array}{cc}
 A_y+O(\|(\eta_{\lambda -})_L\|_{2,1,t}) &O(\|(\eta_{\lambda -})_T\|_{2,1,t})\\
0 & 2C_y\underline{(\eta_{\lambda -})}_L+o(\|(\eta_{\lambda -})_L\|_{2,1,t})
\end{array}\right)
\end{equation}
in terms of the decomposition \(L^2_t(y^*K)=(\ker A_y)^{\perp}\oplus \ker A_y.\)
Using the estimates on $\|(\eta_{\lambda -})_L\|_{2,1,t}, \|(\eta_{\lambda -})_T\|_{2,1,t}$ above, we see that the smallest
eigenvalue is bounded above and below by positive multiples of $|\lambda-\lambda_0|^{1/2}$. Similarly for $A_{y_{\lambda
    +}}$. 
\end{proof}

\section{Existence of Admissible $(J,X)$-Homotopies.}
In this section, we show that given a path \(\{X_\lambda\}_{\lambda\in\Lambda}\) as in the statement of Theorem
\ref{thm:gen-inv}, and an arbitrary path \(\{J_\lambda\}_{\lambda\in
  \Lambda}\) connecting \(J_1, J_2\), the path \(\{J_\lambda, X_\lambda\}_{\lambda\in
  \Lambda}\) can be perturbed into an ``admissible \((J,
X)\)-homotopy'' fixing the end points. (See Definition \ref{def:adm},
Proposition \ref{thm:bifurcation} below). 
  
\subsection{Structure of Parameterized Moduli Spaces.}

In this subsection we discuss structure theorems of moduli spaces 
parameterized by an interval, parallel to the results in \S3.2, 3.3.
Since they follow from routine
modification of section 3, we shall omit most of the proofs, except
for some brief comments on the form of the relevant deformation
operators and configuration spaces, which we shall need in Part II. 

\mysubsubsection[Some terminologies.]
Let \(\Lambda\subset \R\) be an interval.
We shall use the notation \(J^\Lambda\) to denote either a
\(C_\epsilon\) complex structure on the pull-back bundle \(p_2^*K\), where \(p_2:
\Lambda\times T_f\to T_f\) is the projection, or the path of complex
structures \(\{J_\lambda\}_{\lambda\in\Lambda}\subset {\cal J}_K\)
this defines. Similarly, \(X^\Lambda\) denotes either a section in
\(C_\epsilon(p_2^*K)\), or a path
\(\{X_\lambda\}_{\lambda\in\Lambda}\). 
We denote \({\cal X}^\Lambda:=C_\epsilon(p_2^*K)\), and let \({\cal
  J}_K^\Lambda\) be the space of \(C_\epsilon\) complex structures on \(p_2^*K\).

 A {\em \((J,X)\)-homotopy} is an element \((J^\Lambda, X^\Lambda)\in
 {\cal J}^\Lambda_K\times {\cal X}^\Lambda\), or equivalently, the path \(\{(J_\lambda,
X_\lambda)\}_{\lambda\in \Lambda}\) in \({\cal J}_K\times {\cal X}\)
this defines.
\(\Lambda\) is said to be the {\em parameter} of this
\((J, X)\)-homotopy.

A \((J, X)\in {\cal J}_K\times {\cal X}\) is said to {\em belong to} a
\((J, X)\)-homotopy \(\{(J_\lambda,X_\lambda)\}_{\lambda\in \Lambda}\) if
\((J, X)=(J_\lambda, X_\lambda)\) for some \(\lambda\in \Lambda\). If
\(S\subset \Lambda\) is a sub-interval, the \((J, X)\)-homotopy
\(\{(J_\lambda, X_\lambda)\}_{\lambda\in S}\) is said to be a
{\em sub-homotopy (over \(S\))} of the \((J, X)\)-homotopy \(\{(J_\lambda, X_\lambda)\}_{\lambda\in \Lambda}\).

Let \(\Lambda\subset \R\) be an interval, and \(N\subset\Lambda\)
consists of finite points. Then
\[
{\cal H}^{\Lambda}:=C^{\infty}_{\epsilon}(\Lambda\times T_f), \quad
{\cal H}^\Lambda_N:=\{H^\Lambda\, |\, H^\Lambda\in {\cal H}^\Lambda,
\, H_\lambda=0\, \forall \lambda\in N\}.
\]
We shall often call an \(H^\Lambda\in {\cal H}^\Lambda\) a hamiltonian
isotopy, since it corresponds to a path of Hamiltonians
\(\{H_\lambda\}_{\lambda\in\Lambda}\), and hence also a path of
Hamiltonian sympletomorphisms.
Let \(\chi_{H^\Lambda}\in {\cal X}^\Lambda\) denote the path of
symplectic vector fields \(\{\chi_{H_\lambda}\}_{\lambda\in\Lambda}\).

\subsubsection{The structure of \({\cal P}^\Lambda\).}
\begin{definition*}
Let \(J^\Lambda\in {\cal J}_K^\Lambda\) be such that \(J_{\lambda, t}\) is semipositive for all
\(\lambda, t\). (Such shall be called a {\em semipositive path}). An \(X^\Lambda\in {\cal X}^\Lambda\) is
said to be {\em \(J^\Lambda\)-nondegenerate} if the following
hold:\newline
\indent (a) \({\cal P}^{\Lambda}(X^{\Lambda})\) is a (Zariski)
smooth, compact 1-manifold; \newline
\indent (b) for any \(x_\lambda\in {\cal P}^{\Lambda}(X^{\Lambda})\), \(x_\lambda(t)\not\in M_0(J_{\lambda, t})\, \forall t\);\newline
\indent (c) The projection \(\Pi_\Lambda: {\cal
  P}^\Lambda(X^\Lambda)\to \Lambda\) is smooth.
\end{definition*}
\begin{proposition*}\label{PLambda}
Fix a semipositive path \(J^\Lambda\), and an \(X^\Lambda\in
{\cal X}^\Lambda\). Then there is a
Baire set \({\cal H}^{\Lambda, ndg}(J^\Lambda, X^\Lambda)\subset{\cal
  H}^{\Lambda}\) such that \(X^\Lambda+\chi_{H^\Lambda}\) is
$J^\Lambda$-nondegenerate for all
\(H^{\Lambda}\in {\cal H}^{\Lambda, ndg}(J^\Lambda, X^\Lambda)\). 

Suppose in addition that \((J_\lambda, X_\lambda)\) are regular pairs
\(\forall \lambda\in \partial \Lambda\), then \({\cal H}^{\Lambda,
  ndg}_{\partial \Lambda}(J^\Lambda, X^\Lambda):={\cal H}^{\Lambda,
  ndg}(J^\Lambda, X^\Lambda)\cap {\cal H}^\Lambda_{\partial\Lambda}\)
is Baire in \({\cal H}^\Lambda_{\partial\Lambda}\).
\end{proposition*}

Note that \({\cal H}^{\Lambda, ndg}\) is in fact open dense, by the
compactness of \({\cal P}^\Lambda\).

The deformation operator that describes the local structure of \({\cal
  P}^\Lambda\) is an extension of \(A_x\):
Let \(x_\lambda\in {\cal P}^\Lambda\), 
\[\begin{split}
\hat{A}_{x_\lambda}: \, &
\R\oplus L^p_1(x_\lambda^*K)\to L^p(x_\lambda^*K), \,\\
\hat{A}_{x_\lambda}(\alpha,
\xi)&=A_{x_\lambda}\xi+\alpha\partial_{\lambda}\check{\theta}_{X_\lambda}(x_\lambda)+\alpha\p_\lambda
J_\lambda(\p_t x_\lambda-X_\lambda).
\end{split}\]
Note that the last term vanishes when \(x_\lambda\in {\cal P}(X_\lambda)\).

\begin{notation*}
When necessary, we shall insert \(J, X\) or \(J^\Lambda, X^\Lambda\)
as superscripts  to emphasize the dependence of the operators on them.
\end{notation*}

\subsubsection{The structure of \({\cal P}^{\Lambda, deg}\).}

\begin{proposition*}\label{degModuli}
Let \((J^\Lambda, X^\Lambda)\) be as in Proposition \ref{PLambda}. Then
there is a Baire set \[{\cal H}^{\Lambda, rg}(J^\Lambda, X^\Lambda) \subset {\cal H}^{\Lambda, ndg}(J^\Lambda, X^\Lambda)\] such that
\({\cal P}^{\Lambda}_{deg}(X^\Lambda+\chi_{H^{\Lambda}})\) consists of
finitely many (Zariski) smooth points for any \(H^{\Lambda}\in {\cal H}^{\Lambda, rg}(J^\Lambda, X^\Lambda) \). 
\end{proposition*}
\begin{proof}
Let \[
{\cal P}^{{\cal H}^\Lambda}(X^\Lambda):=\bigcup_{H^{\Lambda}\in  {\cal
    H}^{\Lambda}} {\cal P}^{\Lambda}(X^\Lambda
+\chi_{H^\Lambda})\subset \Lambda\times\lc\times {\cal H}^\Lambda,
\] 
and let \[ 
\op{pr}: {\cal P}^{{\cal H}^\Lambda}(X^\Lambda)\to {\cal H}^\Lambda,
\quad
\tilde{\Pi}_{\Lambda}: {\cal P}^{{\cal H}^\Lambda}(X^\Lambda)\to
\Lambda\] denote the projections. 

We denote an element of \({\cal P}^{{\cal H}^\Lambda}(X^\Lambda)\) by
a pair \((x_\lambda, H^\Lambda)\), where \(H^\Lambda\in {\cal
  H}^\Lambda\); \(x_\lambda\in {\cal P}^\Lambda(X^\Lambda+\chi_{H^\Lambda})\).
An element in \begin{equation}\label{nota:TP-univ}
T_{(x_\lambda, H^\Lambda)}{\cal P}^{{\cal H}^\Lambda}(X^\Lambda)\subset T_\lambda\Lambda\times
T_{x_\lambda}\lc\times T_{H^\Lambda}{\cal H}^\Lambda\end{equation} 
shall be written as a triple \((\alpha,
\xi, h^\Lambda)\), the three components belonging to the three factors
in the RHS of (\ref{nota:TP-univ}) respectively. 

The kernel of \(\mbox{pr}_*: T{\cal
  P}^{{\cal H}^\Lambda}(X^\Lambda)\to T{\cal H}^{\Lambda}\), when restricted to \[{\cal P}^{{\cal H}^{\Lambda, ndg}}(X^\Lambda)
:=\mbox{pr}^{-1}{\cal H}^{\Lambda, ndg}(J^\Lambda, X^\Lambda), \]
has constant rank. Thus it defines a real line bundle \({\cal T}\), 
\[
\Pi_{\cal T}: {\cal T}\to {\cal P}^{{\cal H}^{\Lambda , ndg}}(X^\Lambda), 
\]
so that for each \(H^\Lambda\in {\cal H}^{\Lambda, ndg}(J^\Lambda,
X^\Lambda)\), 
\[
{\cal T}\Big|_{{\cal P}^{\Lambda}(J^\Lambda, X^{\Lambda}+\chi_{H^\Lambda})}=T{\cal P}^{\Lambda}(J^\Lambda, X^\Lambda+\chi_{H^\Lambda}).\] 
Regarding \({\cal T}\) as a
subbundle of \(T{\cal P}^{{\cal H}^\Lambda}(X^\Lambda)\), the induced map
$\tilde{\Pi}_{\Lambda*}: {\cal T}\to T\Lambda$ corresponds to a
section \({\mathfrak s}\) of ${\cal T}\otimes \tilde{\Pi}^*_{\Lambda}T^*\Lambda$. 
The zero locus \[{\mathfrak s}^{-1}(0)=\bigcup_{H\in {\cal H}^{\Lambda,
    ndg}}{\cal P}^{\Lambda, deg}(X^\Lambda+\chi_{H^\Lambda})\] is the universal moduli of degenerate critical
points. If ${\mathfrak s}$ is transversal to the zero section, then
the claim of the Proposition follows from the
usual argument via the Sard-Smale theorem.

The transversality is verified as follows.
Let $(x_\lambda, H^{\Lambda})\in {\mathfrak s}^{-1}(0)\subset{\cal P}^{{\cal H}^\Lambda, ndg}(X^\Lambda)$, and let $(\beta, 0, h^{\Lambda})\in T_{(x_\lambda,
  H^{\Lambda})}{\cal P}^{{\cal H}^\Lambda}(X^\Lambda)$, where \(\beta\neq 0\). It satisfies:
\begin{gather}
\beta \partial_{\lambda}(\check{\theta}_{X_\lambda}+\nabla H_{\lambda})(x_\lambda)+\nabla h_\lambda
(x_\lambda)=0; \nonumber\\
\Pi_{{\bf e}_{x_\lambda}}\partial _{\lambda}(\check{\theta}_{X_\lambda}+\nabla H_{\lambda})(x_\lambda)\neq 0.\label{e-nablaV}
\end{gather}
To prove transversality, it suffices to show that
for some \((\beta, 0,
h^\Lambda)\) as above, 
\[(\hat{A}_{x_\lambda}^{J^\Lambda , X^\Lambda+\chi_{H^{\Lambda}}})^{-1}\left(\delta_{(\beta,
  0, h^{\Lambda})}\hat{A}_{x_\lambda}^{J^\Lambda , X^\Lambda+\chi_{H^{\Lambda}}}\right)(0, {\bf e}_{x_\lambda})\in \R\oplus L^p_1(x_\lambda^*K)\] 
has nontrivial \(\R\)-component. (The expression \(\delta_{(\beta, \eta,
  h^{\Lambda})}\hat{A}_{x_\lambda}^{J^\Lambda ,
  X^\Lambda+\chi_{H^{\Lambda}}}\) above denotes the variation of
\(\hat{A}_{x_\lambda}^{J^\Lambda , X^\Lambda+\chi_{H^{\Lambda}}}\),
regarded as an operator-valued function on \({\cal P}^{{\cal H}^\Lambda}(X^\Lambda)\)).

According to (\ref{e-nablaV}), this is equivalent to requiring
\begin{equation}\label{eq:P-rg}
\begin{split}
\Pi_{{\bf e}_{x_\lambda}}&\left(\delta_{(\beta,
  0, h^{\Lambda})}\hat{A}_{x_\lambda}^{J^\Lambda ,
  X^\Lambda+\chi_{H^{\Lambda}}}\right)(0, {\bf e}_{x_\lambda})\\
&=\Big\langle {\bf e}_{x_\lambda}, \nabla_{{\bf e}_{x_\lambda}}\left(\nabla
h_{\lambda}(x_\lambda)+\beta\partial_\lambda\check{\theta}_{X_{\lambda}}(x_\lambda)+\beta\partial _\lambda \nabla H_{\lambda
}(x_\lambda)\right)\Big\rangle_{2,t}{\bf e}_{x_\lambda}\\
&\neq 0.
\end{split}\end{equation}
Notice that ${\bf e}_{x_\lambda}$ is nowhere vanishing, being a
nontrivial solution of a linear ODE.
Thus, one may easily choose a pair of $\beta, h^{\Lambda}$ satisfying
(\ref{eq:P-rg}).
Transversality is now verified.

To verify the compactness, 
suppose the opposite, that there are infinitely many degenerate critical
points  of \({\cal P}^{\Lambda, deg}(J^\Lambda,
X^\Lambda+\chi_{H^{\Lambda}})\) for an \(H^\Lambda\in {\cal
  H}^{\Lambda, rg}(J^\Lambda, X^\Lambda)\). Since \({\cal P}^{\Lambda, deg}(J^\Lambda,
X^\Lambda+\chi_{H^{\Lambda}})\) lies in the compact 1-manifold \({\cal P}^{\Lambda}(J^\Lambda,
X^\Lambda+\chi_{H^{\Lambda}})\), it must contain a subset accumulating
at a point \(x_\lambda\in {\cal P}^{\Lambda}(J^\Lambda,
X^\Lambda+\chi_{H^{\Lambda}})\). This \(x_\lambda\) cannot be in \({\cal P}^{\Lambda}(J^\Lambda,
X^\Lambda+\chi_{H^{\Lambda}})\backslash {\cal P}^{\Lambda, deg}(J^\Lambda,
X^\Lambda+\chi_{H^{\Lambda}})\), since nondegeneracy (as a critical
point) is an open condition. Thus, \(x_\lambda\) is an accumulation
point in \({\cal P}^{\Lambda, deg}(J^\Lambda,
X^\Lambda+\chi_{H^{\Lambda}})\). On the other
hand, the assumption \(H^\Lambda\in {\cal H}^{\Lambda, rg}\) means 
that \(H^\Lambda\) is a regular value of the projection map
\(\op{pr}\Big|_{{\mathfrak s}^{-1}(0)}\), hence the linearization
 of \({\mathfrak s}\Big|_{{\cal P}^{\Lambda}(J^\Lambda,
X^\Lambda+\chi_{H^{\Lambda}})}\) at \(x_\lambda\) must be
surjective. This implies that \(x_\lambda\) is an isolated point in
\({\cal P}^{\Lambda, deg}(J^\Lambda,
X^\Lambda+\chi_{H^{\Lambda}})\), and we have arrived at a contradiction. 
\end{proof}
\begin{remarks*}
(a) Notice that the linearization of \({\mathfrak s}\Big|_{{\cal
    P}^{\Lambda}}\) corresponds to the second derivative of the
projection map \(\Pi_\Lambda: {\cal P}^\Lambda\to \Lambda\). Thus, a (Zariski)
smooth point of \({\cal P}^{\Lambda, deg}\) is a local extremum of
\(\Pi_\Lambda\). In another word, 
(RHFS1) holds for the CHFS generated by \((J^\Lambda,
X^\Lambda+\chi_{H^\Lambda})\), for any \(H^\Lambda \in {\cal
  H}^{\Lambda, rg}(J^\Lambda, X^\Lambda)\).

(b) A simple computation shows that \(x_\lambda\) is a (Zariski) smooth point of
\({\cal P}^{\Lambda, deg}\) if it is in a standard d-b neighborhood. 
\end{remarks*}

\mysubsubsection[The structure of \(\cm_P^\Lambda\), \(\cm_O^\Lambda\). ]
Let \((J^\Lambda, X^\Lambda)\) be a \((J, X)\)-homotopy such that
\(J^\Lambda\) is semipositive, and \(X^\Lambda\) is $H^1$-codirectional.
Furthermore, 
suppose (RHFS1) and (RHFS2d) hold for the CHFS associated to
\((J^\Lambda, X^\Lambda)\). Our next goal is show that (RHFS2, 3) hold for generic Hamiltonian perturbations in the
following set: Given integers \(k, \kappa\), \(k>1, \kappa\geq0\), and
a small positive real number \(\delta\), let
\begin{equation}\label{V-Lambda}
\begin{split} 
V_\delta^{\Lambda; k}(J^\Lambda,
X^\Lambda):=
\Big\{H^{\Lambda}\, \Big|\,  \|H^{\Lambda}\|_{C_{\epsilon}}\leq \delta;
\nabla ^i H_{\lambda,t} (x_{\lambda}(t))=0 \,\, &\forall x_{\lambda}\in {\cal
  P}(J_\lambda, X_{\lambda}), \\
& \qquad \lambda\in \Lambda, 0\leq i\leq k\Big\};  \\
V_\delta^{\Lambda; k, \kappa}(J^\Lambda,X^\Lambda):=
\Big\{H^{\Lambda}\, \Big|\, H^\Lambda\in V_\delta^{\Lambda; k}(J^\Lambda,
X^\Lambda),\, \partial_{\lambda}^j\nabla^i
H_{\lambda}\Big|_{\lambda=\lambda_0}=0\, \,&
\forall \lambda_0\in \Lambda_{db},\\
& \, 0\leq i\leq k, \,
0\leq j \leq\kappa\Big\};\\
V_{N,\delta}^{\Lambda; k, \kappa}(J^\Lambda,
X^\Lambda) :=
\Big\{H^{\Lambda}\, \Big| \, H^\Lambda\in V_\delta^{\Lambda; k, \kappa}(J^\Lambda,
X^\Lambda), \, H_\lambda=0\, \, & \forall\lambda\in N\Big\},
\end{split}\end{equation}
where \(N\subset\Lambda\) consists of finitely many points. 

Given a semi-positive path \(J^\Lambda\), an element \(u_\lambda\in
\cm^\Lambda_P\) or \(\cm^\Lambda_O\) is said to be
{\em \(J^\Lambda\)-regular} if it is nondegenerate and satisfies:
\[
u_\lambda(\cdot, t)\cap
M_0(J_{\lambda, t})=\emptyset\, \, \forall t.
\]
A moduli space \(\cm_P^\Lambda\) or \(\cm_O^\Lambda\) is said to be
{\em \(J^\Lambda\)-regular} if it consists of \(J^\Lambda\)-regular
elements. 
\begin{proposition*}\label{moduliLambda}
Let \((J^\Lambda, X^\Lambda)\) be a \((J, X)\)-homotopy 
such that
\(J^\Lambda\) is semipositive, and \(X^\Lambda\) is \(H^1\)-codirectional.
Suppose furthermore that
\(X^\Lambda\) is \(J^\Lambda\)-nondegenerate, and (RHFS1), (RHFS2d)
hold for the CHFS generated by \((J^\Lambda, X^\Lambda)\).
Let \(k, \kappa, \delta\) be numbers specified before (\ref{V-Lambda}).
Then there is a Baire set \[V^{\Lambda; k, \kappa, reg}_{\delta}(J^\Lambda,
X^\Lambda)\subset V^{\Lambda; k, \kappa}_{\delta}(J^\Lambda,
X^\Lambda)\] such that for any \(H^{\Lambda}\in V^{\Lambda; k, \kappa, reg}_{\delta}(J^\Lambda,
X^\Lambda)\) and any subinterval \(S\subset \Lambda\):
\newline
\indent(a) the properties (RHFS2, 3) hold for
the CHFS associated to \((J^\Lambda, X^\Lambda+\chi_{H^\Lambda})\);
\newline
\indent(b) for any integer \(i\leq 1\), the parameterized moduli spaces
\(\hat{\cm}^{S, i}_P(J^\Lambda, X^\Lambda+\chi_{H^\Lambda})\),
\(\hat{\cm}^{S, i, sim}_O(J^\Lambda, X^\Lambda+\chi_{H^\Lambda})\) are
\(J^\Lambda\)-regular.

Suppose, in addition, that there is a subset \(N\subset \Lambda\) of finitely
many elements such that \((J_\lambda,
X_\lambda)\) is a regular pair for any \(\lambda\in N\). 
Then statements (a) and (b) above hold for \(H^\Lambda\) in a Baire
subset \[
V^{\Lambda; k, \kappa, reg}_{N;\delta}(J^\Lambda,
X^\Lambda)\subset V^{\Lambda; k, \kappa}_{N;\delta}(J^\Lambda,
X^\Lambda).\]
\end{proposition*}
Notice that if \((J^\Lambda,
X^\Lambda)\) satisfies the conditions of this proposition, then 
any of its sub-homotopy satisfies the same conditions. 

A \((J,X)\)-homotopy is
said to be {\em regular} if it can be written in the form
\((J^\Lambda, X^\Lambda+\chi_{H^\Lambda})\), where \((J^\Lambda, X^\Lambda)\) satisfies the
conditions of the above proposition, and \(H^\Lambda\in V^{\Lambda; k, \kappa, reg}_{\delta}(J^\Lambda,
X^\Lambda)\).
If  we only require statements (a) and (b) in the Proposition to hold
for weight-truncated versions of moduli spaces \(\hat{\cm}^{S,
  i}_P(J^\Lambda, X^\Lambda+\chi_{H^\Lambda}; \op{wt}_{-{\cal Y},
  \op{Lf}}\leq \Re)\),
\(\hat{\cm}^{S, i}_O(J^\Lambda, X^\Lambda+\chi_{H^\Lambda}; \op{wt}_{-{\cal Y},
  \op{Lf}}\leq \Re)\), \(\forall \Re <R\), then \((J^\Lambda,
X^\Lambda+\chi_{H^\Lambda})\) is said to be a {\em \(R\)-regular}
\((J, X)\)-homotopy.

\mysubsubsection[Configuration spaces and deformation operators:
parameterized versions.]
We introduce here the configuration space \({\cal B}_P^\Lambda\) and
the deformation operator \(\hat{E}_{u_\lambda}\). These notions will
be needed in Part II. 

Let \((J^\Lambda, X^\Lambda), S\) be as in Proposition
\ref{moduliLambda}, and let \({\bf x}, {\bf y}\in
\aleph_\Lambda\). 
The topology of \(\cm_P^{S, (\sigma_1, \sigma_2)}({\bf x}, {\bf y})\) is given by its embedding into:
\[
{\cal B}_{P}^{S, (\sigma_1,\sigma_2)}({\bf x}, {\bf
  y}):=\bigcup_{\lambda\in S\cap \Lambda_{\bf x}\cap \Lambda_{\bf y}}
{\cal B}_P^{(\sigma_1, \sigma_2)}(x_\lambda, y_\lambda),
\]
which is a Banach manifold, according to the argument of
\cite{floer.cpam} Theorem 3a. 

As usual, we omit the superscript \((\sigma_1, \sigma_2)\) when it is \((0,0)\).
\begin{notation*}
An element of \({\cal B}_{P}^{S}({\bf x}, {\bf
  y})\) will be denoted either by \(u_\lambda\), as before, or as a pair
\((\lambda, u)\), where \(\lambda\in S\) and \(u\in {\cal
  B}_P(x_\lambda, y_\lambda)\). The second notation is better suited
for the discussion that follows. 
\end{notation*}
 The local model of a neighborhood of 
\((\lambda, u)\in{\cal B}_{P}^{S, (\sigma_1,\sigma_2)}({\bf x}, {\bf
  y})\) is \[
T_{(\lambda, u)}{\cal B}_{P}^{S, (\sigma_1,\sigma_2)}({\bf x}, {\bf
  y})=\R\times
L^p_{1:(\sigma_1,\sigma_2)}(u^*K),\] given via the map sending \((\mu,
\xi)\in\R\times
L^p_{1:(\sigma_1,\sigma_2)}(u^*K) \) to \(
(\lambda+\mu, e(\lambda, u; \mu,\xi))\in {\cal B}_{P}^{S, (\sigma_1,\sigma_2)}({\bf x}, {\bf
  y})\), where \(e(\lambda, u;\mu, \xi)\in {\cal
  B}_P^{(\sigma_1, \sigma_2)}(x_{\lambda+\mu}, y_{\lambda+\mu})\) is:
\begin{equation}\label{def:e}
\begin{split}
e(\lambda, u; \mu, \xi )(s)&=e_{R, R'}(\lambda, u; \mu, \xi )(s)\\
&:=\exp \Big(u(s),\xi(s)+\beta(-R-s)\bar{x}_\mu^{\lambda,u}(s)+\beta(s-R')\bar{y}_\mu^{\lambda,u}(s)\Big),
\end{split}
\end{equation}
\(R, R'\) being two large positive constants, and
\(\bar{x}_{\mu}^{\lambda,u}(s), \bar{y}_{\mu}^{\lambda, u}(s)\in
T_{u(s)}\lc\) are defined by
\[\exp(u(s), \bar{x}_{\mu}^{\lambda,u}(s))=x_{\lambda+\mu};\quad \exp(u(s), \bar{y}_{\mu}^{\lambda,u}(s))=y_{\lambda+\mu}.\]

Note that the manifold
structure of \({\cal B}_{P}^{S, (\sigma_1,\sigma_2)}({\bf x}, {\bf
  y})\) does not depend on the
choice of \(R, R'\), though the local coordinates certainly do.
For later applications, we often
  choose $R, R'$ to depend on $u,\lambda$. See e.g. II.2.2.2.

There is a Banach bundle over \({\cal B}_{P}^{S, (\sigma_1,\sigma_2)}({\bf x}, {\bf
  y})\), of which 
the parameterized moduli space \(\cm_{P}^{S, (\sigma_1,\sigma_2)}({\bf x}, {\bf
  y}; J^\Lambda, X^\Lambda)\) is the zero locus of the Fredholm section
\(\bar{\p}_{J^\Lambda , X^{\Lambda}}\),
\[
\bar{\p}_{J^\Lambda , X^{\Lambda}}(\lambda, u):=\bar{\p}_{J_\lambda , X_{\lambda}}u.
\]
The linearization of $\bar{\p}_{J^\Lambda,
X^{\Lambda}}$ with respect to the above local coordinates of ${\cal B}_{P}^{S}({\bf x}, {\bf
  y})$ has the form:
\begin{equation}\label{hat-E}
\hat{E}_{(\lambda, u)}(\alpha, \xi)=E^{J_\lambda ,
  X_{\lambda}}_u\xi+\alpha      Y_{(\lambda,u)},
\end{equation}
where $Y_{(\lambda,u)}\in L^p(u^*K)$ has the following properties:
\begin{itemize}
\item  It is supported on
$(-R-1, R'+1)\times S^1$, $R,R'$ being the real numbers in
(\ref{def:e}).
\item Over \((-R, R')\times S^1\), it agrees with \(\partial_{\lambda}\check{\theta}_{X_\lambda}(u)+\p_\lambda
J_\lambda(\p_t u-X_\lambda)\).
\item The difference \(Y_{(\lambda, u)}-\left(\partial_{\lambda}\check{\theta}_{X_\lambda}(x_\lambda)+\p_\lambda
J_\lambda(\p_t x_\lambda-X_\lambda)\right)\) consists of terms
supported on $(-R-1, R'+1)\backslash(-R, R')\times S^1$, that either involves
a product of a \(i\)-th order derivative of \(\beta(s-R')\) and a
\(j\)-th order derivative of \(\bar{y}^{\lambda, u}_\alpha\), or in parallel,
a product of a \(i\)-th order derivative of \(\beta(-R-s)\) and a
\(j\)-th order derivative of \(\bar{x}^{\lambda, u}_\alpha\), \(i, j\)
being 0 or 1.
\end{itemize}
We shall see later that for practical purposes, the difference
described in item 3 above is usually ignorable.
Finally, note that the form of \(\hat{E}_{(\lambda, u)}\) depends
on the choice of local coordinates, and hence on the numbers \(R, R'\).

\subsection{From a CHFS to an Admissible Homotopy.}
The purpose of this subsection is to use the structure theorems in
\S6.1 to establish the existence of admissible \((J,X)\)-homotopies. 
We first state the definition:
\begin{*definition}\label{def:adm}
A \((J, X)\)-homotopy \((J^\Lambda, X^\Lambda)\) is said to be {\em
  admissible}, if the following three conditions hold:
\begin{enumerate}
\item   \(J^\Lambda\) is a semi-positive path in \({\cal J}_K\), and 
  \(X^\Lambda\) is $H^1$-codirectional. Plus, for \(\lambda\in \partial \Lambda\),
  \((J_\lambda, X_\lambda)\) is regular.
\item All the properties (RHFS*) hold except for (RHFS2c, 3c, 4).
\item Any \(x_\lambda\in {\cal P}^{\Lambda, deg}(J^\Lambda,
  X^\Lambda)\) lies in a standard d-b neighborhood. 
\end{enumerate}
\end{*definition}

\begin{*proposition}\label{thm:bifurcation}
Let \(M\) be \(w^+\)-monotone, and  
\((J_1, X_1)\), \((J_2, X_2)\) be two regular pairs. Suppose
there is a $H^1$-codirectional path \(X_{(0)}^\Lambda\),
\(\Lambda=[1,2]\),
connecting \(X_1, X_2\).
Then there exists an admissible \((J, X)\)-homotopy connecting \((J_1,
X_1)\), \((J_2, X_2)\). 

\end{*proposition}
\begin{remark*}
The Morse-theoretic picture leads one to expect that a generic
\((J,X)\)-homotopy satisfying condition 1 of Definition \ref{def:adm}
generates an RHFS. In particular, condition 3 of Definition
\ref{def:adm} should be unnecessary. Since we only need the {\em
  existence}, not the genericity of RHFSs, there is no harm in
imposing this condition to simplify the estimates in Part II. 
\end{remark*}
The rest of this subsection is devoted to the proof.
\medbreak

First, we give a linear ordering of the properties of an admissible
\((J,X)\)-homotopy. We shall find a sequence of \((J, X)\)-homotopies,
\(\{(J^\Lambda_{(k)},
X^\Lambda_{(k)})\}_{k\in\{0, 1, \ldots, K\}}\), so that the \(k\)-th step \((J^\Lambda_{(k)},
X^\Lambda_{(k)})\) satisfies Properies (0)--(\(k\)).
Thus, the last of the sequence will be an admissible \((J, X)\)-homotopy.

\mysubsubsection[Ordering the properties of admissibility.]
\begin{description}\itemsep -1pt
\item[(0)] Definition \ref{def:adm}, item 1. 
(Semipositivity of
  \(J^\Lambda\) and \(H^1\)-codirection of \(X^\Lambda\)).
\item[(1)] (RHFS1). 
(Smoothness and compactness properties of \({\cal
    P}^\Lambda\), \({\cal P}^{\Lambda, deg}\)).
\item[(2)] (RHFS1i). (Injectivity of \(\Pi_\Lambda\Big|_{{\cal
      P}^{\Lambda, deg}}\)).
\item[(3)] Condition (2a) of Definition \ref{def:std} \(\forall y\in
  {\cal P}^{\Lambda, deg}\). (Integrability and
  \(\lambda\)-independence of \(J^\Lambda\) near \({\cal P}^{\Lambda,
    deg}\)).
\item[(4)] Conditions (1a)--(1c) of Definition \ref{def:std} \(\forall y\in
  {\cal P}^{\Lambda, deg}\). (Constraints on higher
  derivatives of \(n_y\) for degenerate critical points).
\item[(5)] (RHFS2d). (Smoothness and compactness properties of
  \(\cm_P(X_\lambda)\), for \(\lambda\in \Lambda_{db}\)). 
\item[(6)] Conditions (2b)--(2d) of Definition \ref{def:std} \(\forall y\in
  {\cal P}^{\Lambda, deg}\). (Constraints on \(\partial_\lambda\check{\theta}_{X_\lambda}\) and its
  derivatives for \(\lambda\in \Lambda_{db}\)). 
\item[(7)] (RHFS2, 3). (Smoothness and compactness properties of
  \(\cm_P^\Lambda\), \(\cm_O^\Lambda\)). 
\item[(8)] (RHFS2i). (Injectivity of \(\Pi_\Lambda\Big|_{\hat{\cm}_P^{\Lambda,0}}\)).
\end{description}

To achieve these properties, we allow \(J^\Lambda\) to vary among
semipositive paths, and \(X^\Lambda\) to vary within the set
\(\{X_{(0)}^\Lambda+\chi_{H^\Lambda}|\, H^\Lambda\in {\cal
  H}^\Lambda\}\). Notice that such \(X^\Lambda \) is
$H^1$-codirectional if \(X_{(0)}^\Lambda\) is.
This ordering is chosen so that \(J^\Lambda_{(k)}\), \({\cal P}^\Lambda_{(k)}\)
are fixed from Step (3) on. (So that notions such as \(\Lambda_{db}\)
is independent of \(k\) for sufficiently large \(k\), and we are
free to suppress the superscript or parenthetical reference to
\((J^\Lambda_{(k)}, X^\Lambda_{(k)})\) from the notations in this case). The
Hamiltonian perturbation \(H^\Lambda\) is fixed step by step: The low
order derivatives of \(H_\lambda\) at \(y_\lambda\) for \(y_\lambda\in
{\cal P}^{\Lambda, deg}\) are fixed in Step (4). The functions
\(H_\lambda\) for \(\lambda\in \Lambda_{db}\) are fixed in Step
(5). The low order derivatives of \(\partial _\lambda H_\lambda\) at
\(y_\lambda\) \(\forall y_\lambda\in {\cal P}^{\Lambda, deg}\) are
fixed in Step (6). Steps (7) and (8) use the remaining freedom of
\(H^\Lambda\).

\subsubsection{Achieving Properties (0)--(7).}

\paragraph{Step (0).}
By assumption, \(X^\Lambda_{(0)}\) and \((J_1, X_1), (J_2, X_2)\) already
satisfy the requirements of Property (0).  To obtain the 0-th
\((J,X)\)-homotopy, we just need to find a semipositive path
\(J_{(0)}^\Lambda\) connecting \(J_1, J_2\). 

This is easy under the
assumption that \(M\) is \(w^+\)-monotone:
in this case, 
any two regular elements in \({\cal J}_K\)  can be connected by a
semi-positive path. To see this, notice that \({\cal J}_K\) is path-connected,
and by an analog of Lemma 3.1.2 (b), any path connecting \(J_1, J_2\)
may be perturbed into a ``regular path'', namely, a path
\(\{J_\lambda\}\) such that the space \(\bigcup_\lambda S(A,
J_\lambda)/G\) is a smooth manifold of expected dimension
(\(2n+2c_1(A)-4\)). On the other hand, the analog of Lemma 3.1.2 (a)
says that regular paths are semi-positive. 

Let \(J^\Lambda_{(0)}\) be one such regular path; \((J^\Lambda_{(0)},
X^\Lambda_{(0)})\) is our 0-th \((J, X)\)-homotopy.

\paragraph{Step (1).}
Let \((J^\Lambda_{(1)}, X^\Lambda_{(1)})=(J^\Lambda_{(0)},
    X^\Lambda_{(0)}+\chi_{H_{(1)}^\Lambda})\), where \(H_{(1)}^\Lambda\in
    {\cal H}^{\Lambda, rg}(J^\Lambda_{(0)}, X^\Lambda_{(0)})\).
By Propositions 6.1.2, \ref{degModuli}, $(J^\Lambda_{(1)},
X^\Lambda_{(1)})$ satisfies Properties (0) and (1).

\paragraph{Step (2).} There might be two distinct \(x,y\in {\cal
  P}^{\Lambda}_{deg}(J_{(1)}^\Lambda, X_{(1)}^{\Lambda})\) with 
\(\Pi_{\Lambda}x=\Pi_{\Lambda}y\). However, a simple perturbation to
\(X_{(1)}^\Lambda\) may distinguish their values under the projection map
\(\Pi_\Lambda\). For example, one may vary \(X^\Lambda_{(1)}\) locally
near the image of \(x\), setting \(X_{(2)}^\Lambda\) so that
\[
X_{(2),\lambda}=
\begin{cases}
X_{(1),f(\lambda)} &\text{over \(U_x\)}\\
X_{(1), \lambda} &\text{outside \(V_x\)},
\end{cases}
\]
where \(U_x\subset V_x\) are two small regular neighborhoods of the image of
the section \(x: S^1\to T_f\), and \(f: \Lambda\to \Lambda\) is a
diffeomorphism that agrees with the identity outside a neighborhood
of \(\Pi_\Lambda(x)\), and that \(f(\Pi_\Lambda(x))\neq
\Pi_\Lambda(x)\). 

Notice that since \(x\) and \(y\) have disjoint image (being distinct
solutions of a first-order ODE), \(U_x, V_x\) may be chosen so that
the image of \(y\) lies outside of \(V_x\). On the other hand, 
since the difference between \(X_{(2),\lambda}\) and
\(X_{(1),\lambda}\) is supported on a contractible space, it is
necessarily Hamiltonian. Thus, such \(X_{(2)}^\Lambda\) is still
$H^1$-codirectional. Now, \((J_{(1)}^\Lambda, X_{(2)}^\Lambda)\)
might not no longer satisfy Property (1), but one may add a further small
Hamiltonian perturbation to regain Property (1), using 
again Propositions 6.1.2, \ref{degModuli}.
As long as the perturbation is sufficiently small, Property (2) is
preserved. 
Setting \(J_{(2)}^\Lambda=J_{(1)}^\Lambda\), we have obtained our 2nd
\((J,X)\)-homotopy, \((J_{(2)}^\Lambda, X_{(2)}^\Lambda)\), satisfying
Properties (0), (1), (2).

\paragraph{Step (3).} For each \(\lambda\in \Pi_\Lambda{\cal P}^{\Lambda,
  deg}(X_{(2)}^\Lambda)\), we choose
\(J_{(3);\lambda}\) such that:
\begin{itemize}
\item \(J_{(3);\lambda}\in {\cal J}_K^{reg}\);
\item \(J_{(3);\lambda,
  t}\Big|_{U_{y_\lambda}\cap \pi_S^{-1}\{t\}}\) is integrable \(\forall
t\in S^1\), where \(y_\lambda\in {\cal P}(X_\lambda)\) is the unique
degenerate critical point, \(U_{y_\lambda}, V_{y_\lambda}\) are as in Step (2) above,
and \(\pi_S: T_f\to S^1\) is the projection.
\item\(J_{(3); \lambda}\) is close to \(J_{(2); \lambda}\) in
  \(C_\epsilon\)-norm over \(T_f\backslash V_{y_\lambda}\).
\end{itemize}
This is possible because of Lemma 3.1.2 (c).
Moreover, arguing as in Step (0), such \(\{J_{(3); \lambda}\}_{\lambda\in \Lambda_{db}}\) may
be extended into a semi-positive (in fact, regular) path
\(J_{(3)}^\Lambda\). If one so desires, one may also take
\(J_{(3)}^\Lambda\) to \(C_\epsilon\)-approximate \(J_{(2)}^\Lambda\),
away from the image of all \(x\in {\cal P}^{\Lambda,
  deg}(X_{(2)}^\Lambda)\) in \(\Lambda\times T_f\). By a
reparametrization of \(\Lambda\), one may assume that 
\(J_{(3)}^\Lambda\) is constant in \(\lambda\) in a
small neighborhood of \(\Pi_\Lambda{\cal P}^{\Lambda,
  deg}(X_{(2)}^\Lambda)\subset \Lambda\). Noticing that
reparametrizing \(\Lambda\) preserves semipositivity, we now have a semipositive
path \(J_{(3)}^\Lambda\) satisfying Property (3).

The \((J, X)\)-homotopy \((J_{(3)}^\Lambda, X_{(2)}^\Lambda)\) now
satisfies Properties (0), (2), (3), but Property (1) may be lost. 
In this case, one again apply Propositions 6.1.2, \ref{degModuli} to
find a \(X^\Lambda_{(3)}=X^\Lambda_{(2)}+\chi_{H_{(3)}^\Lambda}\),
so that \((J_{(3)}^\Lambda, X_{(3)}^\Lambda)\) satisfies Properties
(0), (1), (3). If \(H_{(3)}^\Lambda\) is sufficiently small, Property
(2) will still be preserved. 

\paragraph{Step (4).} Property (4) may be obtained by simply setting
\(X_{(4)}^\Lambda=X_{(3)}^\Lambda+\chi_{H^\Lambda_{(4)}}\), with
\(H^\Lambda _{(4)}\in V_{\partial \Lambda,\delta}^{\Lambda; 2}(J_{(3)}^\Lambda,
X_{(3)}^\Lambda)\), varying
\(\nabla^2\check{\theta}_{X_\lambda}(y_\lambda)\),
\(\nabla^3\check{\theta}_{X_\lambda}(y_\lambda)\) for \(y_\lambda\in
{\cal P}^{\Lambda , deg}(X_{(3)}^\Lambda)\), but leaving the lower
order derivatives unchanged. Notice that the fact that \(H^\Lambda
_{(4)}\in V_{\partial \Lambda,\delta}^{\Lambda; k}\) for \(k>1\) implies that \({\cal
  P}^\Lambda(X_{(4)}^\Lambda)={\cal P}^\Lambda(X_{(3)}^\Lambda)\) in
\(\Lambda\times \lc\), and the deformation operators
\(A_{y_\lambda}\), \(\hat{A}_{y_\lambda}\) also remain the same. Thus,
Properties (0)--(3) are still valid for \((J_{(4)}^\Lambda, X_{(4)}^\Lambda):=(J_{(3)}^\Lambda, X_{(4)}^\Lambda)\).

\paragraph{Step (5).} Let \((J_{(5)}^\Lambda,
X_{(5)}^\Lambda)=(J_{(4)}^\Lambda,
X_{(4)}^\Lambda+\chi_{H_{(5)}^\Lambda})\), where \(H_{(5)}^\Lambda\)
is such that \(H_{(5); \lambda}\in V_{\delta} ^{k, Preg}(J_{(4); \lambda},
X_{(4), \lambda})\), \(\forall \lambda\in \Lambda_{db}\). By
Proposition \ref{y-moduli}, Property (5) then holds for \((J_{(5)}^\Lambda,
X_{(5)}^\Lambda)\). By taking \(H_{(5)}^\Lambda\in V_{\partial \Lambda, \delta}^{\Lambda;
  k, \kappa}(J_{(4)}^\Lambda, X_{(4)}^\Lambda)\) for \(k>2\),
Properties (0)--(4) remain valid.

\paragraph{Step (6).}  Set \((J_{(6)}^\Lambda,
X_{(6)}^\Lambda)=(J_{(5)}^\Lambda,
X_{(5)}^\Lambda+\chi_{H_{(6)}^\Lambda})\), where \(H_{(6)}^\Lambda\in
V_{\partial \Lambda, \delta}^{\Lambda;
  k, 0}(J_{(5)}^\Lambda, X_{(5)}^\Lambda)\) for \(k>2\),
choosing \(\partial _\lambda\nabla H_{(6),
  \lambda}\) near the neighborhood of all \(y_\lambda\in {\cal
  P}^{\Lambda, deg}(J_{(5)}^\Lambda,
X_{(5)}^\Lambda)\) so that Conditions (2b)--(2d) of Definition
\ref{def:std} are met. The smallness of \(H_{(6)}^\Lambda\) and the
constraints on its low order derivatives
imply that Properties (0)--(5) remain valid.

\paragraph{Step (7).} Set \((J_{(7)}^\Lambda,
X_{(7)}^\Lambda)=(J_{(6)}^\Lambda,
X_{(6)}^\Lambda+\chi_{H_{(7)}^\Lambda})\), where \(H_{(7)}^\Lambda\in
V_{\partial \Lambda, \delta}^{\Lambda;
  k, \kappa, reg}(J_{(6)}^\Lambda, X_{(6)}^\Lambda)\) for \(k>2\),
\(\kappa\geq1\). By Proposition \ref{moduliLambda}, Property (7) holds for \((J_{(7)}^\Lambda,
X_{(7)}^\Lambda)\); on the other hand, Properties (0)---(6) remain
valid by the constraints on \(k, \kappa\).

\subsubsection{The final step: injectivity of
\(\Pi_\Lambda\Big|_{\hat{\cm}_P^{\Lambda, 0}}\).}
More care is required for this last step, 
because unlike the case of Step (2), 
\(\hat{\cm}_P^{\Lambda, 0}\) consists of possibly infinitely
many elements, while we are also under the extra constraints 
to preserve the many 
Properties already established above. 
  
 To obtain Property (8), we shall again set \((J_{(8)}^\Lambda,
X_{(8)}^\Lambda)=(J_{(7)}^\Lambda,
X_{(7)}^\Lambda+\chi_{H_{(8)}^\Lambda})\), where \(H_{(8)}^\Lambda\in
V_{\partial \Lambda, \delta}^{\Lambda;
  k, \kappa, reg}(J_{(7)}^\Lambda, X_{(7)}^\Lambda)\) for \(k>2\),
\(\kappa\geq1\). The Hamiltonian perturbation \(H_{(8)}^\Lambda\) will
be chosen iteratively below, via a refinement of the standard
transversality argument.

First, observe the following
\begin{lemma*}\label{pert:handle2}
Suppose Properties (0)--(7) above hold for the \((J,X)\)-homotopy
\((J^\Lambda, X^{\Lambda})\), and let 
\[\{(\lambda_1, u_1), (\lambda_2,
u_2,), \ldots (\lambda_n, u_n)\}\subset
\hat{\cm}^{\Lambda,0}_P(J^\Lambda, X^{\Lambda})\] 
be such that $\lambda_1, \lambda_2, \ldots,
\lambda_{n-1}\in \Lambda$ are mutually distinct, and $\lambda_n$ agrees with one
$\lambda_i$, for some $i\in \{1, \ldots, n-1\}$. Let \(k, \kappa\) be
as in Step (7). Then for any
sufficiently small \(\delta\), there exists 
an $H^{\Lambda}\in
V^{\Lambda,k, \kappa, reg}_{\partial \Lambda, \delta}(J^\Lambda, X^\Lambda)$, such that:

Denoting by \((\tilde{\lambda}_j, \tilde{u}_j)\) the unique element in \(\hat{\cm}^{\Lambda,0}_P(J^\Lambda, X^{\Lambda}+\chi_{H^\Lambda})\) close to \((\lambda_j, u_j)\) in \({\cal
  B}_P^\Lambda\), the points \(\tilde{\lambda}_1, \tilde{\lambda}_2, \ldots , \tilde{\lambda}_{n-1},
\tilde{\lambda}_n\) are mutually distinct in \(\Lambda\). 
\end{lemma*}
\begin{proof}
Ideally, we would like the Hamiltonian perturbation to: \newline
\indent (i) shift \(\lambda_n\), but meanwhile\newline
\indent (ii) leave \(\lambda_1, \ldots, \lambda_{n-1}\) unchanged. 
\medbreak

To guarantee (ii), one may try an \(h^\Lambda\in V^{\Lambda,k,
  \kappa}_{\partial \Lambda, \delta}(J^\Lambda, X^\Lambda)\), such
that: 
\begin{gather}\label{h:a}\qquad
\hbox to 10pt{\hss\hbox to 
  \hsize{\noindent\hphantom{$\bullet$\enspace}\parbox{0.8\hsize}{\hspace{0pt}\llap{$\bullet$\enspace}$h_{\lambda}$ vanishes when \(\lambda\in \Lambda\backslash
  S_{\lambda_n}\), where \(S_{\lambda_n}\subset \Lambda\) is a small
  interval about \(\lambda_n\), so that 
\(
S_{\lambda_n}\cap (\Lambda_{db}\cup \{\lambda_1, \ldots, \lambda_{n-1}\})=\emptyset;
\)
}}\hss}\\\label{h:b}\qquad
\hbox to 10pt{\hss\hbox to \hsize{$\bullet$\enspace\parbox{0.8\hsize}
{\(h_\lambda\) is supported away from the image of \(u_i\) in
\(T_f\), $\forall i\in \{1, \ldots, n-1\}$.}}\hss}
\end{gather}
To explain how to achieve (i), some preliminary discussion is
required.

Let \({\cal B}_P^\Lambda=\coprod_{{\bf x}, {\bf
    y}\in\aleph_\Lambda}{\cal B}_P^\Lambda({\bf x}, {\bf y})\), and
\[
\cm_P^{1, V_{\partial\Lambda, \delta}^{\Lambda; k, \kappa}}(J^\Lambda,
  X^\Lambda):=\bigcup_{H^\Lambda\in V_{\partial\Lambda, \delta}^{\Lambda; k, \kappa}(J^\Lambda,
  X^\Lambda)}\cm_P^{\Lambda,1}(J^\Lambda,
X^\Lambda+\chi_{H^\Lambda})\subset {\cal B}_P^\Lambda
\]
be the universal moduli space. Let \((\lambda, u)\in \cm_P^{\Lambda,1}(J^\Lambda,
X^\Lambda )\), and let \(((\lambda, u), 0)\) denote the associated
element in \(\cm_P^{1, V_{\partial\Lambda, \delta}^{\Lambda; k, \kappa}}(J^\Lambda,
  X^\Lambda)\).
An element in its tangent space \[(\alpha,\xi,  h^{\Lambda})\in\R\times
L^p_1(u^*K)\times V_{\partial\Lambda, \delta}^{\Lambda; k, \kappa}(J^\Lambda,
  X^\Lambda)=
T_{((\lambda, u), 0)}\cm_P^{1, V_{\partial\Lambda, \delta}^{\Lambda; k, \kappa}}(J^\Lambda,
  X^\Lambda)\] satisfies 
\begin{equation}\label{Tmod}
E_u^{J_\lambda, X_\lambda}(\xi)+\alpha Y_{(\lambda,u)}+\nabla
h_{\lambda}(u)=0,
\end{equation}
\(Y_{(\lambda,u)}\) being as in (\ref{hat-E}).

Let \({\frak f}_u\) be a unit vector in \(\cok E_u\).
Since \((J^\Lambda, X^\Lambda)\) is assumed to satisfy
(RHFS2), \((\lambda_j, u_j)\) is an nondegenerate point of \( \cm_P^{\Lambda,1}(J^\Lambda,
X^\Lambda )\) for any \(j\); hence \(\cok
E_{u_j}=\op{Span}\{\mathfrak{f}_{u_j}\}\) and 
\[\Pi_{{\frak f}_{u_j}}Y_{(\lambda, u)}\ne 0.\]

Thus, by (\ref{Tmod}), a small $h^{\Lambda}$ would satisfy (i) if 
\begin{equation}\label{h:c}
\Pi_{{\frak f}_u}\nabla h_{\lambda}(u_n)\neq 0.
\end{equation}

We now show the existence of \(h^\Lambda \) satisfying all the three
conditions (\ref{h:a}), (\ref{h:b}), (\ref{h:c}).

As in the standard transversality argument (as in the proof of
Theorem 5.1 (i) of \cite{trans}), consider the following three
possibilities for \({\frak f}_{u_n}\): 
\underline{\sc Case 1}: The set
\[
Q_1:=\Big\{ (s,t)\, \Big|\, (s,t)\in \R\times S^1_1, \,{\frak f}_{u_n}(s,t)\neq
Cu'_{n}(s,t)\, \forall C\in \R\Big\}
\]
is nonempty. 
\underline{\sc Case 2}: ${\frak f}_{u_n}=g(s,t)u'_n$, and the set
\[
Q_2:=\{(s,t)\, |\, \p_s g(s,t)\neq 0\}
\]
is nonempty. 
\underline{\sc Case 3}. ${\frak f}_{u_n}=g(t)u'_n$. 
\medbreak

Case 3 was shown
to be impossible in the proof of Theorem 5.1 (i) of \cite{trans}. 

For the other two cases, apply unique continuation (cf. \cite{trans}) as in
section 3 to see that the image of $u_i, u_n$ at most intersect at 
discrete points. On the other hand, 
both sets \(Q_1\) and \(Q_2\) are open. Thus, in both cases 1 and 2, we can
find a neighborhood \(B\) in \(Q_1\) or \(Q_2\),
such that $u_n(B)$ does not intersect the image of $u_i$.

In either case, we can choose a small 
$\underline{H}_{\lambda_n}\in V^{k}_\delta(J_{\lambda_n}, X_{\lambda_n})$
supported in a small neighborhood ${\mathfrak B}\subset T_f$, so that
\(u_n^{-1}({\mathfrak B})\subset B\), and $\langle {\frak f}_{u_n}, \nabla \underline{H}_{\lambda_n}
(u_n)\rangle_2\neq 0$, as in the proof of Theorem 5.1 (i) of \cite{trans}.
This $\underline{H}_{\lambda_n}$ can be extended smoothly to get a small
$\underline{H}^{\Lambda}\in V_{\partial\Lambda, \delta}^{\Lambda; k, \kappa}(J^\Lambda,
  X^\Lambda)$ satisfying (\ref{h:a}).

\(\underline{H}^\Lambda\) now satisfies (\ref{h:a}), (\ref{h:b}),
(\ref{h:c}) by construction. 
Since
\(V^{\Lambda,k, \kappa, reg}_{\partial \Lambda, \delta}(J^\Lambda,
X^\Lambda)\) is dense in \(V_{\partial\Lambda, \delta}^{\Lambda; k, \kappa}(J^\Lambda,
  X^\Lambda)\), we
can approximate \(\underline{H}^{\Lambda}\) by an element
\(H^\Lambda\in V_{\partial\Lambda, \delta}^{\Lambda; k, \kappa, reg}(J^\Lambda,
  X^\Lambda)\), still keeping \(\Pi_\Lambda\) injective on
\(\{(\tilde{\lambda}_1, \tilde{u}_1), \ldots, (\tilde{\lambda}_{n-1},
\tilde{u}_{n-1}), (\tilde{\lambda}_n, \tilde{u}_n)\}\subset \hat{\cm}_P^{\Lambda, 0}(J^\Lambda,
X^\Lambda+\chi_{H^\Lambda})\).
\end{proof}
\medbreak

We now return to Step (8) of the proof of Proposition \ref{thm:bifurcation}.

Since \((J^\Lambda_{(7)},X^{\Lambda}_{(7)})\) satisfies (RHFS2), 
there are countably many elements in
$\hat{\cm}^{\Lambda,0}_P(J^\Lambda_{(7)},X^{\Lambda}_{(7)})$; let's enumerate
them as $(\lambda_1, u_1), (\lambda_2, u_2), \ldots(\lambda_n, u_n), \ldots$, such that the
weight \(\op{wt}_{-\langle {\cal Y}\rangle, e_{\cal P}}\) increases monotonically in
\(n\). An small Hamiltonian perturbation to $X^\Lambda_{(7)}$ 
will only change them slightly.

We now apply the previous Lemma recursively to the first \(n\)
elements of \(\cm_P^{\Lambda, 0}(J^\Lambda, X^\Lambda)\)
(in the order of \(\op{wt}_{-\langle {\cal Y}\rangle, e_{\cal P}}\)) for each \(n\), labeling
the $(J^\Lambda, X^\Lambda)$ used in the $n$-th step by
$(\tilde{J}^\Lambda_{[n]}, \tilde{X}^\Lambda_{[n]})$, the number
$\delta$ used in the $n$-th step by $\tilde{\delta}_{[n]}$, and the small
Hamiltonian pertubation $H^\Lambda$ obtained in the $n$-th step by $\tilde{H}^\Lambda_{[n]}$.
Set \begin{gather*}
\tilde{J}^\Lambda_{[n]}=J^\Lambda_{(7)}\quad \forall n\in \Z^+,\\
\tilde{X}^\Lambda_{[1]}:=X^\Lambda_{(7)}, \quad 
\tilde{X}^\Lambda_{[n]}=\tilde{X}^\Lambda_{[n-1]}+\chi_{\tilde{H}^\Lambda_{[n-1]}}\,
\, \text{for \(n>1\). }
\end{gather*}
Notice that \[V_{\partial\Lambda, \delta}^{\Lambda; k, \kappa}(J^\Lambda,
  X^\Lambda)=V_{\partial\Lambda, \delta}^{\Lambda; k, \kappa}(J^\Lambda,
  X^\Lambda+\chi_{H^\Lambda})\quad \text{if \(H^\Lambda\in V_{\partial\Lambda, \delta}^{\Lambda; k, \kappa}(J^\Lambda,
  X^\Lambda)\)}.\]
Thus, all the various \(\tilde{H}^\Lambda_{[n-1]}\), or any sum of
them, are actually in the same space, namely \(V_{\partial\Lambda, \delta}^{\Lambda; k, \kappa}(J^\Lambda_{(7)},
  X^\Lambda_{(7)})\) for certain \(\delta>0\).

To explain our choice of \(\tilde{\delta}_{[n]}\), we need to first
give an explicit description of the Baire set
\(V_{\partial\Lambda, \delta}^{\Lambda; k, \kappa, reg}(J^\Lambda,
  X^\Lambda)\) as a countable intersection of certain open dense sets. 

Let \(V_{\partial\Lambda, \delta}^{\Lambda; k, \kappa, R-reg}(J^\Lambda,
  X^\Lambda)\subset V_{\partial\Lambda, \delta}^{\Lambda; k, \kappa}(J^\Lambda,
  X^\Lambda)\) be the subset consisting of all \(H^\Lambda\) such that
\((J^\Lambda, X^\Lambda+\chi_{H^\Lambda})\) is \(R\)-regular. By the
compactness of truncated moduli spaces, 
this is an open dense set. By construction, 
\[
V_{\partial\Lambda, \delta}^{\Lambda; k, \kappa, reg}(J^\Lambda,
  X^\Lambda)=\bigcap_{R>0} V_{\partial\Lambda, \delta}^{\Lambda; k, \kappa, R-reg}(J^\Lambda,
  X^\Lambda).
\]
When \((J^\Lambda, X^\Lambda)\) is regular, given \(R\in \R^+\), there is an
\(\varepsilon_R(J^\Lambda, X^\Lambda)\in \R^+\) such that 
\[
V_{\partial\Lambda, \delta}^{\Lambda; k, \kappa, R-reg}(J^\Lambda,
  X^\Lambda)=V_{\partial\Lambda, \delta}^{\Lambda; k, \kappa}(J^\Lambda,
  X^\Lambda)\quad \forall \delta\leq \varepsilon_R(J^\Lambda, X^\Lambda).
\]

Now we are ready to state:
\begin{equation}\label{def:delta-n}
\tilde{\delta}_{[n]}:=\op{min}\, 
\Big\{2^{-n+j-1}\varepsilon_{R_j}(\tilde{J}^\Lambda_{[j]},
\tilde{X}^\Lambda_{[j]}) \, \Big| \, 1\leq j\leq n\Big\},
\end{equation}
where \(R_j=\op{wt}_{-\langle {\cal Y}\rangle, e_{\cal
    P}}(u_j)\). Since the weight is defined by cohomological pairing,
this is independent of hamiltonian
perturbations. Namely, \(R_j=\op{wt}_{-\langle {\cal Y}\rangle,
  e_{\cal P}}(\tilde{u}_j)\) also, in the notation of the above lemma. 

Finally, set
\begin{gather*}
J^\Lambda_{(8)}=J^\Lambda_{(7)};\\
X^\Lambda_{(8)}=X^\Lambda_{(7)}+\sum_{n=1}^\infty\chi_{\tilde{H}_{[n]}^\Lambda}
=\lim_{n\to\infty}\tilde{X}^\Lambda_{[n]}.
\end{gather*}

Notice that the limit exists by (\ref{def:delta-n}). The limit of a
sequence of regular \((J,X)\)-homotopies might not be regular;
however, (\ref{def:delta-n}) also implies:
\[
\Big\|\sum_{j=n}^\infty \tilde{H}^\Lambda_{[j]}\Big\|_{C_\epsilon}\leq \varepsilon_{R_n}(\tilde{J}^\Lambda_{[n]},
\tilde{X}^\Lambda_{[n]}).
\]
Thus, 
\[
X_{(8)}^\Lambda:=\tilde{X}_{[n]}^\Lambda+\chi_{\sum_{j=n}^\infty\tilde{H}_{[n]}^\Lambda}
\]
is \(R_n\)-regular, for any \(n\in \Z^+\).

Depending on whether
$\hat{\cm}^{\Lambda,0}_P(J^\Lambda_{(7)},X^{\Lambda}_{(7)})$ consists
of finitely many points, either \(n\) has a maximum, or \(R_n\to
\infty\). 
Either way, \((J^\Lambda_{(8)}, X_{(8)}^\Lambda)\) will be a regular
  $(J, X)$-homotopy: in the first case,
  \((J^\Lambda_{(8)}, X_{(8)}^\Lambda)=(\tilde{J}^\Lambda_{[n]},
  \tilde{X}_{[n]}^\Lambda)\) for certain \(n\), and is regular by
  definition; in the second case, the above argument shows that it is
  \(R\)-regular for all \(R>0\).

\((J^\Lambda_{(8)}, X_{(8)}^\Lambda)\) is the admissible
$(J,X)$-homotopy we desire.
End of the proof of Proposition \ref{thm:bifurcation}.

{\small
\paragraph{Acknowledgments.} This work began in 1997 as a joint project with
Michael Hutchings. The author is indebted
to M. H. for his early input during 1997--1998; 
especially, the finite-dimensional picture he
worked out in \cite{Hut} was indispensable. She also thanks
P. Seidel for suggesting this problem, M. Sullivan for explaining his
related work, D. Burghelea
for communicating their relevant preprint when this work was
first circulated in 2001. 
The author was supported by MSRI, IAS, and NSF grant DMS \# 0333163
during the multiple re-writings of this work.

}

\begin{thebibliography}{FuOOO}\itemsep -1.5pt
\bibitem[BCOV]{BCOV} M. Bershadsky, S. Cecotti, H. Ooguri, C. Vafa, {\it
    Kodaira-Spencer theory of gravity and exact results for quantum
    string amplitudes}, Commun. Math. Phys. 165, 311--427 (1994). 
\bibitem[Ba]{banyaga} A. Banyaga, {\it Sur la structure du groupe des
    difféomorphismes qui préservent une forme symplectique}, 
Comment. Math. Helv. 53 (1978), no. 2, 174--227.
\bibitem[BH]{BH} D. Burghelea, S. Haller, {\em Non-contractible
    periodic trajectories of symplectic vector fields, Floer
    cohomology and symplectic torsion}. math.SG/0104013 
\bibitem[CFH]{CFHW} K. Cieliebak; A. Floer; H. Hofer, {\it Symplectic homology. II. A general constructio}n. Math. Z. 218 (1995), no. 1, 103--122.
\bibitem[CJS]{CJS}R. Cohen, J. Jones, G. Segal, {\it Morse theory and
    classifying space}, preprint available from {\tt http://math.stanford.edu/\~{}ralph/papers.html}. 
\bibitem[DoS]{dostoglou-salamon} S. Dostoglou and D. Salamon, {\it
    Self-dual instantons and holomorphic curves},  Ann. of
Math. (2) 139 (1994), no. 3, 581--640.

\bibitem[EES]{EES}T. Ekholm, J. Etnyre, M. Sullivan, {\em Submanifolds
    in $R^{2n+1}$ and Contact Homology}, math.SG/0210124.
\bibitem[El]{elias}Y. Eliashberg, {\it Topological characterization of Stein manifolds of dimension $>2$.} Internat. J. Math. 1 (1990), no. 1, 29--46. 
\bibitem[ElG]{Gromov-Eliashberg} Y. Eliashberg, M. Gromov, {\it Lagrangian intersection theory: finite-dimensional approach.} 
 in Geometry of differential equations, AMS, 1998. 

\bibitem[F88a]{floer.cpam}A. Floer, {\it The unregularized gradient flow of the symplectic action.} Comm. Pure Appl. Math. 41
(1988), no. 6, 775--813.


\bibitem[F88b]{floer.jdg} A. Floer, {\it Morse theory for Lagrangian intersections}, J. Differential Geom. 28 (1988), no. 3,
513--547. 


\bibitem[F89]{floer.cmp} A. Floer, {\it Symplectic fixed points and holomorphic spheres}. Comm. Math. Phys. 120 (1989),
no. 4, 575--611.


\bibitem[FH]{FH} A. Floer and H. Hofer, {\it Coherent orientations for periodic orbit problems in symplectic geometry}, Math.
Z. 212 (1993), no. 1, 13--38.


\bibitem[FHS]{trans} A. Floer, H. Hofer, D. Salamon, {\it Transversality in elliptic
Morse theory for the symplectic action}, Duke Math. J 80, 251, 1995.

\bibitem[Fu]{fukaya} K. Fukaya, {\it The symplectic s-cobordism conjecture: a
summary}, in Geometry and physics (Aarhus, 1995), Dekker.
\bibitem[FuOh]{fukaya-oh} K. Fukaya, Y. Oh, {\it Zero-loop open strings in
    cotangent bundles and Morse homotopy}, Asian J. Math. {\bf 1}
  (1998), 96.
\bibitem[FuOno]{fukaya-ono} K. Fukaya and K. Ono, {\em Arnold conjecture and
Gromov-Witten invariant}, Topology 38 (1999), no. 5, 933--1048.
\bibitem[FuOOO]{FOOO} Fukaya, Oh, Ohta, Ono, {\it Lagrangian
    intersection Floer theory\allowbreak---Anomaly and obstruction}, preprint available from\\
{\tt http://www.math.kyoto-u.ac.jp/\~{}fukaya/fukaya.html}.
\bibitem[G]{gromov}M. Gromov, {\it Pseudo-holomorphic curves in
    symplectic manifolds}, Invent. Math. 82 (1985), no. 2, 307--347.

\bibitem[HS]{HS} H. Hofer and D. Salamon, {\it Floer homology and Novikov
    rings}, in The Floer memorial volume, Progr. Math., 133, Birkhäuser, Basel, 1995.

\bibitem[HuL1]{hl1} M. Hutchings and Y-J. Lee, {\it Circle-valued Morse theory, Reidemeister torsion, and
Seiberg-Witten invariants of $3$-manifolds}, Topology 38 (1999), no. 4, 861--888. 

\bibitem[HuL2]{hl2} M. Hutchings and Y-J. Lee, {\it Circle-valued Morse theory and Reidemeister torsion}, Geom.
Topol. 3 (1999), 369--396.

\bibitem[Hu]{Hut} M. Hutchings, {\it Reidemeister torsion in generalized Morse theory}. Forum Math. 14 (2002), no. 2, 209--244.


\bibitem[Ig]{igusa} K. Igusa, {\it Higher Franz-Reidemeister torsion},
 AMS/IP Studies in Advanced Mathematics, 31, 2002. 

\bibitem[IP]{IP} E. Ionel, T. Parker, {\it Gromov invariants and symplectic maps}. Math. Ann. 314 (1999), no. 1, 127--158.
\bibitem[K]{kont-icm}M. Kontsevitch, {\it Homological algebra in mirror
    symmetry}, in Proceedings of 1994 ICM.
\bibitem[Pt2]{part2}Y-J. Lee, {\it Reidemeister torsion in symplectic
    Floer theory and counting pseudo-holomorphic tori II}, preprint
  available from\\ {\tt http://www.math.purdue.edu/\~{}yjlee}.
\bibitem[L]{L} Y-J. Lee, {\it Seiberg-Witten theory on 3-manifolds with
  euclidean ends}, Commun. Analysis and Geometry, to appear. 
\bibitem[LS]{LS} Y-J. Lee, M. Sullivan, {\it Reidemeister in the Floer of
  lagrangian intersections}, unfinished manuscript.
\bibitem[L2]{L2} Y-J. Lee, {\it Noncontractible periodic orbits,
    Gromov invariants, and Floer-theoretic torsions}, math.SG/0308185.
\bibitem[L3]{ogw} Y-J. Lee, 
{\it Corrected open Gromov invariants}, in preparation. 
\bibitem[LeOno]{le-ono} Le Hong Van and K. Ono, {\it Symplectic fixed points, the Calabi invariant and Novikov homology},
Topology 34 (1995), no. 1, 155--176. 

\bibitem[LiuT]{liu-tian} G. Liu and G. Tian, {\it Floer homology and Arnold conjecture}, J. Differential Geom. 49 (1998), no.
1, 1--74.

\bibitem[Mi]{Mil}J. Milnor, {\em Whitehead torsion},
  Bull. Amer. Math. Soc. 72 (1966) 358-426. 
\bibitem[McS]{MS} D. McDuff, D. Salamon, J-holomorphic curves and quantum
  cohomology, University Lecture Series, 6. AMS 1994.

\bibitem[MMR]{MMR} J. Morgan, T. Mrowka, D. Ruberman, The \(L^2\)-moduli
  space and a vanishing theorem for Donaldson polynomial invariants,
  IP 1994.

\bibitem[N]{novikov} S. Novikov, {\it Multivalued functions and functionals. An analogue of the Morse theory.} (Russian) Dokl.
Akad. Nauk SSSR 260 (1981), no. 1, 31--35.






\bibitem[R]{ruan} Y. Ruan, {\it Virtual neighborhoods and pseudo-holomorphic curves}. Proceedings of 6th
G\:{o}kova Geometry-Topology Conference. Turkish J. Math. 23 (1999), no. 1, 161--231.



\bibitem[S]{salamon.park} D. Salamon, Lectures on Floer homology, in
  Symplectic geometry and topology (Park City, UT, 1997), AMS 1999.

\bibitem[SZ]{salamon-zehnder} D. Salamon and E. Zehnder, {\it Morse Theory for Periodic
Solutions of Hamiltonian Systems and the Maslov
index}, Comm. Pure Appl. Math, 17, 1303 (1992).
\bibitem[Sch]{schwarz} M. Schwarz, Morse homology, Birkhauser 1993.
\bibitem[Sei]{Sei} P. Seidel, {\it Symplectic Floer homology and the mapping
  class group},  Pacific J. Math. 206 (2002), no. 1, 219--229.


\bibitem[Sul]{Sul}M. Sullivan, {\it $K$-theoretic invariants for Floer
  homology}. Geom. Funct. Anal. 12 (2002), no. 4, 810--872.
.
\bibitem[T]{T2} C.H. Taubes, {\it Counting Pseudo-holomorphic submanifolds in
  dimension 4}, J. Diff. Geom., 44 (1996), no. 4, 818--893.
\bibitem[Tu]{turaev:spinc}V. Turaev, {\it Torsion invariants of ${\rm Spin}\sp c$-structures on $3$-manifolds}. Math. Res.
Lett. 4 (1997), no. 5, 679--695.
\bibitem[V]{Vit}C. Viterbo, {\it Functors and computations in Floer homology with applications. I.} Geom. Funct. Anal. 9 (1999), no. 5, 985--1033. 
\end{thebibliography}
\end{document}